\documentclass[twoside,11pt]{article}

\usepackage{jmlr2e}

\usepackage{mathrsfs}
\usepackage{amsfonts}
\usepackage{amssymb,amsmath}
\usepackage{bm}

\usepackage{graphics}
\usepackage{graphicx}
\usepackage{epsfig}
\usepackage{algorithm} 
\usepackage{algorithmic} 
\usepackage{multirow} 
\usepackage{xcolor}
\usepackage{makecell}
\usepackage{setspace}
\usepackage{tikz}
\usetikzlibrary{backgrounds}
\numberwithin{equation}{section}
\newtheorem{remark}[theorem]{Remark}
\newtheorem{assumption}[theorem]{Assumption}

\usepackage{changes}

\DeclareMathOperator*{\argmin}{\textup{arg\,min}}
\ShortHeadings{Variational Inverting Network}{Jia, Wu, Li and Meng}
\firstpageno{1}

\begin{document}

\title{Variational Inverting Network for Statistical Inverse Problems of Partial Differential Equations}

\author{\name Junxiong Jia \email jjx323@mail.xjtu.edu.cn \\
       \addr School of Mathematics and Statistics, \\
       Xi'an Jiaotong University,\\
       Xi'an, 710049, China
       \AND
       \name Yanni Wu \email wuyanni@stu.xjtu.edu.cn \\
       \addr School of Mathematics and Statistics, \\
       Xi'an Jiaotong University,\\
       Xi'an, 710049, China
   	   \AND
   	   \name Peijun Li \email lipeijun@math.purdue.edu \\
   	   \addr Department of Mathematics, \\
   	   Purdue University, \\
   	   West Lafayette, Indiana, 47907, USA
   	   \AND
       \name Deyu Meng\footnote{Corresponding author} \email dymeng@mail.xjtu.edu.cn \\
       \addr School of Mathematics and Statistics,\\
       Xi'an Jiaotong University,\\
       Xi'an, 710049, China; \\
   	   Pazhou Lab, Guangzhou, 510330, China}

\editor{}

\maketitle

\begin{abstract}
To quantify uncertainties in inverse problems of partial differential equations (PDEs), we formulate them into statistical inference problems using Bayes' formula. Recently, well-justified infinite-dimensional Bayesian analysis methods have been developed to construct dimension-independent algorithms. However, there are three challenges for these infinite-dimensional Bayesian methods: prior measures usually act as regularizers and are not able to incorporate prior information efficiently; complex noises, such as more practical non-i.i.d. distributed noises, are rarely considered; and time-consuming forward PDE solvers are needed to estimate posterior statistical quantities. To address these issues, an infinite-dimensional inference framework has been proposed based on the infinite-dimensional variational inference method and deep generative models. Specifically, by introducing some measure equivalence assumptions, we derive the evidence lower bound in the infinite-dimensional setting and provide possible parametric strategies that yield a general inference framework called the Variational Inverting Network (VINet). This inference framework can encode prior and noise information from learning examples. In addition, relying on the power of deep neural networks, the posterior mean and variance can be efficiently and explicitly generated in the inference stage. In numerical experiments, we design specific network structures that yield a computable VINet from the general inference framework. Numerical examples of linear inverse problems of an elliptic equation and the Helmholtz equation are presented to illustrate the effectiveness of the proposed inference framework.
\end{abstract}

\begin{keywords}
Infinite-dimensional variational inference, inverse problems, Bayesian analysis for functions, partial differential equations, deep neural networks
\end{keywords}

\section{Introduction}\label{introduction}

Motivated by their significant applications in seismic exploration, radar imaging, and many other domains, inverse problems of partial differential equations (PDEs) have undergone enormous development over the past few decades \citep{Engl1996book}. As computational power keeps increasing, researchers are not satisfied with obtaining just an estimated solution but pursue performing some statistical analysis based on uncertain information, which is essential for some applications, such as artifact detection \citep{Zhou2020SIIMS}. The Bayesian inverse approach transforms inverse problems into statistical inference problems and has provided a framework for analyzing the uncertainties of the interested parameters in inverse problems of PDEs \citep{Stuart2010AN}.

As is known, PDEs are usually defined on some infinite-dimensional spaces of functions \citep{Evans2010PDEbook}, which makes it difficult to directly use the finite-dimensional Bayes’ formula. To resolve this issue, a straightforward approach is to discretize PDEs in order to approximate the original problem in some finite-dimensional spaces and then solve the approximate problems by applying finite-dimensional Bayes’ methods. This strategy makes nearly all of the parametric Bayesian inference methods developed in the statistical literature applicable \citep{Berger1980book,Kaipio2004Book}. Recently, under the finite-dimensional setting, machine learning methods, e.g., Wasserstein GAN \citep{Adler2018arXiv}, have been employed to learn the prior and construct fast sampling algorithms. However, given that the original problems are defined in infinite-dimensional spaces, two critical issues are inevitably encountered. 
\begin{enumerate}
	\item Model consistency: The elements in a finite-dimensional model (e.g., prior probability measure) possess different intrinsic characteristics from their intuitive infinite-dimensional counterparts. One typical example is the total variation prior measure, which has been comprehensively analyzed in \citep{Lassas2004IP}.
	\item Algorithm applicability: To ensure that the algorithms can preserve the structures of the original infinite-dimensional problems, the algorithms designed on finite-dimensional space need to be reconstructed on infinite-dimensional space. Through re-design, the algorithms well-defined on infinite-dimensional space will have consistent behavior for different discretizations; see \citep{Cotter2013SS,Beskos2015SC,Jia2021SVGD} for examples. 
\end{enumerate}

To overcome these obstacles, Bayes' formula defined on some separable Banach space is employed to handle the inverse problems of PDEs \citep{Dashti2014,Stuart2010AN}, which is an active research topic in the field of nonparametric Bayesian inference. We refer to \citep{Ghosal2017book} for a general account of the nonparametric Bayesian inference approach. Recently, the theory of well-posedness of the infinite-dimensional Bayesian analysis method has been generalized to incorporate some machine learning problems \citep{Latz2020SIAMUQ}. Based on this infinite-dimensional Bayesian analysis theory, graph-based Bayesian semi-supervised learning algorithms have been analyzed from many different aspects, such as large data and zero noise limit \citep{Dunlop2020ACHA}, as well as the consistency and scalability of sampling algorithms \citep{Hoffmann2020JMLR,Trillos2020JMLR}. As new surrogate models of PDEs, operator learning methods and theories have been studied in \citep{Bhattacharya2021SMAI-JCM,Zongyi2020ICLR,Kovachki2021JMLR}, which play essential roles in reducing the computational complexity of Bayesian inverse methods. The infinite-dimensional Bayesian inference method proposed for solving inverse problems of PDEs has been validated to be useful for some machine learning problems. Reciprocally, machine learning techniques can also contribute to the generalization of the infinite-dimensional Bayesian theory. Generally speaking, there are three main research directions for infinite-dimensional Bayesian analysis methods: designing prior probability measures on infinite-dimensional spaces, constructing appropriate noise models, and extracting information from posterior measures efficiently. In the following, we review some typical works along these three research directions and point out their essential difficulties in applications.

Regarding prior measures defined on the infinite-dimensional space, the Gaussian measure is prevalent based on its abundant theoretical studies in the field of stochastic PDEs \citep{Prato2014book}. Related studies on employing Gaussian prior measures can be found in \citep{agapiou2014SIAM,Tan2013IP,Thanh2016IPI,Cotter2009IP,Jia2020SIAM,Wang2019SISC}. To characterize the discontinuity of the function parameters, the Besov type prior measure has been proposed \citep{Lassas2009IPI,Dashti2012IPI,Jia2016IP}. Recently, the well-posedness of the Bayesian inverse method under prior measures with exponential tails has been analyzed in \citep{Hosseini2017SIAMASA}.
Besides, starting from the seminal work \citep{Knapik2011AS}, the posterior consistency and contraction rates have been analyzed in detail in a series of works \citep{Agapiou2013IP,Giordano2020IP,Jia2021IPI,Kekkonen2016IP,Szabo2015AS,Vollmer2013IP}, which convey a general understanding of what types of priors make a Bayesian nonparametric method effective. However, as illustrated in \citep{Adler2018arXiv,Arridge2019ActaNumerica}, most of these priors need to be pre-designed by hand-craft and much attention has been paid to characterize families of priors to ensure posterior consistency and good convergence rates.
This makes these priors more often chosen as a regularizer rather than essentially improving the quality of the final output \citep{Arridge2019ActaNumerica}. In addition, although the hand-crafted infinite-dimensional prior measure encodes some intuitive ideas of the prior information, it can only reflect some rough prior information, e.g., the function parameter is smooth, changes slowly along certain directions, and may have sharp edges. Thus, it is desirable to adaptably encode a prior specifically suitable for a particular input in data-driven manners \citep{Adler2018arXiv}.

For the noise model, the Gaussian measure is frequently used in the research literature, e.g., the Gaussian noise setting is employed in \citep{Tan2013IP,Thanh2016IPI,Jia2018JFA,Wang2019SISC}. By characterizing the model error as a Gaussian random variable, a Bayesian method with approximate error can be constructed \citep{Kaipio2004Book,Kaipio2019IP}. A natural generalization was proposed in \citep{Jia2019IP}, which employed a mixture of Gaussian distributions to characterize model noises, making it a better fit for a broader range of noise types due to the universe approximation capability of the mixture of Gaussian for general distributions. Through a learning process, the model error information was encoded into the parameters of Gaussian mixture distributions, which improved the estimated quality by only using an approximate forward solver. However, the method proposed in \citep{Jia2019IP} only provides a maximum a posterior estimate, which can be seen as an incomplete Bayesian approach. Furthermore, they have not considered more general and practical non-independently and non-identically distributed (non-i.i.d.) noises.

Furthermore, extracting information from the posterior measure is one of the essential issues for employing Bayes' method. Sampling algorithms, such as the Markov chain Monte Carlo, are often employed. The preconditioned Crank--Nicolson algorithm defined on some infinite-dimensional space has been proposed \citep{Beskos2008SD} to ensure the robustness of the convergence speed under mesh refinement, see also \citep{Cotter2013SS}. Then, multiple dimension-independent Markov chain Monte Carlo type sampling algorithms have been proposed \citep{Agapiou2017SS,Beskos2015SC,Cui2016JCP,Feng2018SIAM}. However, as illustrated in \citep{Arridge2019ActaNumerica}, the current sampling algorithms still face the critical issue of computational efficiency, especially for large-scale inverse problems.

Under the finite-dimensional setting, variational inference (VI) methods have been widely studied in machine learning \citep{Zhang2018IEEE} to reduce the computational burden of Markov chain Monte Carlo type sampling algorithms. Compared to finite-dimensional problems, however, infinite-dimensional problems have been much less studied for VI methods. When the approximate measures are restricted to be Gaussian, a novel Robbins-Monro algorithm was developed in \citep{Pinski2015SIAMMA,Pinski2015SIAMSC} from a calculus-of-variations viewpoint. It was shown in \citep{Sun2019ICLR} that the Kullback-Leibler (KL) divergence between the stochastic processes is equal to the supremum of the KL divergence between the measures restricted to finite marginals. Meanwhile, they developed a VI method for functions parameterized by Bayesian neural networks. Under the classical mean-field assumption, a general VI framework defined on separable Hilbert space was proposed recently in \citep{Jia2020SIAM}. A function space particle optimization method, including the Stein variational gradient descent, was developed in \citep{Wang2019ICLR} to solve the particle optimization directly in the space of functions. Infinite-dimensional Stein variational gradient descent with preconditioning operators has further been proposed in \citep{Jia2021SVGD}, which provides a detailed mathematical analysis in some separable Hilbert space. These VI methods defined on infinite-dimensional space have been proposed to train deep neural networks (DNNs) or solve inverse problems of PDEs, which indeed ameliorate the computational complexity to a certain extent. However, for solving inverse problems of PDEs, many forward and adjoint PDEs still need to be solved, which is time-consuming. Furthermore, as illustrated in \citep{Arridge2019ActaNumerica}, these VI methods also require manual pre-design of priors with explicit forms, which are not always sufficiently flexible and adaptable to faithfully deliver complicated prior knowledge underlying the investigated problem.

Through the above review and discussions, we summarize the main points as follows:
\begin{enumerate}
	\item To more efficiently and adaptably encode prior information (not just for regularizing), machine learning methods with data-driven manners are expected to be involved in the infinite-dimensional setting.
	\item The mixture of Gaussian noise modeling strategies has been introduced to improve the noise fitting capability of the problem by current research, but they are still restrictive on i.i.d. noise types, while more practical non-i.i.d. noise cases are expected to be considered.
	\item Although infinite-dimensional VI (IDVI) methods have been proposed recently, time-consuming iterations are still inevitably required. More efficient regimes are thus required to be investigated. Besides, it is highly expected for current IDVI methods to flexibly and faithfully specify prior information directly from data rather than commonly adopted hand-crafted manners.
\end{enumerate}

To address these issues, inspired by recent investigations on image denoising and restorations \citep{Yue2019NIPS,Yue2020VariationalIRN}, we aim to propose a general IDVI method that intends to integrate prior information learning, noise information learning, and posterior mean and covariance learning into a unified machine learning framework. Unlike image denoising and restoration problems, it is expensive to generate high-quality training datasets for many inverse problems of PDEs. Therefore, it is essential for us to incorporate model information with the data-driven machine learning framework.

In our opinion, the models contain two aspects: the physical model and the probabilistic model. For the physical model, we mean the PDEs that encode our knowledge of the forward process. For the probabilistic model, we indicate the elements of Bayesian methods, such as the prior and noise measures, which encode the structures we employ to solve the inverse problems of PDEs. In this work, we focus on linear inverse problems and intend to construct a generating model (see \citep{Kingma2017Thesis,Yue2019NIPS,Yue2020VariationalIRN} for examples) that could incorporate the constraint of the two models mentioned above. Explicitly speaking, we employ the mean-field IDVI theory to derive the general forms of the approximate posterior measure that provides the probabilistic model. For using the mean-field IDVI theory, we assume the noises of data are non-i.i.d. in statistics, and the prior measure is Gaussian with the synthetic background truth as the mean function. It should be pointed out that the mean-field IDVI theory here only provides a probabilistic model but not a practical algorithm since the synthetic background truth is not known in the inference stage. The probabilistic model provides the general structure of our methods (see Subsections \ref{VIclassical} and \ref{VINetTheory} for details), then we parameterize the parameters (some parameters are functions) of the approximate posterior measure by some well-designed DNNs that incorporate the physical model in some specially designed layers (see Subsections \ref{VINetTheory} and \ref{NetworkStructure} for details). Through the training process, the prior and noise information will be encoded into the parameters of neural networks, i.e., the general forms of the posterior measure are determined through IDVI theory but the function- and vector-valued parameters are learned from training datasets. During the inference stage, there will be no explicit Bayesian formulation. We only need to directly feed the new data to the trained deep neural network, and then we will easily obtain the parameters of the approximate posterior measure. The prior and noise information are then not required to be encoded in explicit prior and noise measures with specific parameters but encoded in the parameters of the DNNs with probabilistic model constraint. Numerical examples are given to illustrate the flexibility and effectiveness of the proposed approach. In summary, this work mainly contains four contributions: 

\begin{enumerate}
\item Inspired by the investigations on generative models in the finite-dimensional setting \citep{Kingma2017Thesis,Yue2019NIPS,Yue2020VariationalIRN},
we propose a general Bayesian generative model on infinite-dimensional space named as \textbf{V}ariational \textbf{I}nverting \textbf{Net}work (VINet),
which is naturally derived from the analysis on the IDVI method for linear inverse problems with non-i.i.d. noises.
Due to the general assumptions on the noises and priors, the proposed VINet integrates noise information learning, prior
information learning, and posterior statistical estimation into a unified machine learning framework.
\item Instead of only focusing on the interested function parameters as conventional works like deep learning for elastic source imaging \citep{Yoo2018SIAP} and deep Bayesian inversion \citep{Adler2018arXiv}, the proposed VINet provides uncertainty information both on the function parameters and on the random noises, allowing the noise and function parameter estimations to be mutually ameliorated during the inverse computational stage.
\item By imposing equivalence conditions on the prior and the approximate posterior measure, we derive the evidence lower bound under the infinite-dimensional setting. A method of amortized variational inference is developed to train all variables in the generative model. The model can then be explicitly used in the inference stage to readily achieve the posterior information of any noisy data without the time-consuming iteration process, making it more efficient for the inverse problems of PDEs than the conventional paradigms illustrated in \citep{Dashti2014}. 
\item We provide a specific parametric strategy by DNNs, which has been proven to satisfy the requirements of the general inference theory. Based on the proposed parametric strategy and the general abstract theory, we design a rational specific structure of VINet. Two numerical examples of typical inverse problems of PDEs are given to substantiate the effectiveness of the proposed method.
\end{enumerate}

The rest of the paper is outlined as follows. In Section \ref{generalTheory}, the general theory is presented.
Specifically, in Subsection \ref{problemsetup}, the general settings of the employed Bayesian inference model is given.
In Subsection \ref{VIclassical}, the Bayesian inference model is analyzed under the mean-field based infinite-dimensional VI framework. In Subsection \ref{VINetTheory}, we construct the general Bayesian generative model (i.e., VINet), derive the amortized variational inference method, and provide a simple preliminary theoretical analysis.
In Subsection \ref{ParametizationVINet}, we introduce a specific parametric strategy that can parameterize the posterior. In Section \ref{NumericalExamples}, two numerical examples are presented to illustrate the effectiveness of the proposed method. In particular, in Subsection \ref{NetworkStructure}, we design a specific network structure
according to our general theory, and then apply the proposed VINet to a simple toy smoothing model in Subsection \ref{SimpleModel}. We then demonstrate a more realistic example that is an inverse source problem governed by the Helmholtz equation in Subsection \ref{HelmholtzModel}. In Section \ref{Conclusion}, we summarize our results and propose some directions for future research.

\section{Variational Inverting Network}\label{generalTheory}

\begin{figure}[t]
	\centering
	\includegraphics[width=1\textwidth]{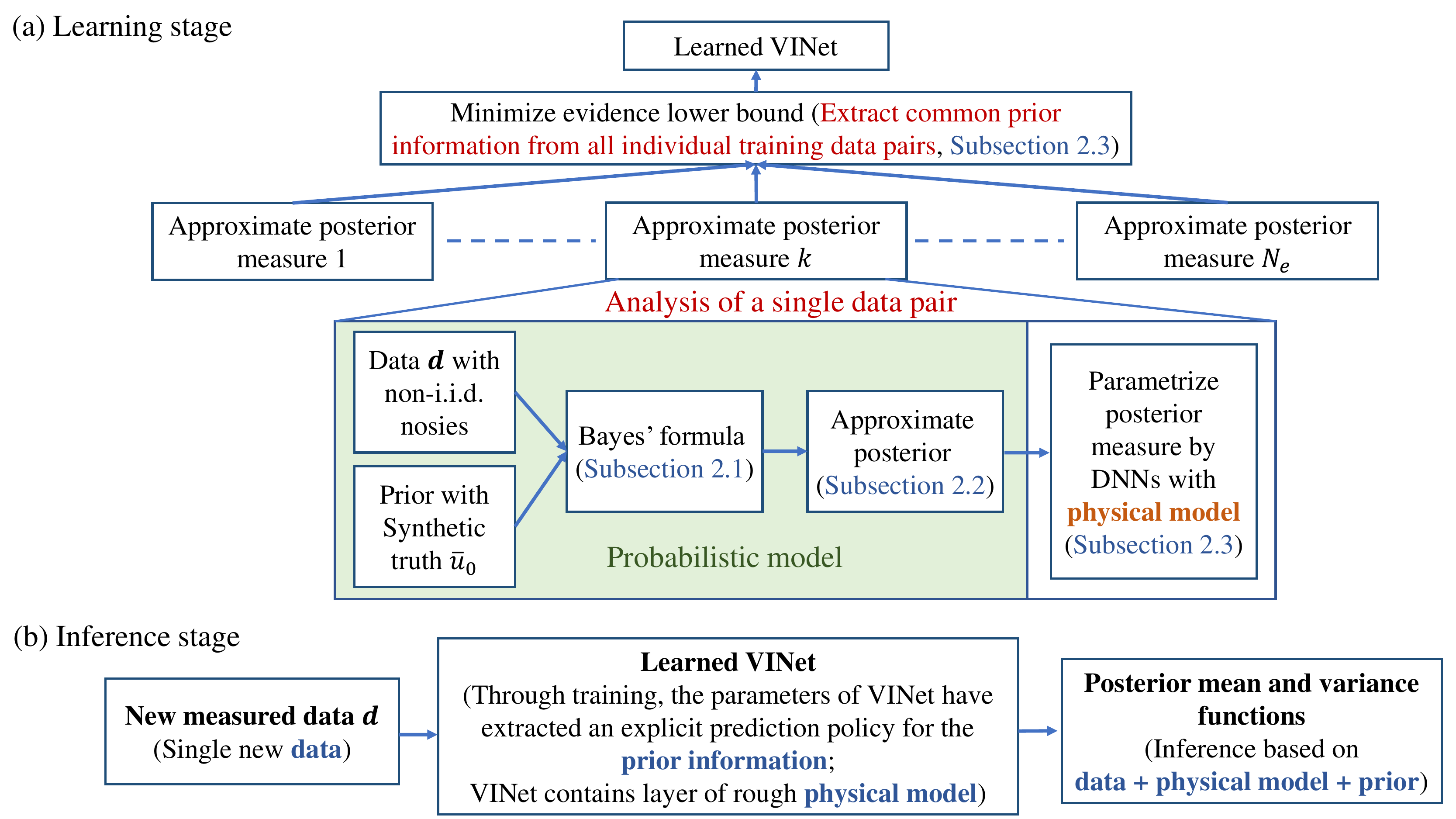}\\
	\vspace*{-0.6em}
	\caption{Overall structure of Section 2.}\label{GeneralFramework}
	\vspace*{0.5em}
\end{figure}

In this section, we construct the general inference model named VINet, which is well-defined on an infinite-dimensional space. Since the overall logical structure is a bit complex, we plan to give a general illustration before going into the details. In Figure \ref{GeneralFramework}, we exhibit the overall structure of the present work. Generally speaking, the method can be divided into two stages: (a) Learning stage, and (b) Inference stage.

In the learning stage, we formulate Bayes' formula (see Subsection 2.1) by assuming the noises are non-i.i.d. random variables and the prior measure of the function parameter is Gaussian. In particular, the mean function of the Gaussian measure of each training sample is set to its specific synthetic truth (these measures are thus different for different training samples), which is available in the training stage but not in the inference stage. Then, relying on the infinite-dimensional VI framework, we derive the approximate posterior measures (see Subsection 2.2), which encode the information of all noisy samples and their own specific synthetic truths. One of the key points of the present work is that the approximate posterior measures of all single training data pairs are then parameterized by a carefully designed deep neural network containing a rough physical model (details can be found in Subsection 2.3). Through training on all of the individual learning data pairs, the parameters of the deep neural network are expected to extract their underlying common prediction principle for the prior information.

Once the learning stage is complete, it will be convenient for us to make inferences when new measurement data is given. As illustrated in Figure \ref{GeneralFramework}, the new data (there will be no synthetic truth) will be fed to the trained VINet. Then, the trained VINet will generate the posterior information, such as the mean and variance functions.

\subsection{The Bayesian inference method}\label{problemsetup}

Let $\mathcal{H}_{u}$ be a separable Hilbert space, $N_d$ be a positive integer, and $H:\mathcal{H}_u\rightarrow\mathbb{R}^{N_d}$ be a bounded linear operator. Consider the  model
\begin{align}\label{addativeModelGe}
\bm{d} = Hu + \bm{\epsilon},
\end{align}
where $\bm{d}$ and $\bm{\epsilon}\in\mathbb{R}^{N_d}$ represent the measurement data and the random noise, respectively.
Denote by $\mathcal{C}_0$ a symmetric, positive definite, and trace class operator defined on $\mathcal{H}_u$. Let 
$(e_k, \alpha_k)$ be an eigen-system of the operator $\mathcal{C}_0$ such that $\mathcal{C}_0 e_k = \alpha_k e_k$. Without loss of generality, we may assume that the eigenvectors $\{ e_k \}_{k=1}^{\infty}$ are orthonormal, the eigenvalues $\{\alpha_k \}_{k=1}^{\infty}$
are sorted in a descending order, and the summation of eigenvalues is finite $\sum_{k=1}^{\infty}\alpha_k < \infty$. 
It follows from Subsection 2.4 in \citep{Dashti2014} that we have
\begin{align}\label{priorDecom}
\mathcal{C}_0 = \sum_{k=1}^{\infty}\alpha_{k}e_k\otimes e_k.
\end{align}
For the parameter $u$ and the noise $\bm{\epsilon}$, we assume
\begin{align}\label{uandepsilon}
u\sim \mu_0^u := \mathcal{N}(\bar{u}_0, \mathcal{C}_0), \quad
\bm{\epsilon} \sim \mathcal{N}(0, \Sigma),
\end{align}
where $\Sigma = \text{diag}(\sigma_1, \sigma_2, \cdots, \sigma_{N_d})$ with $\{\sigma_i\}_{i=1}^{N_d}$ being
a series of positive real numbers. Here, $\mathcal{N}(v, \mathcal{C})$ stands for a Gaussian measure with mean $v$ and covariance operator $\mathcal{C}$. Clearly, $\mathcal{N}(\bar{u}_0, \mathcal{C}_0)$ denotes a Gaussian measure defined on $\mathcal{H}_u$ with mean $\bar{u}_0\in\mathcal{H}_u$ and covariance operator $\mathcal{C}_0$, and $\mathcal{N}(0,\Sigma)$ is a Gaussian measure defined on $\mathbb{R}^{N_d}$. It is worth mentioning that the mean function $\bar{u}_0$ as indicated in Figure \ref{GeneralFramework} is set to be the synthetic truth. By our understanding, the synthetic truth is similar to, but not exactly the same as, the background truth (see also Remark \ref{remarkunknonTruth} in Subsection \ref{ParametizationVINet}). The background truth can hardly be obtained \citep{Yue2019NIPS,Davoudi2019NMI}. Since the synthetic truth is set to be the mean function, we usually take the magnitude of the covariance operator $\mathcal{C}_0$ to be small, which reflects the deviation of the synthetic truth from the background truth. For one specific form of the operator $\mathcal{C}_0$, please see Subsection \ref{ParametizationVINet} and Section \ref{NumericalExamples}.

For the additive noise model (\ref{addativeModelGe}), the data $\bm{d}$ given $u$ are distributed according to the following translated distribution of $\bm{\epsilon}$:
\begin{align}
\bm{d}|u \sim \mathcal{N}(Hu, \Sigma).
\end{align}
Usually, we assume that $\Sigma$ is known a priori when making inferences on the infinite-dimensional parameter $u$. However, similar to studies such as in \citep{Dunlop2017SC,Jin2010JCP}, we assume that the parameters $\{\sigma_i\}_{i=1}^{\infty}$ are unknown random variables, which are called hyper-parameters and also need to be estimated. The key elements of an abstract inverse problem can then be summarized as follows.
\begin{itemize}
	\item \textbf{Noisy data}: The data vector $\bm{d}\in\mathbb{R}^{N_d}$ is obtained from $\bm{d} = Hu^{\dag} + \bm{\epsilon}$, where $u^{\dag}$ is the background truth and $\bm{\epsilon}\sim\mathcal{N}(0,\Sigma)$ with $\Sigma = \text{diag}(\sigma_1, \cdots, \sigma_{N_d})$;
	\item \textbf{Inverse problem}: Based on the noisy data $\bm{d}$, the goal is to find estimates of $u$ and the hyper-parameters
	$\bm{\sigma} := \{\sigma_1,\ldots,\sigma_{N_d}\}$.
\end{itemize}
We assume that the hyper-parameters in $\Sigma$ satisfy
\begin{align}\label{gammaDef}
\bm{\sigma} \sim \mu_0^{\bm{\sigma}}:=\prod_{i=1}^{N_d}\text{IG}\left( \alpha^0_i, \beta^0_i \right),
\end{align}
where $\alpha_i^0 > 1$ and $\beta_i^0 > 0$ ($i=1,\ldots,N_d$) are hyper-parameters that will be specified in numerical examples,
and $\text{IG}(\alpha_i^0,\beta_i^0)$ denotes the inverse Gamma distribution with parameters $\alpha_i^0$ and $\beta_i^0$ ($i=1,\ldots,N_d$).
Here, the inverse Gamma distribution is employed since it is a conjugate prior measure that yields the same form of approximate posterior measure (see Subsection \ref{VIclassical}). Conjugate prior measures are widely used in the studies of inverse problems and machine learning \citep{Kekkonen2016IP,Zhang2018IEEE,Yue2020VariationalIRN}, which makes it convenient to construct computational methods.

In order to construct an appropriate Bayes' formula, let us introduce the prior probability measure
\begin{align}\label{priorMeasure}
\mu_0 := \mu_0^u \otimes \mu_0^{\bm{\sigma}}
\end{align}
and the potential function
\begin{align}\label{defPhi}
\Phi(u,\bm{\sigma}; \bm{d}) := \frac{1}{2} \|\bm{d} - Hu\|_{\Sigma}^2 + \frac{1}{2} \log(\det\Sigma),
\end{align}
where $\|\cdot\|_{\Sigma} := \|\Sigma^{-1/2}\cdot\|$ and $\det\Sigma$ is the determinant of $\Sigma$.
Since $H$ is a bounded linear operator, it is easy to note that $\Phi(u, \bm{\sigma}; \bm{d})$ is continuous with respect to the variables $u, \bm{\sigma}$, and $\bm{d}$. Assume $\bm{d}, \bm{d}_1$, and $\bm{d}_2$ are included in a ball of $\mathbb{R}^{N_d}$ with radius $r$ and centered at the origin. Following simple calculations, we have  
\begin{align}
\begin{split}
\Phi(u,\bm{\sigma};\bm{d}) = & \frac{1}{2}\|\bm{d} - Hu\|_{\Sigma}^2 + \frac{1}{2}\log(\det\Sigma) \geq \frac{1}{2}\log\left( \prod_{k=1}^{N_d}\sigma_k \right)
= \frac{1}{2}\sum_{k=1}^{N_d}\log(\sigma_k)
\end{split}
\end{align}
and
\begin{align}
\begin{split}
|\Phi(u,\bm{\sigma};\bm{d}_1) - \Phi(u,\bm{\sigma};\bm{d}_2)| &=  \frac{1}{2}|\|\bm{d}_1 - Hu\|_{\Sigma}^2 - \|\bm{d}_2 - Hu\|_{\Sigma}^2| \\
&\leq  \frac{1}{2}\|\bm{d}_1 + \bm{d}_2 - 2Hu\|_{\Sigma}\|\bm{d}_1 - \bm{d}_2\|_{\Sigma} \\
&\leq  \left(1 + \sum_{k=1}^{N_d}\sigma_k^{-1/2}\right)\sum_{\ell=1}^{N_d}\sigma_{\ell}^{-1/2}(r + C\|u\|_{\mathcal{H}_u})\|\bm{d}_1 - \bm{d}_2\|,
\end{split}
\end{align}
where $C$ is a constant larger than the operator norm of $H$. Noting
\begin{align}
\exp\left(-\sum_{k=1}^{N_d}\log\sigma_k \right) = \prod_{k=1}^{N_d}\sigma_k^{-1} \in L_{\mu_0}^1(\mathcal{H}_u\times\mathbb{R}^{+};\mathbb{R})
\end{align}
and
\begin{align}
\prod_{k=1}^{N_d}\sigma_k^{-1}\left(1 + \sum_{k=1}^{N_d}\sigma_k^{-1/2}\right)^2\left(\sum_{\ell=1}^{N_d}\sigma_{\ell}^{-1/2}\right)^2(r + C\|u\|_{\mathcal{H}_u})^2
\in L_{\mu_0}^1(\mathcal{H}_u\times\mathbb{R}^{+};\mathbb{R}),
\end{align}
we know that Assumption 1 and the conditions of Theorems 15--16 in \citep{Dashti2014} are satisfied,
which provide the following Bayes' formula defined on infinite-dimensional space:
\begin{align}\label{mainBayesFormula}
\frac{d\mu^{\bm{d}}}{d\mu_0}(u,\bm{\sigma}) = \frac{1}{Z_{\bm{d}}}\text{det}(\Sigma)^{-1/2}\exp\left( -\frac{1}{2}\|\bm{d}-Hu\|_{\Sigma}^2 \right),
\end{align}
where $\mu^{\bm{d}}$ represents the posterior probability measure and
\begin{align}
Z_{\bm{d}} = \int_{(\mathbb{R}^{+})^{N_d}}\int_{\mathcal{H}_u} \text{det}(\Sigma)^{-1/2}
\exp\left( -\frac{1}{2}\|\bm{d}-Hu\|_{\Sigma}^2 \right)
\mu_0^{u}(du)\mu_0^{\bm{\sigma}}(d\bm{\sigma}).
\end{align}
For the reader's convenience, we list Assumption 1 and Theorems 15--16 of \citep{Dashti2014} in Subsection \ref{BayesWellPosed} of the Appendix.

\subsection{The infinite-dimensional variational inference method}\label{VIclassical}

Before constructing the VINet, we discuss the mean-field assumption
based infinite-dimensional variational inference (IDVI) theory developed in \citep{Jia2020SIAM}.
A brief introduction is given in Subsection \ref{AppendixGeneTheoryVI} of the Appendix.

In Subsection \ref{problemsetup}, we have not introduced hyper-parameters for the prior measure $\mu_0^u$ of the function parameter.
Hence, we can choose the prior measure $\mu_0^u$ as the reference probability measure required in Assumption 8 in \citep{Jia2020SIAM}. Using the mean-field approximation, i.e., assuming the parameters $u$ and $\bm{\sigma}$ to be independent random variables, we introduce an approximate probability measure $\nu(du, d\bm{\sigma}) = \nu^{u}(du)\nu^{\bm{\sigma}}(d\bm{\sigma})$ with
\begin{align}\label{nuudef}
\frac{d\nu^{u}}{d\mu_0^{u}}(u)=\frac{1}{Z_{u}}\exp\left( -\Phi_u(u) \right)
\end{align}
and
\begin{align}\label{nusigmadef}
\frac{d\nu^{\bm{\sigma}}}{d\mu_0^{\bm{\sigma}}}(\bm{\sigma}) = \frac{1}{Z_{\bm{\sigma}}}\exp\left( -\Phi_{\bm{\sigma}}(\bm{\sigma}) \right),
\end{align}
where $\Phi_{u}(\cdot)$ and $\Phi_{\bm{\sigma}}(\cdot)$ are two potential functions required to be calculated explicitly, and
\begin{align*}
Z_{u} = \mathbb{E}^{\mu_0^u}\left[ \exp(-\Phi_{u}(u)) \right], \quad
Z_{\bm{\sigma}} = \mathbb{E}^{\mu_0^{\bm{\sigma}}}\left[ \exp(-\Phi_{\bm{\sigma}}(\bm{\sigma})) \right].
\end{align*}
Here, $\mathbb{E}^{\mu}$ represents taking expectation with respect to the probability measure $\mu$.

The VI method essentially needs to solve the optimization problem
\begin{align}\label{optimizationVI}
\argmin_{\nu^u\in\mathcal{A}_u, \nu^{\bm{\sigma}}\in\mathcal{A}_{\bm{\sigma}}}D_{\text{KL}}\left( \nu || \mu^{\bm{d}} \right),
\end{align}
where $\mathcal{A}_u$ and $\mathcal{A}_{\bm{\sigma}}$ are two appropriate spaces of probability measure, and $D_{\text{KL}}\left( \nu || \mu^{\bm{d}}\right)$ is the Kullback-Leibler (KL) divergence defined as follows:
\begin{align*}
D_{\text{KL}}(\nu||\mu^{\bm{d}}) =
\int_{\mathcal{H}}\log\bigg( \frac{d\nu}{d\mu^{\bm{d}}}(x) \bigg)\frac{d\nu}{d\mu^{\bm{d}}}(x)\mu^{\bm{d}}(dx)
= \mathbb{E}^{\mu^{\bm{d}}}\bigg[ \log\bigg( \frac{d\nu}{d\mu^{\bm{d}}}(x) \bigg)\frac{d\nu}{d\mu^{\bm{d}}}(x) \bigg].
\end{align*}
Here $0\log 0 = 0$ is used as a convention. For the finite-dimensional theory (e.g., Chapter 10 in \citep{Bishop2006book}), no additional assumptions need to be made on $\mathcal{A}_u$ and $\mathcal{A}_{\bm{\sigma}}$. However, special attention must be paid to $\mathcal{A}_u$ and $\mathcal{A}_{\bm{\sigma}}$ in the infinite-dimensional theory \citep{Jia2020SIAM}.

Define $T_N^u = \left\{ u \, | \, 1/N \leq \|u\|_{\mathcal{Z}_u} \leq N \right\}$ with $\mathcal{Z}_u$ being a Hilbert space
that is embedded in $\mathcal{H}_u$ and satisfies $\sup_{N}\mu_0^u(T_N^u) = 1$.
Denote $T_N^{\bm{\sigma}} = \left\{ \bm{\sigma} \, | \, 1/N \leq \sigma_i \leq N \text{ for }i=1,\ldots,N_d \right\}$ that
obviously satisfies $\sup_N \mu_0^{\bm{\sigma}}(T_N^{\bm{\sigma}}) = 1$. Let 
\begin{align*}
R_u^1 = \left\{ \Phi_u \, \Big| \,
\sup_{u\in T_N^u}\Phi_u(u) < \infty \text{ for all }N > 0
\right\},
\end{align*}
\begin{align*}
R_u^2 = \left\{ \Phi_u \, \Big| \,
\int_{\mathcal{H}_u}\exp\left(-\Phi_u(u)\right)\max(1, \|u\|_{\mathcal{H}_u}^2)\mu_0^u(du) < \infty
\right\},
\end{align*}
\begin{align*}
R_{\bm{\sigma}}^1 = \left\{ \Phi_{\bm{\sigma}} \, \Big| \,
\sup_{\bm{\sigma}\in T_{N}^{\bm{\sigma}}}\Phi_{\bm{\sigma}}(\bm{\sigma}) < \infty \text{ for all }N > 0
\right\},
\end{align*}
\begin{align*}
R_{\bm{\sigma}}^2 = \left\{ \Phi_{\bm{\sigma}} \, \Big| \,
\int_{(\mathbb{R^+})^{N_d}}  \exp\left( -\Phi_{\bm{\sigma}}(\bm{\sigma}) \right)
\max(1, a(\epsilon, \bm{\sigma}))\mu_0^{\bm{\sigma}}(d\bm{\sigma}) < \infty, \text{ for }\epsilon\in[0,\epsilon_0)
\right\},
\end{align*}
where $a(\epsilon, \bm{\sigma}) := \sum_{i=1}^{N_d}\max\left( \sigma_i^{\epsilon}, \exp(\epsilon/\sigma_i) \right)$
and $\epsilon_0$ is a fixed positive number. Usually, we may choose $\epsilon_0 \ll \min_{1\leq i\leq N_d}\{\alpha^0_i,\beta^0_i\}$, where $\{\alpha_i^0,\beta_i^0\}_{i=1}^{N_d}$ are given in (\ref{gammaDef}). Define 
\begin{align*}
\mathcal{A}_u = \left\{ \nu^u \in \mathcal{M}(\mathcal{H}_{u}) \, \Big| \,
\begin{array}{l}
\nu^u \text{ is equivalent to }\mu_0^u\text{ with (\ref{nuudef}) holding true,}\\
\text{and }\Phi_u\in R_u^1\cap R_u^2
\end{array}
\right\},
\end{align*}
\begin{align*}
\mathcal{A}_{\bm{\sigma}} = \left\{ \nu^{\bm{\sigma}} \in \mathcal{M}(\mathcal{H}_{\bm{\sigma}}) \, \Big| \,
\begin{array}{l}
\nu^{\bm{\sigma}} \text{ is equivalent to }\mu_0^{\bm{\sigma}}\text{ with (\ref{nusigmadef}) holding true,}\\
\text{and }\Phi_{\bm{\sigma}}\in R_{\bm{\sigma}}^1\cap R_{\bm{\sigma}}^2
\end{array}
\right\}.
\end{align*}

\begin{remark}\label{R1_def_shuoming}
	The definitions of $R_u^1$ and $R_{\bm{\sigma}}^1$ are slightly different from the corresponding definitions
	of $R_j^1(j=1,\ldots,M)$ introduced in \citep{Jia2020SIAM}.
	The sets $R_k^1(k=u,\bm{\sigma})$ defined here enlarge the original definition (without uniform bound), which makes the theory more appropriate. For completeness, we provide a short illustration of the modified general theory in Subsection \ref{AppendixGeneTheoryVI} of the Appendix.
\end{remark}

Next, we present the key theorem that leads to an iterative algorithm.

\begin{theorem}\label{keyVItheorem}
	Assume that the prior measure $\mu_0$, noise measure, and posterior measure are defined in (\ref{priorMeasure}), (\ref{gammaDef}), and (\ref{mainBayesFormula}), respectively. Let $\Phi(u,\bm{\sigma};\bm{d})$ be defined in (\ref{defPhi}).
	If $\mathcal{A}_u$ and $\mathcal{A}_{\bm{\sigma}}$ are defined as above, then the problem (\ref{optimizationVI}) possesses a solution $\nu(du,d\bm{\sigma}) = \nu^u(du)\nu^{\bm{\sigma}}(d\bm{\sigma})$,  which satisfies 
	\begin{align}
	\frac{d\nu}{d\mu_0}(u,\bm{\sigma}) \propto \exp\left( -\Phi_u(u) - \Phi_{\bm{\sigma}}(\bm{\sigma}) \right),
	\end{align}
	where
	\begin{align}
	\Phi_u(u) & = \int_{(\mathbb{R}^{+})^{N_d}}\Phi(u, \bm{\sigma}; \bm{d})\nu^{\bm{\sigma}}(d\bm{\sigma}) + \text{Const},  \label{approUC} \\
	\Phi_{\bm{\sigma}}(\bm{\sigma}) & = \int_{\mathcal{H}_u}\Phi(u, \bm{\sigma}; \bm{d})\nu^{u}(du) + \text{Const}. \label{approSC}
	\end{align}
	Here, ``Const'' denotes some general constant. Futhermore, we have
	\begin{align*}
	\nu^{u}(du) \propto \exp\left( -\Phi_u(u) \right)\mu_0^u(du), \quad
	\nu^{\bm{\sigma}}(d\bm{\sigma}) \propto \exp\left( -\Phi_{\bm{\sigma}}(\bm{\sigma}) \right)\mu_0^{\bm{\sigma}}(d\bm{\sigma}).
	\end{align*}
\end{theorem}

To avoid a possible distraction from the presentation of the work, the proof of this theorem is given in Subsection \ref{AppendixProof} of the Appendix. In the following, we denote $\Phi(u, \bm{\sigma})$ instead of $\Phi(u, \bm{\sigma}; \bm{d})$ when there is no ambiguity from the context.

\begin{remark}\label{transformRemark}
	In the general theory, the space of parameters is assumed to be a separable Hilbert space. The parameter $\bm{\sigma}$ resides in $(\mathbb{R}^{+})^{N_d}$, which is not a Hilbert space. However, as stated in Remark 15 in \citep{Jia2020SIAM}, a simple transformation can be adopted, e.g., $\sigma_k' = \log\sigma_k$ ($k=1,\ldots,N_d$). Then, the new parameter $\bm{\sigma}' = \{ \sigma_1',\ldots,\sigma_{N_d}' \}$ belongs to $\mathbb{R}^{N_d}$, which is a Hilbert space. All of the statements and verifications of $\bm{\sigma}$ can be equivalently transformed to the new variable $\bm{\sigma}'$. Therefore, we still use the parameter $\bm{\sigma}$ for simplicity of notation. By the above discussion, the reader may intuitively understand the connection between the finite- and infinite-dimensional theories.
\end{remark}

Using Theorem \ref{keyVItheorem}, we can present a classical iterative algorithm when the parameter $u$ belongs to some infinite-dimensional separable Hilbert space $\mathcal{H}_u$.

\textbf{Calculate $\Phi_u(u)$}. A direct application of (\ref{approUC}) yields
\begin{align}
\begin{split}
\Phi_u(u) &=  \frac{1}{2}\int_{(\mathbb{R}^{+})^{N_d}}
\left( \|\bm{d} - Hu\|_{\Sigma}^2 + \sum_{k=1}^{N_d}\log\sigma_k \right)\nu^{\bm{\sigma}}(d\bm{\sigma}) + \text{Const} \\
&=  \frac{1}{2}\int_{(\mathbb{R}^{+})^{N_d}} (\bm{d}-Hu)^T\Sigma^{-1}(\bm{d}-Hu)\nu^{\bm{\sigma}}(d\bm{\sigma}) + \text{Const}.
\end{split}
\end{align}
Denoting
\begin{align}\label{221}
\Sigma_{\text{inv}}^* = \text{diag}\left( \mathbb{E}^{\nu^{\bm{\sigma}}}[\sigma_1^{-1}], \ldots,
\mathbb{E}^{\nu^{\bm{\sigma}}}[\sigma_{N_d}^{-1}] \right),
\end{align}
we have
\begin{align}
\Phi_u(u) = \frac{1}{2}\|\bm{d}-Hu\|_{\Sigma_{\text{inv}}^*}^2 + \text{Const}.
\end{align}
Hence, 
\begin{align}
\frac{d\nu^{u}}{d\mu_0^u}(u) \propto \exp\left( -\frac{1}{2}\|\bm{d}-Hu\|_{\Sigma_{\text{inv}}^*}^2 \right),
\end{align}
which implies that the approximate posterior measure of $u$ is a Gaussian measure
\begin{align}\label{CpostU}
\nu^{u} = \mathcal{N}(\bar{u}_p, \mathcal{C}_p),
\end{align}
where
\begin{align}\label{CposU1}
\mathcal{C}_p^{-1} = H^*\Sigma_{\text{inv}}^* H + \mathcal{C}_0^{-1}, \quad
\bar{u}_p = \mathcal{C}_p(H^* \Sigma_{\text{inv}}^* \bm{d} + \mathcal{C}_0^{-1}\bar{u}_0).
\end{align}

\textbf{Calculate $\Phi_{\bm{\sigma}}(\bm{\sigma})$}.
Following (\ref{approSC}), we have
\begin{align}
\begin{split}
\Phi_{\bm{\sigma}}(\bm{\sigma}) & =  \frac{1}{2}\int_{\mathcal{H}_u}
\left( \|\bm{d} - Hu\|_{\Sigma}^2 + \sum_{k=1}^{N_d}\log\sigma_k \right)\nu^{u}(du) + \text{Const} \\
& =  \frac{1}{2}\int_{\mathcal{H}_u} (\bm{d}-Hu)^T\Sigma^{-1}(\bm{d}-Hu)\nu^{u}(du) +
\frac{1}{2}\sum_{k=1}^{N_d}\log\sigma_k + \text{Const}  \\
& =  \sum_{k=1}^{N_d}\frac{1}{2}\left( \int_{\mathcal{H}_u} (d_k - (Hu)_k)^2 \nu^u(du)\frac{1}{\sigma_k} + \log\sigma_k \right) + \text{Const}.
\end{split}
\end{align}
Hence, 
\begin{align}
\frac{d\nu^{\bm{\sigma}}}{d\mu_0^{\bm{\sigma}}}(\bm{\sigma}) \propto
\exp\left( -\sum_{k=1}^{N_d}\frac{1}{2}\left( \mathbb{E}^{\nu^u}[(d_k - (Hu)_k)^2]\frac{1}{\sigma_k} + \log\sigma_k \right) \right).
\end{align}
Recalling 
\begin{align}
\mu_0^{\bm{\sigma}}(d\bm{\sigma}) \propto \prod_{k=1}^{N_d}\sigma_{k}^{-\alpha_k^0 - 1}\exp\left( -\frac{\beta_k^0}{\sigma_k} \right)d\bm{\sigma},
\end{align}
we obtain
\begin{align}\label{CpostS}
\nu^{\bm{\sigma}}(d\bm{\sigma}) \propto \prod_{k=1}^{N_d} \sigma_k^{-\alpha_k^0 - \frac{3}{2}}\exp\left(
-\frac{\beta_k^0 + \frac{1}{2}\mathbb{E}^{\nu^u}[(d_k - (Hu)_k)^2]}{\sigma_k} \right)d\bm{\sigma},
\end{align}
which implies that the posterior distribution of each component of $\nu^{\bm{\sigma}}$ takes the form
\begin{align}\label{noiseCEst}
\sigma_k \sim \text{IG}\left( \alpha_k^0 + \frac{1}{2}, \beta_k^0 + \frac{1}{2}\mathbb{E}^{\nu^u}[(d_k - (Hu)_k)^2] \right), \quad
k=1,\ldots,N_d.
\end{align}

The above illustrations confirm the feasibility of Algorithm \ref{classicalVI}, which is similar to the finite-dimensional case. However, it is worth mentioning that the infinite-dimensional formulation provides a general framework for conducting appropriate discretizations. To keep dimension-independent properties, the discretization of Algorithm \ref{classicalVI} should be done carefully; for example, the adjoint operator $H^*$ is usually not equal to the transpose of the discrete approximation of $H$. For more detailed discussions on this topic, please refer to \citep{Jia2020SIAM,Tan2013IP,Petra2014SIAM,Wang2019SISC,Thanh2016IPI}. A recent study on the Stein variational gradient descent algorithm, designed on the functional space for training DNNs \citep{Wang2019ICLR}, further indicates the infinite-dimensional formulation.

\begin{algorithm}
	\setstretch{1}
	\caption{A classical VI algorithm}
	\label{classicalVI}
	\begin{algorithmic}
		\STATE {1: Give an initial guess $\mu_{0}$ ($\bar{u}_0$, $\mathcal{C}_0$, $p$, and $\{\alpha_i^0,\beta_i^0\}_{i=1}^{N_d}$). Set the tolerance \emph{tol}  \\
			\quad\,\!\! and let $k=1$.
		}
		\STATE {2: \textbf{Do}
		}
		\STATE {3: \qquad Let $k = k + 1$
		}
		\STATE {4: \qquad Calculate $\Sigma_{\text{inv}}^{k*}$ by $\mathbb{E}^{\nu_{k-1}^{\bm{\sigma}}}[\sigma_i^{-1}]$($i=1,\ldots,N_d$)
		}
		\STATE {5: \qquad Update $\nu_k^u$ by
			\begin{align*}
			\mathcal{C}_p^{-1} = H^* \Sigma_{\text{inv}}^{k*} H + \mathcal{C}_0^{-1}, \quad
			u_k = \mathcal{C}_p (H^* \Sigma_{\text{inv}}^{k*} \bm{d} + \mathcal{C}_0^{-1}\bar{u}_0)
			\end{align*}
		}
		\STATE {6: \qquad Calculate $\nu_k^{\bm{\sigma}}$ by
			\begin{align*}
			\nu_k^{\bm{\sigma}} = \prod_{i=1}^{N_d}\text{IG}(\alpha_i^k,\beta_i^k),
			\end{align*}
			\quad\qquad\,\!\! where
			\begin{align*}
			\alpha_i^k = \alpha_i^0 + \frac{1}{2}, \quad \beta_i^k = \beta_i^0 + \frac{1}{2}\mathbb{E}^{\nu_k^{u}}[(d_i-(Hu_k)_i)^2]
			\end{align*}
		}
		\STATE{7: \textbf{Until} $\max\left( \|u_k-u_{k-1}\|/\|u_k\|, \|\bm{\sigma}_k-\bm{\sigma}_{k-1}\|/\|\bm{\sigma}_k\| \right) \leq tol$
		}
		\STATE{8: Return $\nu_{k}^{u}(du)\nu_{k}^{\bm{\sigma}}(d\bm{\sigma})$ as the solution.
		}
	\end{algorithmic}
\end{algorithm}

\begin{remark}\label{compuBurdenClassVI}
	In the implementation of Algorithm \ref{classicalVI}, it is required to calculate $\mathbb{E}^{\nu_k^{u}}[(d_i-(Hu_k)_i)^2]$ for each $i=1,\ldots,N_d$. Regarding these terms, we have
	\begin{align}\label{updateBeta}
	\begin{split}
	\mathbb{E}^{\nu_k^{u}}[(d_i-(Hu_k)_i)^2]& =  \int (\bm{d} - Hu_k)_i^2 \nu_k^u(du) +
	\int [H(u-u_k)]_i^2 \nu_k^u(du)  \\
	& =  (\bm{d}-Hu_k)_i^2 + \bm{e}_i^T H\mathcal{C}_pH^* \bm{e}_i,
	\end{split}
	\end{align}
	where $\bm{e}_i$, $i = 1, \ldots, N_d$, represents the standard vector basis in $\mathbb{R}^{N_d}$ and Proposition 1.18 in \citep{Prato2006book} is employed to derive the second equality. It can be seen from (\ref{updateBeta}) that the adjoint equations, forward equations, and inverse of the operator $\mathcal{C}_p$ are needed with $N_d$ iterations to update the parameters $\{\beta_i^k\}_{i=1}^{N_d}$. Especially, when calculating $\mathcal{C}_pu$ with $u$ being a function (e.g., $H^*\bm{e}_i$), we need to evaluate $(H^*\Sigma_{\text{inv}}^{k*}H + \mathcal{C}_0^{-1})^{-1}u$. The operator $(H^*\Sigma_{\text{inv}}^{k*}H + \mathcal{C}_0^{-1})^{-1}$ cannot be calculated explicitly since, in discrete form, it will be a large dense matrix that is difficult to store and time-consuming to find its explicit inverse operator. Instead, we usually solve a linear equation $(H^*\Sigma_{\text{inv}}^{k*}H + \mathcal{C}_0^{-1}) \tilde{u} = u$. The solution $\tilde{u}$ is the function we expected. By employing some iterative linear equation solvers, e.g., the conjugate gradient algorithm, we need to solve a lot of PDEs, i.e., for each iteration, we need to solve one forward and one adjoint equation \citep{Jia2020SIAM,Jin2010JCP}. Overall, the computational complexity will be very large when the number of measurement points $N_{d}$ increases.
\end{remark}

\subsection{Variational inverting network}\label{VINetTheory}

Let us provide the following observations of the classical IDVI algorithm:
\begin{itemize}
	\item The entire prior and noise probability measure (general forms and all of the function- and vector-valued parameters) must be specified manually. This makes the use of the algorithm highly dependent on the subjective experiences of the practitioners.
	\item Under the settings of Subsections \ref{problemsetup} and \ref{VIclassical}, the derived Algorithm \ref{classicalVI} actually still cannot be employed for practical problems since the synthetic truth (i.e., prior mean function) is not known for new inverse problems. 
	\item When the practitioner acquires new data, the entire iterative process should be started from the beginning to obtain a new estimate, which can be time-consuming for many applications.
\end{itemize}

Inspired by the recent work on denoising tasks under the finite-dimensional setting \citep{Yue2019NIPS,Yue2020VariationalIRN}, we introduce some parametric forms of posterior probability measure by adopting DNNs. As illustrated in the introduction and Figure \ref{GeneralFramework}, we will rely on the general form of the posterior measures derived from the IDVI theory, which provides the probabilistic model (or called probabilistic constraint) for DNNs-based data-driven inverse methods. However, the function- and vector-valued parameters in the posterior measures will be replaced by well-designed DNNs that could incorporate the physical model into our method. Relying on this methodology, it is then hopeful that we can alleviate the problems of the classical IDVI algorithms mentioned above.

Based on the discussions in Subsection \ref{VIclassical}, we have obtained the general form of the approximate posterior measures, as seen in formulas (\ref{CpostU}) and (\ref{CpostS}). Therefore, we introduce the following parametric form of the approximate posterior probability measures:
\begin{align}\label{postPM}
	\nu^{u} = \mathcal{N}(\bar{u}_p(\bm{d};W_I), \mathcal{C}_p(\bm{d};W_I)), \quad
	\nu^{\bm{\sigma}} = \prod_{i=1}^{N_d}\text{IG}\left( \alpha_i(\bm{d};W_S), \beta_i(\bm{d};W_S) \right),
\end{align}
where $\{\alpha_i(\bm{d};W_S)\}_{i=1}^{N_d}$, $\{\beta_i(\bm{d};W_S)\}_{i=1}^{N_d}$ are represented by neural networks with parameters denoted by $W_S$, and $\{\bar{u}_p(\bm{d};W_I), \mathcal{C}_p(\bm{d};W_I)\}$ are represented by neural networks with parameters denoted by $W_I$ (where $\mathcal{C}_p(\bm{d};W_I)$ is a symmetric, positive definite, and trace class operator). The neural networks with parameters $W_S$ and $W_I$ are referred to as SNet (sigma network) and INet (inverting network), respectively.

Let us provide intuitive ideas of how to incorporate physical models into the INet, which rely on some analysis of the posterior probability measures given by (\ref{CpostU}) and (\ref{CpostS}). It follows from the formula as shown in Example 6.23 of \citep{Stuart2010AN} that the expression of the mean function (\ref{CposU1}) can be written into the following form:
\begin{align}\label{meanFun1}
\begin{split}
\bar{u}_p = & \bar{u}_p^1 + \bar{u}_p^2 = (H^*\Sigma_{\text{inv}}^*H + \mathcal{C}_0^{-1})^{-1} H^* \Sigma_{\text{inv}}^* \bm{d} +
(H^*\Sigma_{\text{inv}}^*H + \mathcal{C}_0^{-1})^{-1}\mathcal{C}_0^{-1}\bar{u}_0.
\end{split}
\end{align}
Obviously, the posterior mean function contains two parts $\bar{u}_p^1$ and $\bar{u}_p^2$. We analyze these two parts separately and then join them to provide specific structures of the INet.

Now, let us illustrate some intuitive ideas about the data-informed part $\bar{u}_p^1$. Denote $\{e_j\}_{j=1}^{\infty}$ in (\ref{priorDecom}) and let $\{\lambda_j\}_{j=1}^{\infty}$ be a positive sequence belonging to $\ell^2$. In addition, we assume that the forward operator $H$ can be decomposed into $S\circ F$, where
\begin{align}\label{ForE}
(Fu)(x) = \sum_{j=1}^{\infty}\lambda_j\langle u, e_j\rangle e_j(x) = \sum_{j=1}^{\infty}\lambda_j u_je_j(x)  \quad \forall \,\, u\in\mathcal{H}_u
\end{align}
and
\begin{align}\label{ForS}
(S\circ F) u = ((Fu)(x_1), \ldots, (Fu)(x_{N_d}))^T \quad \forall \,\, (x_1,\ldots, x_{N_d})\in\mathbb{R}^{N_d}.
\end{align} This is only needed in Theorem \ref{vanishNoiseTheo}.

\begin{theorem}\label{vanishNoiseTheo}
	Assume that the forward operator $H$ can be decomposed into $S\circ F$, where $F$ and $S$ are defined in (\ref{ForE}) and (\ref{ForS}), respectively.
	Then, we have
	\begin{align}
	\lim_{\max_{1\leq i\leq N_d}\beta_i^*\rightarrow 0}\|(H^*\Sigma_{\text{{\rm inv}}}^*H + \mathcal{C}_0^{-1})^{-1}\mathcal{C}_0^{-1}\bar{u}_0\|_{\mathcal{H}_u} = 0,
	\end{align}
	where
	\begin{align}
	\Sigma_{\text{{\rm inv}}}^{*} = \text{diag}\left( \frac{\alpha_1^*}{\beta_1^*}, \ldots, \frac{\alpha_{N_d}^*}{\beta_{N_d}^*} \right)
	\end{align}
	with $\alpha_i^*, \beta_i^*$ ($i=1,\ldots,N_d$) being positive real numbers.
\end{theorem}

The proof is given in Subsection \ref{AppendixProof} of Appendix. It can be seen from  Theorem \ref{vanishNoiseTheo} that the synthetic truth (prior mean) related term $\bar{u}_p^2$ will disappear as the noise goes to zero, i.e., the data $\bm{d}$ becomes to be the clean data. Hence, in order to fully extract the information of noisy data,
we introduce a neural network named as DNet (denoising network) that aims to transform the noisy data into nearly clean data (data with low noise). The term $\bar{u}_p^1$ will be transformed into 
\begin{align}
	\bar{u}_p^1 = (H^*\Sigma^*_{\text{inv}}H + \mathcal{C}_0)^{-1}H^*\Sigma^*_{\text{inv}}\text{DNet}(\bm{d}).
\end{align}
Since $(H^*\Sigma^*_{\text{inv}}H + \mathcal{C}_0)^{-1}H^*\Sigma^*_{\text{inv}}$ can be seen as a classical inverse method, we can introduce a \textbf{c}omputationally \textbf{e}fficient \textbf{c}lassical \textbf{inv}erse method (which may provide a quick but noisy estimate) combined with a deep neural network called ENet (enhancing neural network). For simplicity, we call this the CECInv layer, which incorporates the physical model into our inversion method. 
In summary, the network structure is defined as follows: 
\begin{align}\label{up1}
	\bar{u}_p^1 = \text{ENet}(\text{CECInv}(\text{DNet}(\bm{d}))).
\end{align}
Practically, the noisy data can hardly be transformed into clean data due to the limited training datasets.
Considering the ill-posedness of the inverse problems, the clean data can only provide limited information to yield estimates with limited accuracy. Therefore, the synthetic truth (prior mean) related term $\bar{u}_p^2$ should also be taken into consideration.

\begin{figure}[t]
	\centering
	\includegraphics[width=0.8\textwidth]{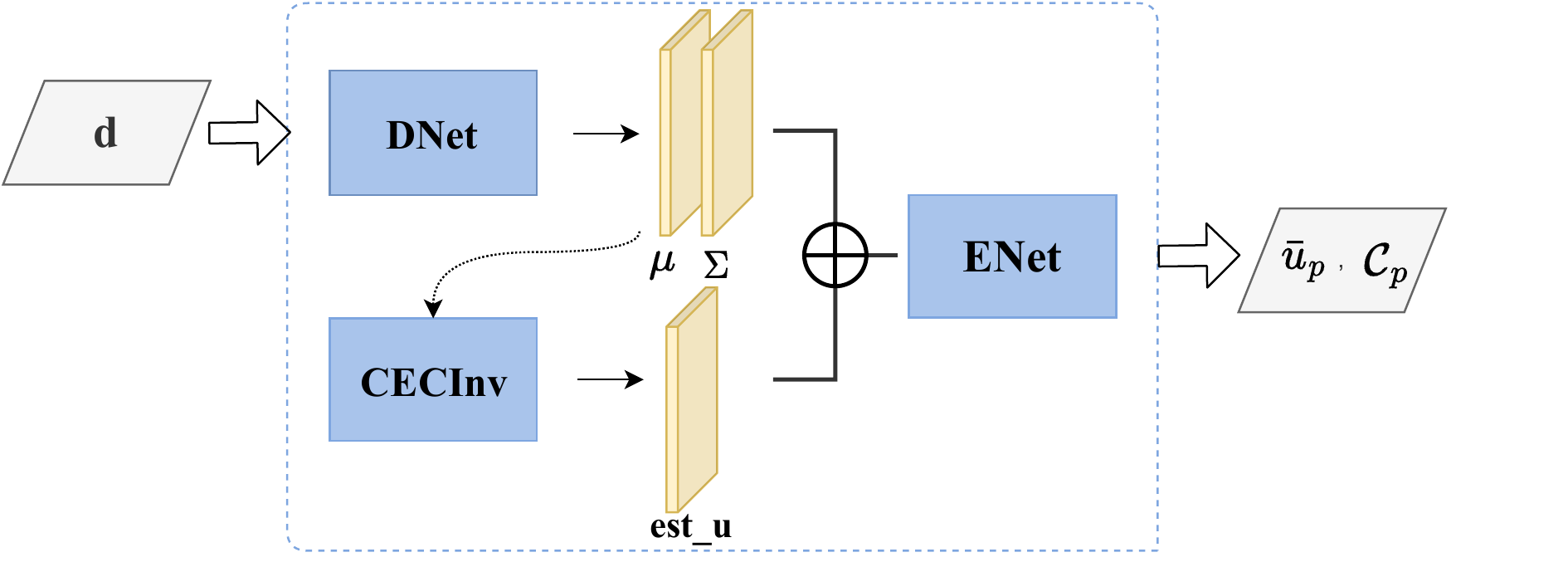}\\
	\vspace*{-0.6em}
	\caption{General structure of INet.}\label{INet}
	\vspace*{0.5em}
\end{figure}

During the training stage, we have the synthetic truth, i.e., the prior mean function $\bar{u}_0$. However,
the synthetic truth is not available during the inference stage, and only the noisy data is available. Based on these 
considerations, we need to design a mechanism that can encode the prior information of the synthetic truth
into the network parameters through the training process on these data with ground truth knowledge. 
Remembering the form of the term $\bar{u}_p^2$, we see that if the synthetic truth is much similar to the truth, then 
the prior variance should be very small. That is to say, we have
\begin{align}\label{up2}
	\bar{u}_p^2 \approx \bar{u}_0. 
\end{align}
Remembering the formula (\ref{up1}), during the training stage, the noisy data is transformed into near-clean data. Then, $\text{ENet}(\text{CECInv}(\cdot))$ provides an estimate that should be as similar as possible to the synthetic truth parameters in the training datasets. Obviously, we may recognize
\begin{align}\label{u0prior}
	\bar{u}_0 \approx \text{ENet}(\text{CECInv}(\text{DNet}(\bm{d}))).
\end{align} 
Finally, considering (\ref{meanFun1}), (\ref{up1}), (\ref{up2}), and (\ref{u0prior}), 
the network for predicting the posterior mean functions could rationally be designed as follow: 
\begin{align}
	\bar{u}_p = \bar{u}_p^1 + \bar{u}_p^2 \approx \text{ENet}(\text{CECInv}(\text{DNet}(\bm{d}))) + \text{ENet}(\text{CECInv}(\text{DNet}(\bm{d}))),
\end{align}
which can be reduced to 
\begin{align}\label{upfinal}
	\bar{u}_p = \text{ENet}(\text{CECInv}(\text{DNet}(\bm{d}))).
\end{align}
We incorporate the variances by adding an additional channel to the related neural networks to represent the function parameters involved in the covariance operator. From the above statements, we can see that the final structure of INet achieves our aim: the general probability model is specified by the variational inference method, the physical model is encoded in the CECInv layer, and the prior information of the training datasets is encoded in both the ENet and DNet (see formula (\ref{u0prior})). For clarity, we depict a schematic diagram of INet in Figure \ref{INet}. 
The general framework illustrated here requires no restrictions on the specific forms of the DNet, ENet, and SNet. In Section \ref{NumericalExamples}, we provide two numerical examples and the particular DNNs used in the present work.

After introducing the general inverse framework, named VINet (variational inverting network), we need to focus on designing algorithms that can efficiently train the involved DNNs. To formulate an optimization problem for the network parameters of INet and SNet, we need to find the variational lower bound under the infinite-dimensional setting as the loss functional to guide the DNN training. This requires the following basic assumption.

\begin{assumption}\label{assVINet}
The approximate probability measure $\nu$ with components $\nu^{u}$ and $\nu^{\bm{\sigma}}$ is assumed to be equivalent to the prior probability measure $\mu_0$, as well as the posterior probability measure $\mu^{\bm{d}}$.
\end{assumption}

Regarding the logarithm of the marginal likelihood, we have
\begin{align*}
\log Z_{\bm{d}} & =  \int_{\mathcal{H}_u \times (\mathbb{R}^{+})^{N_d}} \log Z_{\bm{d}} \, \nu(du,d\bm{\sigma}) \\
& =  \int_{\mathcal{H}_u \times (\mathbb{R}^{+})^{N_d}} \log\left\{
\frac{d\mu_0}{d\nu}\frac{d\nu}{d\mu^{\bm{d}}}\frac{1}{\sqrt{\det(\Sigma)}} \exp\left( -\frac{1}{2}\|\Sigma^{-1/2}(\bm{d}-Hu)\|^2 \right)
\right\}\nu(du,d\bm{\sigma}).
\end{align*}
Expanding the terms in the integrand of the above equality yields 
\begin{align}\label{exe1}
\log Z_{\bm{d}} = \mathcal{L}(u,\bm{\sigma};\bm{d}) + D_{\text{KL}}(\nu||\mu^{\bm{d}}),
\end{align}
where
\begin{align}
\mathcal{L}(u,\bm{\sigma};\bm{d}) = \int \log\frac{d\mu_0}{d\nu} - \frac{1}{2}\log(\det(\Sigma))
- \frac{1}{2}\|\Sigma^{-1/2}(\bm{d}-Hu)\|^2\nu(du,d\bm{\sigma}).
\end{align}
Noting that the KL divergence between the variational approximate measure $\nu$ and the true posterior measure $\mu^{\bm{d}}$ is non-negative, 
we can thus get
\begin{align}
\log Z_{\bm{d}} \geq \mathcal{L}(u,\bm{\sigma};\bm{d}).
\end{align}
Since $\mathcal{L}(u,\bm{\sigma};\bm{d})$ constitutes a lower bound of $\log Z_{\bm{d}}$, we call it evidence lower bound (ELBO)
as in the finite-dimensional setting \citep{Bishop2006book}.
Under Assumption \ref{assVINet}, we know that the optimization problems
\begin{align}\label{exe2}
\min_{W_I, W_S}D_{\text{KL}}(\nu || \mu^{\bm{d}})
\end{align}
and
\begin{align}\label{exe3}
\min_{W_I, W_D}-\mathcal{L}(u,\bm{\sigma};\bm{d})
\end{align}
are equivalent to each other. The ELBO can be decomposed into the following three terms:
\begin{align}\label{FinalEstI123}
\mathcal{L}(u,\bm{\sigma};\bm{d}) = I_1 + I_2 + I_3, 
\end{align}
where
\begin{align}
I_1 = \int\frac{d\mu_0^u}{d\nu^u}\nu^u(du), \quad I_2 = \int\frac{d\mu_0^{\bm{\sigma}}}{d\nu^{\bm{\sigma}}}\nu^{\bm{\sigma}}(d\bm{\sigma}),
\end{align}
\begin{align}
I_3 = \int - \frac{1}{2}\log(\det(\Sigma)) - \frac{1}{2}\|\Sigma^{-1/2}(\bm{d}-Hu)\|^2\nu(du,d\bm{\sigma}).
\end{align}

\textbf{Term $I_1$}.
Different from the finite-dimensional case, it is difficult to reduce the form of $I_1$. However, for some specific problems, the following reparameterization trick can be used to approximate it by Monte Carlo estimation:
\begin{align}
I_1 \approx \frac{d\mu_0^u}{d\nu^u}(\tilde{u}),
\end{align}
where $\tilde{u} = \bar{u}_p(\bm{d};W_I) + \mathcal{C}_p(\bm{d};W_I)^{1/2}\eta$ with $\eta\sim\mathcal{N}(0,\text{Id})$.

\textbf{Term $I_2$}. Since $\bm{\sigma}$ belongs to the finite-dimensional space, we may follow \citep{Yue2019NIPS} and obtain
\begin{align}
\begin{split}
& \int\log\frac{d\nu^{\bm{\sigma}}}{d\mu_0^{\bm{\sigma}}}\nu^{\bm{\sigma}}(d\bm{\sigma}) =
\sum_{i=1}^{N_d}\Bigg\{
\left( \alpha_i - \alpha_i^0 + 1 \right)\psi(\alpha_i) + \alpha_i\left( \frac{\beta_i^0}{\beta_i}-1 \right)  \\
& \qquad\quad
+ \left( \alpha_i^0-1 \right)\left( \log\beta_i - \log\beta_i^0 \right) +
\left[ \log\Gamma\left( \alpha_i^0 - 1 \right) - \log\Gamma(\alpha_i) \right]
\Bigg\},
\end{split}
\end{align}
where $\Gamma(\cdot)$ and $\psi(\cdot)$ represent the Gamma and Degamma functions, respectively.

\textbf{Term $I_3$}. For the first term in $I_3$, we have
\begin{align}\label{termI3_1}
\begin{split}
-\frac{1}{2}\int \log(\det(\Sigma))\nu(du,d\bm{\sigma}) &=  -\frac{1}{2}\sum_{i=1}^{N_d}\int \log\sigma_i^2 \nu^{\bm{\sigma}}(d\bm{\sigma}) \\
& =  -\frac{1}{2}\sum_{i=1}^{N_d}\left( \log\beta_i - \psi(\alpha_i) \right).
\end{split}
\end{align}
For the second term in $I_3$, a simple calculation yields 
\begin{align}\label{termI3_2}
\begin{split}
-\frac{1}{2}\int \|(\bm{d}-Hu)\|_{\Sigma}^2 \nu(du,d\bm{\sigma}) = &
-\frac{1}{2}\int (\bm{d}-Hu)^T \Sigma_{\text{inv}} (\bm{d}-Hu)\nu^u(du),
\end{split}
\end{align}
where
\begin{align}
\begin{split}
\Sigma_{\text{inv}} = \text{diag}\left( \mathbb{E}^{\nu^{\bm{\sigma}}}[\sigma_1^{-1}], \ldots, \mathbb{E}^{\nu^{\bm{\sigma}}}[\sigma_{N_d}^{-1}] \right)
= \text{diag}\left( \frac{\alpha_1}{\beta_1},\ldots,\frac{\alpha_{N_d}}{\beta_{N_d}} \right).
\end{split}
\end{align}
Using (\ref{termI3_1})--(\ref{termI3_2}) and the Monte Carlo estimate, we get
\begin{align}
I_3 \approx -\frac{1}{2}\sum_{i=1}^{N_d}\left( \log\beta_i - \psi(\alpha_i) \right) -
\frac{1}{2}(\bm{d}-H\tilde{u})^T\Sigma_{\text{inv}}(\bm{d}-H\tilde{u}),
\end{align}
where $\tilde{u} = \bar{u}_p(\bm{d};W_I) + \mathcal{C}_p(\bm{d};W_I)^{1/2}\eta$ with $\eta\sim\mathcal{N}(0,\text{Id})$.

For a given training data set $\{\bm{d}_i, u_i\}_{i=1}^{N_e}$, where $N_e$ is a positive integer, the following objective function, the negative ELBO on the entire training set, can be constructed to optimize the parameters of INet and SNet: 
\begin{align}\label{FinalOpProb}
\min_{W_I,W_S} - \sum_{i=1}^{N_e}\mathcal{L}(u^{(i)}, \bm{\sigma}^{(i)}; \bm{d}^{(i)}),
\end{align}
where $u^{(i)}$, $\bm{\sigma}^{(i)}$, and $\bm{d}^{(i)}$ denote the corresponding variables for the $i$-th training pair.

\begin{remark}
Here, we provide some intuitive explanations of the VINet. From formula (\ref{exe1}), we know that the problems (\ref{exe2}) and (\ref{exe3}) are equivalent. If there is only one piece of training data, the optimization problem (\ref{exe3}) will find the parameters $W_I$ and $W_S$ of DNNs that make the parameterized measure (\ref{postPM}) similar to the approximate posterior measure calculated from the mean-field assumption based IDVI method. We can expect that the learned model from one piece of training data will have poor generalization properties. During the training stage, we actually have $N_e$ pieces of training data and solve the minimization problem (\ref{FinalOpProb}) that can be written as 
\begin{align}\label{FinalOpProb11}
	\min_{W_I,W_S} - \frac{1}{N_e}\sum_{i=1}^{N_e}\mathcal{L}(u^{(i)}, \bm{\sigma}^{(i)}; \bm{d}^{(i)}).
\end{align}
Hence, our training process intuitively finds common network parameters $W_I$ and $W_S$ that perform well on all of the training datasets on average. In this sense, we can intuitively recognize that the DNet, ENet, and SNet capture the information shared by all training datasets (prior and noise information) in some average sense. During the inference stage, there will be no synthetic truth. The trained VINet will combine the information extracted from the new data $\bm{d}$, the physical model encoded in the CECInv layer, and the prior information encoded in the DNet, SNet, and ENet to yield estimates of the function- and vector-valued parameters of the probability measure (\ref{postPM}).
\end{remark}

At the end of this subsection, we present a preliminary analysis of the estimation error of the mean function, which provides some intuitions for future study. The proof is given in Subsection \ref{AppendixProof} of the Appendix.

\begin{theorem}\label{errorSimpleAnalysis}
	Let $u^{\dag}$ be the background truth function,
	$\bar{u}_p$ be the estimated function given by (\ref{meanFun1}),
	and $\mathbb{E}_0$ be the expectation corresponding to the distribution of the data.
	Then the following estimate holds: 
	\begin{align}\label{simpleErrorEst}
	\begin{split}
	\mathbb{E}_0[\|\bar{u}_p(\bm{d};W_I) - u^{\dag}\|_{\mathcal{H}_u}] & \leq 
	\mathbb{E}_0[\|\bar{u}_p(\bm{d};W_I) - \bar{u}_p\|_{\mathcal{H}_u}] +
	\|\bar{u}_0 - u^{\dag}\|_{\mathcal{H}_u}  \\
	&\quad   + \text{tr}\left( \mathcal{C}_pH^*\Sigma_{\text{{\rm inv}}}^*\Sigma \Sigma_{\text{{\rm inv}}}^*H \mathcal{C}_p^* \right)^{1/2},
	\end{split}
	\end{align}
	where $\Sigma$ is the background truth noise covariance matrix and
	$\Sigma_{\text{{\rm inv}}}^*$ is the diagonal matrix
	$\text{diag}\left( \alpha_1^*/\beta_1^*, \ldots, \alpha_{N_d}^*/\beta_{N_d}^* \right)$ defined in (\ref{221}).
\end{theorem}

We have the following observations on the three terms appearing in (\ref{simpleErrorEst}). The first term on the right-hand side of (\ref{simpleErrorEst}) reflects the approximation error arising from learning by DNNs, which can be analyzed in detail, e.g., by introducing the concepts of sample error and approximation error \citep{Cucker2001BAMS}. The DNNs employed here need to learn a nonlinear operator between two separable Hilbert spaces, which, to the best of our knowledge, can hardly be analyzed under the current learning theory. However, there are some works on the learnability of the nonlinear operators \citep{Chen1995IEEE} that should be useful for future investigations.

The second term accounts for the error induced by the synthetic function parameter used for training and the true background function. This term is difficult to analyze in its current form since $u^{\dag}$ and $\bar{u}_0$ cannot be obtained in the inverse computational stage. However, we can assume that the residual of the synthetic function and the background truth can be bounded by some constant with high probability. The third term on the right-hand side of (\ref{simpleErrorEst}) represents the errors induced by noise and the prior covariance operator.

\subsection{Parametric strategies}\label{ParametizationVINet}

In the implementation of the algorithm, it is necessary to choose an appropriate prior probability measure $\mu_0$ and approximate probability measure $\nu$ to ensure Assumption \ref{assVINet}. This subsection focuses on a possible parametric strategy.

For the finite-dimensional parameter $\bm{\sigma}$, there are no singular issues with the prior and posterior probability measures. Therefore, we focus on the more problematic function parameter $u$. The major difficulty is that Gaussian measures defined on infinite-dimensional space are singular with each other when the covariance operators are scaled by different constants. It remains an open question as to how to resolve such a singularity problem in terms of theory and algorithms \citep{agapiou2014SIAM,Dunlop2017SC}.

Let us introduce some assumptions proposed in \citep{Dashti2014} that inspire us to define appropriate covariance operators in the prior and posterior probability measures.

\begin{assumption}\label{Ass_A}\citep{Dashti2014}
	The operator $A$, densely defined on the Hilbert space $\mathcal{H}=L^2(\Omega;\mathbb{R})$ where 
	$\Omega\subset\mathbb{R}^d$ is a bounded domain with smooth boundary, is assumed to satisfy the following properties:
	\begin{enumerate}
		\item $A$ is positive definite, self-adjoint and invertible;
		\item The eigenfunctions $\{e_j\}_{j\in\mathbb{N}}$ of $A$ form an orthonormal basis for $\mathcal{H}$;
		\item The eigenvalues of $A$ satisfy $\lambda_j \asymp j^{2/d}$;
		\item There is $C>0$ such that
		\begin{align*}
		\sup_{j\in\mathbb{N}}\left( \|e_j\|_{L^{\infty}} + \frac{1}{j^{1/d}}\text{Lip}(e_j) \right) \leq C,
		\end{align*}
		where $\text{Lip}(e_j)$ denotes the Lipschitz constant of $e_j, j=1,\ldots,\infty$.
	\end{enumerate}
\end{assumption}

With Assumption \ref{Ass_A}, we now restrict the prior probability measure of $u$ to be $\mu_0^u = \mathcal{N}(\bar{u}_0, \mathcal{C}_{\epsilon_0})$ with $\mathcal{C}_{\epsilon_0} = (\epsilon_0^{-1}\text{Id} + \delta A^{\alpha/2})^{-2}$. Here, the parameters $\epsilon_0,\delta >0$ are positive constants satisfying $\delta \ll \epsilon_0^{-1}$, and the parameter $\alpha > d/2$ controls the regularity properties of the generated random samples. The mean function $\bar{u}_0$ is taken as the simulated background truth, which is available at the training stage but cannot be obtained at the inference stage.

\begin{remark}\label{remarkunknonTruth}
When solving inverse problems, it is usually difficult to obtain absolutely correct background true parameters as a training dataset. For example, in geophysical applications, the ``Marmousi'' model, as a typical synthetic model \citep{Haber2003IP}, is frequently used for training. Therefore, we may choose the simulated background truth as the prior mean function, and the covariance operator can account for the uncertainty information due to the inaccuracy of the training dataset.
\end{remark}

For the posterior probability measure $\nu^{u} = \mathcal{N}(\bar{u}_p(\bm{d};W_I), \mathcal{C}_p(\bm{d};W_I))$,
we specify the covariance operator as follows:
\begin{align}
\mathcal{C}_p(\bm{d}; W_I) = (a(\bm{d};W_I)\text{Id} + \delta A^{\alpha/2})^{-2},
\end{align}
where $\alpha$ and $\delta$ are taken the same as the prior measure, and $a(\bm{d};W_I)$ is a positive function parameterized by the neural network with parameters $W_I$. In the following, we let $\mathcal{H}_u = \mathcal{H}$ which is defined in Assumption \ref{Ass_A}.

The following result illustrates that $\mu_0^{u}$ and $\nu^{u}$ defined above are equivalent under certain general conditions.

\begin{theorem}\label{equivalentTheorem}
	Assume that $\Omega\subset\mathbb{R}^{d}$ ($d\leq 3$) is a bounded domain with a smooth boundary, and the function $a(\bm{d};W_I)\in L^{\infty}(\Omega)$ has positive lower and upper bounds, i.e., $0 < \underline{a}\leq a(\bm{d};W_I)(x) \leq \bar{a}$ for all $x\in\Omega$ with $\underline{a},\bar{a}\in\mathbb{R}^{+}$. Assume further that the mean functions $\bar{u}$ and $\bar{u}_p(\bm{d};W_I)$ introduced in $\mu_0^u$ and $\nu^{u}$ belong to the Hilbert scale $\mathcal{C}_{\epsilon_0}^{1/2}\mathcal{H}_u$. Then, the approximate posterior measure $\nu^u$ is equivalent to the prior measure $\mu_0^u$.
\end{theorem}

The proof is provided in Subsection \ref{AppendixProof} of the Appendix.

\section{Numerical Examples}\label{NumericalExamples}

In this section, we provide the specific settings of INet and SNet for inverse problems of PDEs. We then apply the VINet to two examples: a simple smoothing model, which is widely used for benchmark testing \citep{agapiou2014SIAM}, and an inverse source problem of the Helmholtz equation.

\subsection{Network structures}\label{NetworkStructure}

In Section \ref{generalTheory}, we introduce a general inverse framework called VINet, wherein any type of rational neural network can be chosen. Well-constructed neural networks typically perform better and require a smaller training dataset, which is crucial for practical applications.

The INet proposed in the general inference framework has three functionalities:
\begin{enumerate}
	\item The network should be able to accurately transmit the information of the noises $\bm{\epsilon}$,
	which is assumed to be distributed according to some non-i.i.d. Gaussian distribution.
	\item The network should be able to learn the prior information from the training data set, which can be used to provide high-quality estimations.
	\item The forward evaluation of the network can be calculated more efficiently than the classical iterative algorithm.
\end{enumerate}

\begin{figure}[t]
	\centering
	\includegraphics[width=0.7\textwidth]{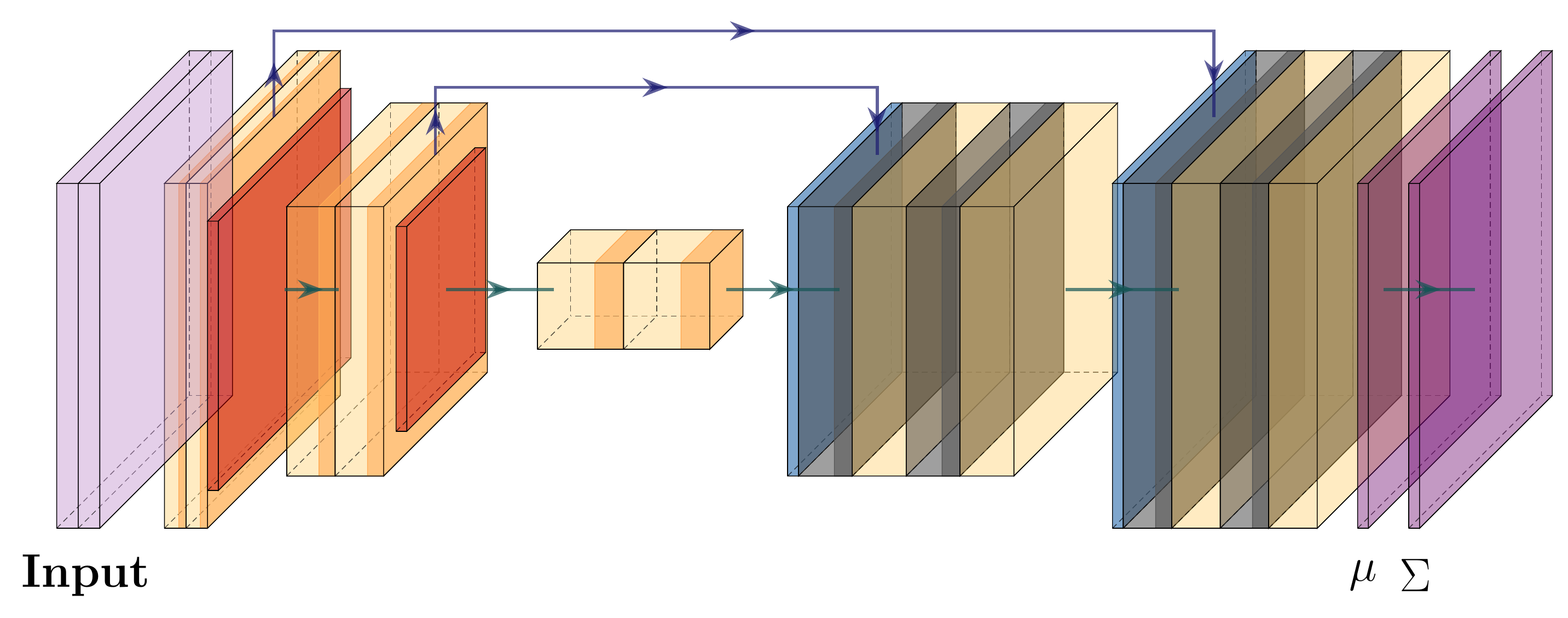}\\
	\vspace*{-0.6em}
	\caption{The tiny U-Net architecture employed for constructing DNet.}\label{UNet}
	\vspace*{0.5em}
\end{figure}

Based on the three criteria above, we designed a specific INet structure in the previous section, which is illustrated in Figure \ref{INet}. The network contains three components. The first component, named DNet, undertakes the denoising task, which is expected to learn the noise structure and provide estimated clean data. The second component, named CECInv, is a small-scale classical inversion, in which every inverse of sparse matrices can be explicitly calculated, and the results are some dense matrices. Clearly, only matrix multiplications are needed during the computations, which can be more efficiently implemented than solving linear equations. The third component, named ENet, can be viewed as an enhancing layer, which learns the residual of the inexact inversion obtained by the CECInv layer with the background truth. More details of the three components are listed below:

\begin{enumerate}
	\item DNet is chosen to be a small U-Net architecture, as illustrated in Figure \ref{UNet}, containing three encoder blocks, two decoder blocks, and symmetric skip connections under each scale.
	\item A small-scale classical inversion (which may be different for each problem) is employed as the CECInv layer, providing a low-resolution estimate (which may contain some artifacts).
	\item For the ENet, we adopt two types of neural networks: U-Net with rescale layers (details are given in Remark \ref{DesNetRemark}) and Fourier neural network (FNO) \citep{Zongyi2020ICLR}. We will provide specific settings later, since different parameters are employed for the two numerical examples.
\end{enumerate}

As for the SNet, it is adopted to infer the variational posterior parameters $\{\alpha_i(\bm{d};W_S)\}_{i=1}^{\infty}$, $\{\beta_i(\bm{d};W_S)\}_{i=1}^{\infty}$ from the noisy measurement data $\bm{d}$. Similar to the finite-dimensional denoising case \citep{Yue2019NIPS}, we use the DnCNN \citep{Zhang2017IEEE} architecture with five layers, and the feature channels of each layer are set to be 64.

\begin{figure}[t]
	\centering
	\includegraphics[width=1.0\textwidth]{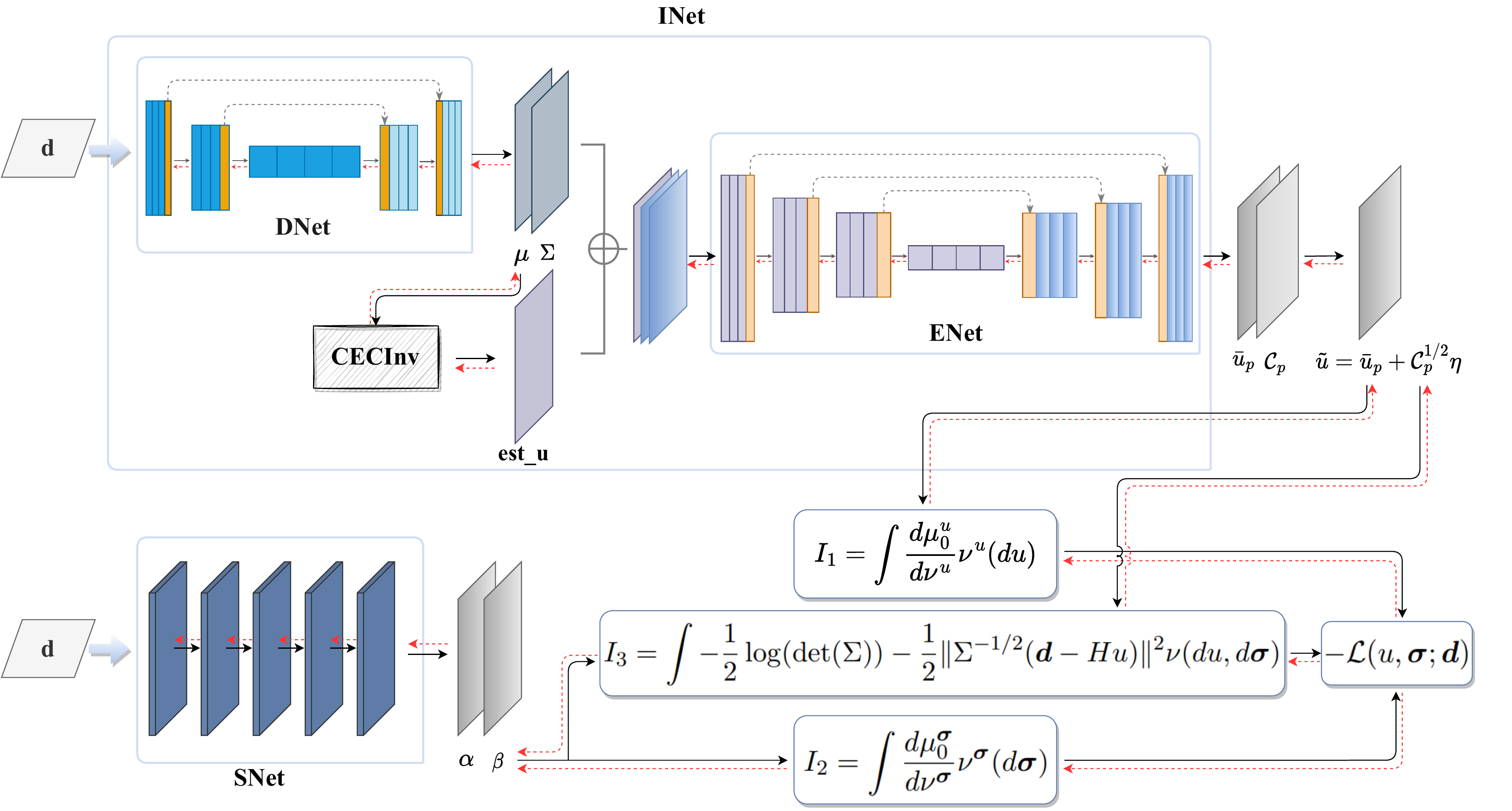}\\
	\vspace*{-0.6em}
	\caption{A particular example of the VINet framework (ENet is illustrated by a U-Net architecture).}\label{FullNet}
	\vspace*{0.5em}
\end{figure}

Combining INet and SNet, we illustrate the whole framework in Figure \ref{FullNet} with ENet has
chosen to be the U-Net that contains rescale layers.
It is easy to train our model based on the objective function presented in (\ref{FinalOpProb}) by using the stochastic gradient descent algorithm or any variant optimization techniques. As for the finite-dimensional case, each term of (\ref{FinalEstI123}) can be intuitively explained. The INet and SNet are updated from the terms $I_1$ and $I_2$, respectively, by controlling the discrepancy between the variational posteriors and the priors. Corresponding to the noise in the training dataset, the term $I_3$ couples INet and SNet together, generating gradients simultaneously during backpropagation.

For any new noisy data $\bm{d}$, the trained INet can be used to obtain the final estimated result directly and explicitly. Additionally, SNet can be used to easily infer the noise information, even with complicated non-i.i.d. configurations, by inputting the noisy data.

Numerically, the related PDEs are solved using the finite element method, which is implemented with the open software FEniCS (Version 2019.1.0). Additional information on FEniCS can be found in \citep{Logg2012Book}. The related neural networks are implemented with the open software PyTorch (Version 1.12.1+cu115). The mesh of the finite element method is generated automatically by the FEniCS command \emph{UnitSquareMesh($\cdot$)}. Since the mesh is regular, the finite difference discretization can also be adopted based on the same grid points. Therefore, we use the finite difference method to compute the integral formulas involved in Terms 1 and 3 of the loss functional (\ref{FinalEstI123}). All programs were run on a computer with an Intel(R) Core(TM) i9-12900KF (CPU), a GeForce RTX 3080 (GPU) with 12 GB memory, and Ubuntu 20.04.1 LTS (OS).

\begin{remark}\label{DesNetRemark}
It should be mentioned that the ENet is a mapping between two infinite-dimensional spaces, which is typically specified as a neural network with discretization-invariant property, e.g., the neural network proposed in \citep{Bhattacharya2021SMAI-JCM,Nelsen2021SISC,Zongyi2020arXiv,Zongyi2020ICLR}. In the present work, we employed the FNO \citep{Zongyi2020ICLR} as one possible choice of the ENet, which learns a fixed number of Fourier modes for each layer, making the whole network process discretization-invariant. Besides the FNO, we also used the U-Net as one option for the ENet. As is well-known, the U-Net structure does not have discretization-invariant property. To overcome this obstacle, we added rescale layers as the input and output layers. The core idea of the rescale layer is that the input function is viewed as being discretized on a mesh by the finite element method, so we can compute the function values on any different discretized meshes. With rescale layers, we projected the input functions on a fixed mesh (e.g., a mesh generated by UnitSquareMesh($\cdot$) command in FEniCS with mesh size $150\times150$) and then projected the output function on a desired mesh by a rescale layer. Numerical examples show that the U-Net with rescale layers processes the discretization-invariant property to some extent.
\end{remark}

\subsection{A simple smoothing model}\label{SimpleModel}

\subsubsection{Basic settings}\label{SimpleModel1}

Consider an inverse source problem of the elliptic equation 
\begin{align}\label{simpleSmoothModel}
\begin{split}
-\alpha\Delta w + w & = u  \quad \text{in }\Omega, \\
\frac{\partial w}{\partial \bm{n}} & = 0 \quad \text{on }\partial\Omega,
\end{split}
\end{align}
where $\Omega=[0,1]^2\subset\mathbb{R}^2$, $\alpha>0$ is a positive constant, and $\bm{n}$ denotes the outward normal vector.
The forward operator $H$ is defined as follows:
\begin{align}
Hu = (w(\bm{x}_1), w(\bm{x}_2), \cdots, w(\bm{x}_{N_d}))^T,
\end{align}
where $u\in \mathcal{H}_u := L^{2}(\Omega)$, $w$ denotes the solution of (\ref{simpleSmoothModel}), and
$\bm{x}_i\in\Omega$ for $i=1,\ldots,N_d$. In our implementations, the measurement points $\{\bm{x}_i\}_{i=1}^{N_d}$ are taken at the coordinates
$\left\{ \left(\frac{k}{20}, \frac{\ell}{20}\right) \right\}_{k,\ell=1}^{20}$.
To avoid the inverse crime \citep{Kaipio2004Book}, we discretize the elliptic equation by the finite element method on a regular
mesh (the grid points are uniformly distributed on the rectangular domain) with the number of grid points being equal to $400\times 400$. All of the training and testing data pairs are generated based on this discretization.
For revealing the impact of the discretization resolution, 
the related functions are discretized by finite element method on a series of regular mesh
with the number of grid points $n=\{30\times30, 50\times50, 100\times100, 150\times150, 200\times200, 250\times250, 300\times300\}$
at the inverse stage.

For the operator $A$ in Assumption \ref{Ass_A} used in the prior probability measure, it is chosen to be the negative Laplace operator with Neumann boundary condition, which is similar to \citep{Tan2013IP}. Hence, the operators in Subsection \ref{ParametizationVINet} can be specified explicitly as follows:
\begin{align*}
\mathcal{C}_{\epsilon_0} = (\epsilon_0^{-1}\text{Id}+\delta(-\Delta)^{\alpha/2})^{-2}, \quad 
\mathcal{C}_{p}(\bm{d};W_I) = (a(\bm{d};W_I)\text{Id}+\delta(-\Delta)^{\alpha/2})^{-2}.
\end{align*}
To simplify the numerical implementations, we choose $\alpha=2$, which makes the operators $\mathcal{C}_{\epsilon_0}$ and $\mathcal{C}_{p}(\bm{d};W_I)$ to be discretized as illustrated in \citep{Tan2013IP}. Since the fractional Laplace operator needs to be discretized differently \citep{Thanh2016IPI,Lischke2020JCP}, we leave the case $\alpha\neq 2$ for our future research.

Concerned with the CECInv layer, we employed a classical truncated singular value method \citep{Engl1996book}. Specifically, the inverse method can be written in the form
\begin{align}\label{SVDDef}
R_{\lambda_{m}}\hat{\bm{d}} = q_{\lambda_{m}}(H^*H)H^*\hat{\bm{d}},
\end{align}
where $\hat{\bm{d}}$ is the output of the DNet, and
\begin{align}\label{SVDdef}
q_{\lambda_m}(\lambda) = \left\{\begin{aligned}
& \lambda^{-1}, \quad \lambda \geq \lambda_m, \\
& 0, \qquad\, \lambda < \lambda_m,
\end{aligned}\right.
\end{align}
with $\lambda_m$ denoting a predefined truncate level. In the above, the operator $q_{\lambda_{m}}(H^*H)$ appearing in (\ref{SVDDef}) should be understood as the functional calculus for bounded self-adjoint operators, as described in Chapter 2 of \citep{Engl1996book}. The PDEs were discretized on a regular mesh with the number of grid points being equal to $30\times30$, which is relatively small compared with the forward discretization. Based on such a small mesh, we can calculate the related inversion of the sparse matrices generated by the finite element method explicitly, i.e., there is no need to solve a lot of sparse linear equations.

In our experiments, we generate $1000$ training function parameters $u$ based on a Gaussian probability measure with mean function equal to constant $2$ and covariance operator given by 
\begin{align}
	\mathcal{C}^{-1} = \alpha^{2}(\text{Id} - \nabla\cdot(\Theta\nabla \cdot))^2,
\end{align}
where $\Theta := \text{diag}(10,1)$ and $\alpha=0.01$. 
For testing the proposed methods, we generate two types of testing datasets with $100$ function parameters $u$ for each one. Specifically, type 1 test datasets are generated by a Gaussian probability measure with the same settings as the training datasets. Type 2 test datasets are generated by a Gaussian probability measure with the mean function given by  $2 + 0.2\sin(2\pi x)\sin(2\pi y)$ and the same covariance operator as the training datasets. 
Since, as mentioned in Remark \ref{remarkunknonTruth}, we may not know the ideal background truth practically, type 2 test datasets are designed to test the generalization abilities of the proposed VINet framework.

Previous studies \citep{Yue2019NIPS,Yue2020VariationalIRN} have implied that multiplicative noises can be removed efficiently by a non-i.i.d. noise model. Therefore, multiplicative noises are added to the clean data $\bm{d}_{c}$ to generate noisy data $\bm{d}$ in order to test the performance of the proposed VINet on non-i.i.d. noise cases. In the training stage, we add noises according to $\bm{d} = \bm{d}_{c} + a(\bm{\eta}\odot\bm{d}_{c})$, where $\odot$ is an element-wise product, $\bm{\eta}\sim\mathcal{N}(0,\text{Id})$, and $a\sim U[-0.1,0.1]$. In the testing stage, we take $a=0.1$, i.e., $\bm{d} = \bm{d}_{c} + 0.1(\bm{\eta}\odot\bm{d}_{c})$, where $\bm{\eta}\sim\mathcal{N}(0,\text{Id})$.

\begin{figure}[thbp]
	\centering
	\includegraphics[width=1.0\textwidth]{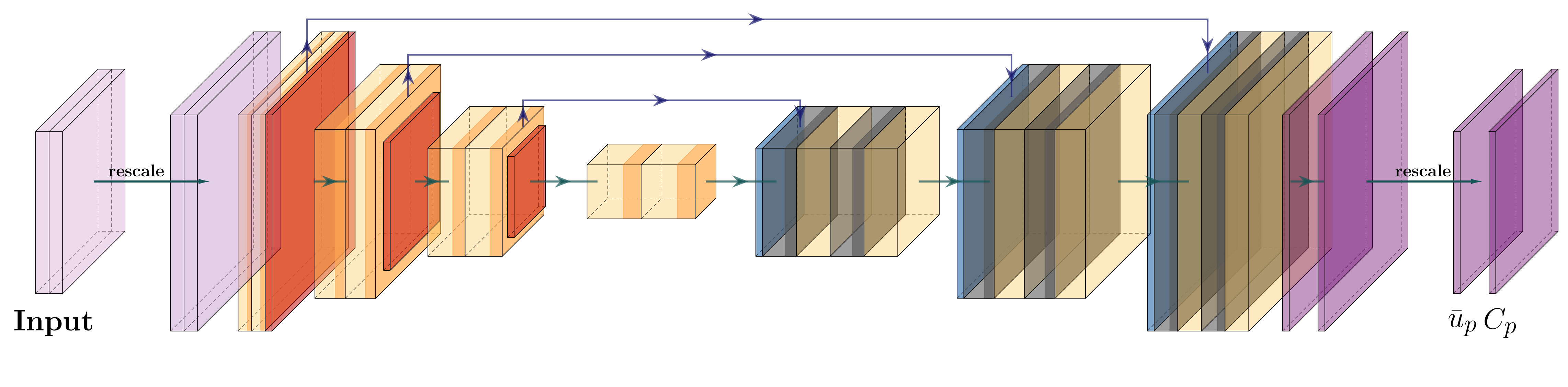}\\
	\vspace*{-0.6em}
	\caption{The U-Net architecture is employed for constructing ENet. The input discretized functions are projected onto a fixed mesh by a rescaling layer. Then, through a conventional U-Net, the output functions are obtained on the fixed mesh. Finally, the output functions are projected onto the required mesh by a rescaling layer. }\label{ENetU}
	\vspace*{0.5em}
\end{figure}

\begin{figure}[t]
	\centering
	\includegraphics[width=1.0\textwidth]{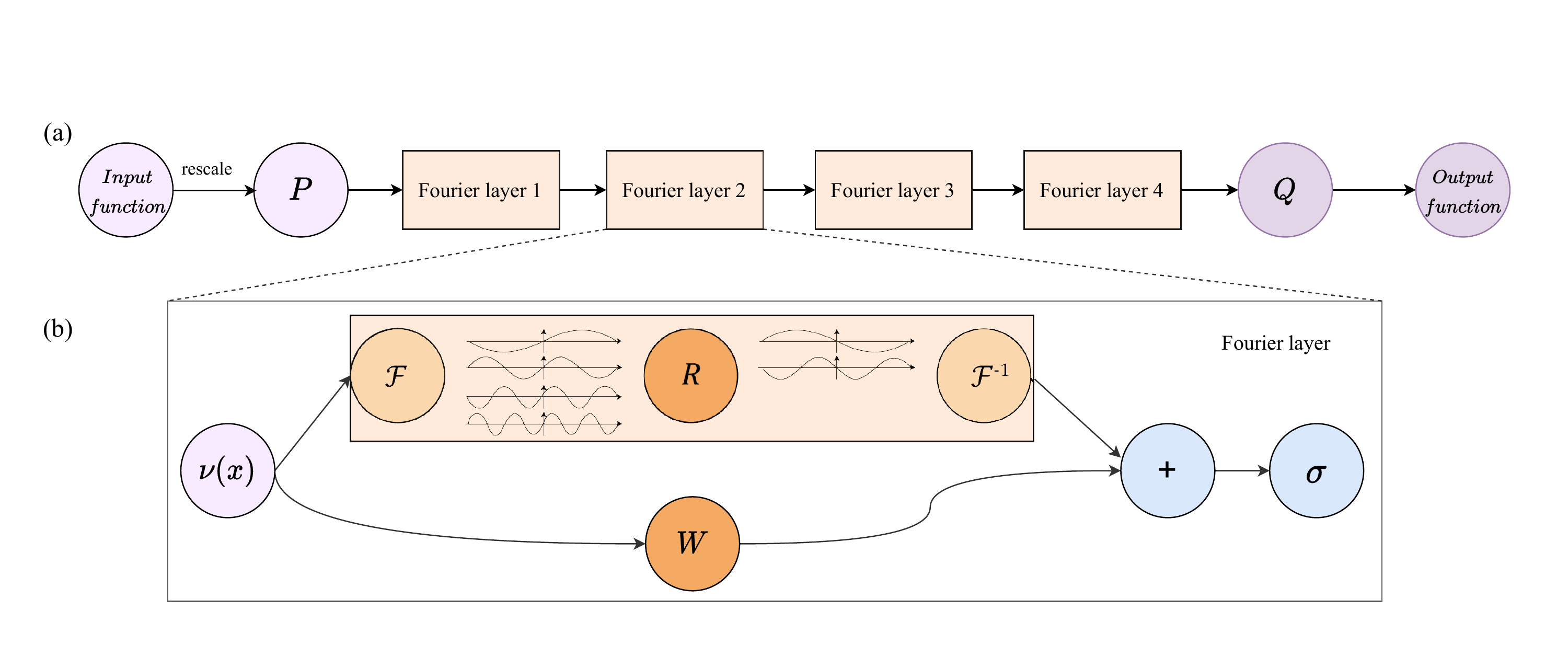}\\
	\vspace*{-0.6em}
	\caption{The FNO architecture is employed for constructing ENet (this figure is nearly the same as Figure 2 in \citep{Zongyi2020ICLR}). (a) \textbf{The full architecture of FNO}: Input functions are lifted to a higher-dimensional channel space by a neural network $P$. Four Fourier layers are employed to enhance the results of the CECInv layer. Finally, the functions are projected back to the target dimensions (2 in our setting) by a neural network $Q$. 
	(b) \textbf{Fourier layers}: Upper part: the Fourier transform $\mathcal{F}$ is applied; a linear transform $R$ is implemented on the lower Fourier modes; then the inverse Fourier transform $\mathcal{F}^{-1}$ is applied. Lower part: a local linear transform $W$ is applied.
	}\label{ENetFNO}
	\vspace*{0.5em}
\end{figure}

On the one hand, the ENet has been chosen to be a U-Net, as illustrated in Figure \ref{ENetU}, which contains four encoder blocks, three decoder blocks, symmetric skip connections under each scale, and two rescale blocks as the input and output layers. For simplicity, we refer to this as U-ENet. On the other hand, we also take ENet to be an FNO, as illustrated in Figure \ref{ENetFNO}, which contains one lift layer, four Fourier integral operator layers, and one project layer. We set the maximal number of Fourier modes to be twelve for each Fourier integral operator, and the original inputs are lifted to a higher dimensional representation with a dimension equal to thirty-two. These settings are the same as the original paper of FNO \citep{Zongyi2020ICLR} for two-dimensional problems. For simplicity, we refer to this as FNO-ENet. Lastly, we should mention that the input layer of ENet contains two channels: the output of the CECInv layer and the variance function obtained from DNet. The output layer of ENet also contains two channels: the mean function and the function parameter in the covariance operator.

For training the VINet, we first train the DNet and SNet, then train the ENet and SNet with the trained DNet (parameters are fixed), which makes the two components more interpretable and reduces the training time. For this example of DNet, we use the Adam optimizer to train for 70 epochs with an initial learning rate of 0.0001, halved every 10 epochs. For the U-ENet and FNO-ENet, we also use the Adam optimizer to train for 15 epochs with an initial learning rate of 0.0001, halved every 10 epochs.

\subsubsection{Numerical results}\label{SimpleModel2}

We compare the results of VINet, the mean-field based variational inference method (MFVI) exhibited in Algorithm \ref{classicalVI}, and the truncated singular value decomposition (TSVD) method shown in formula (\ref{SVDDef}). 
For the MFVI method, the estimated function $u$ highly depends on the particular choices of the prior measure.
Here, we assume that the distribution generating the truth is known, which means the statistical properties of the truth are known exactly for the MFVI method. For comparison, we define the relative error as follows: 
\begin{align}\label{DefRelError}
	\text{relative error} = \frac{\|\tilde{u} - u_{\text{truth}}\|_{L^2}}{\|u_{\text{truth}}\|_{L^2}},
\end{align} 
where $\tilde{u}$ is the estimated function.

For the TSVD method, its performance depends on the predefined truncation level. To determine the optimal truncation level $\lambda_m$, we set the truncation level equal to $2, 1.5, 1, 0.5$, and $0.1$. Then we applied the TSVD method to the type 1 test datasets, which provided the optimal truncation level based on the average relative errors. Detailed results are given in Table \ref{TSVG_est}. From the table, we can see that the optimal truncation level is $1.0$, no matter which discrete levels are employed. Therefore, we always specify the truncation level to be $1.0$ in the following.
In our opinion, the reason the discrete levels are not highly influential is that the measured noisy data $\bm{d}$ can only provide limited information, and the background true parameters are smooth in some sense. For the next numerical example, the discrete level will be more influential due to the poor regularity of the parameter $u$. One point should be mentioned: we employed the discrete strategy illustrated in \citep{Tan2013IP} to calculate the eigen-system, i.e., the operator $H$ will never be constructed (only $Hu$ can be calculated), so a randomized algorithm named double pass \citep{Halko2011SIAMReview,Villa2021ACMTMS} has been employed that does not require $H$ and $H^*$ to be constructed explicitly.

\begin{table}[htbp]
	\caption{Average relative errors of the estimates obtained by the TSVD.}
	\vspace{1mm}
	\begin{center}
		\begin{spacing}{1.5}
			\begin{tabular}{c|ccccc}
				\Xhline{1.1pt}
				Truncate level & 2.0 & 1.5 & 1.0 & 0.5 & 0.1 \\
				\hline
				Relative error ($n=30^2$) & 0.1527 & 0.1523 & \textbf{0.1484} & 0.2020 & 1.1417  \\
				\hline
				Relative error ($n=50^2$) & 0.1484 & 0.1484 & \textbf{0.1439} & 0.1945 & 1.1743	\\
				\hline
				Relative error ($n=100^2$) & 0.1489 & 0.1488 & \textbf{0.1448} & 0.1939 & 1.3359	\\
				\hline 
				Relative error ($n=150^2$) & 0.1493 & 0.1493 & \textbf{0.1455} & 0.1942 & 1.3299	\\
				\hline
				Relative error ($n=200^2$) & 0.1494 & 0.1494 & \textbf{0.1456} & 0.1942 & 1.4759	\\
				\hline 
				Relative error ($n=250^2$) & 0.1495 & 0.1495 & \textbf{0.1482} & 0.1942 & 1.4752	\\
				\hline 
				Relative error ($n=300^2$) & 0.1495 & 0.1495 & \textbf{0.1482} & 0.1944 & 1.4752	\\
				\Xhline{1.1pt}
			\end{tabular}
		\end{spacing}
	\end{center}\label{TSVG_est}
\end{table}

From the iteration procedure of Algorithm \ref{classicalVI}, it is easily seen that the MFVI method is computationally expensive (see Remark \ref{compuBurdenClassVI}). Let us denote $N_{ite}$ as the iterative number of linear equation solvers employed to solve $(H^*\Sigma_{\text{inv}}^{k*}H + \mathcal{C}_0^{-1}) \tilde{u} = u$ mentioned in Remark \ref{compuBurdenClassVI}. Let $k_{max}$ be the maximum iterative number of Algorithm \ref{classicalVI}.
In order to reduce the computational time, we employ the methods given in \citep{Jia2020SIAM}, which calculate the mean function $u_k$ in Algorithm \ref{classicalVI} with a fine mesh and then project the mean function to a rough mesh (mesh size $n=50\times 50$) to reduce the computational complexity for updating $\beta$. Although the computational complexity has been reduced, we still need to solve approximately $(2N_{ite}+2)k_{max}$ number of PDEs, e.g., if $N_{ite}=100$ and $k_{max}=20$, there are $4040$ number of PDEs that need to be solved. 
Actually, the above illustration is just a rough estimate, and we only control the errors for calculating $u_k$ but not the iteration number $N_{ite}$ since the accuracy of $u_k$ has a large impact on the performance of MFVI. Hence the number $N_{ite}$ may be larger than $100$. Considering the large computational complexity, we only apply the MFVI method to the first 50 types 1 and 2 testing examples. Since MFVI is not a learning based method, we can expect that the average performance on the first 50 testing examples is similar to the average performance on the whole testing examples.

Before presenting the numerical results, we should mention that the VINets with U-ENet and FNO-ENet were all trained when the related function parameters $u$ were discretized on a regular mesh with a mesh size of $100\times 100$. All other inversion results on different meshes were produced by the same trained VINet.  In the following, let us visually compare the results obtained by TSVD, MFVI, and VINet when the mesh size is equal to $150\times 150$. Then, we will provide a detailed comparison of these methods with different mesh sizes.

\begin{figure}[t]
	\centering
	\includegraphics[width=0.9\textwidth]{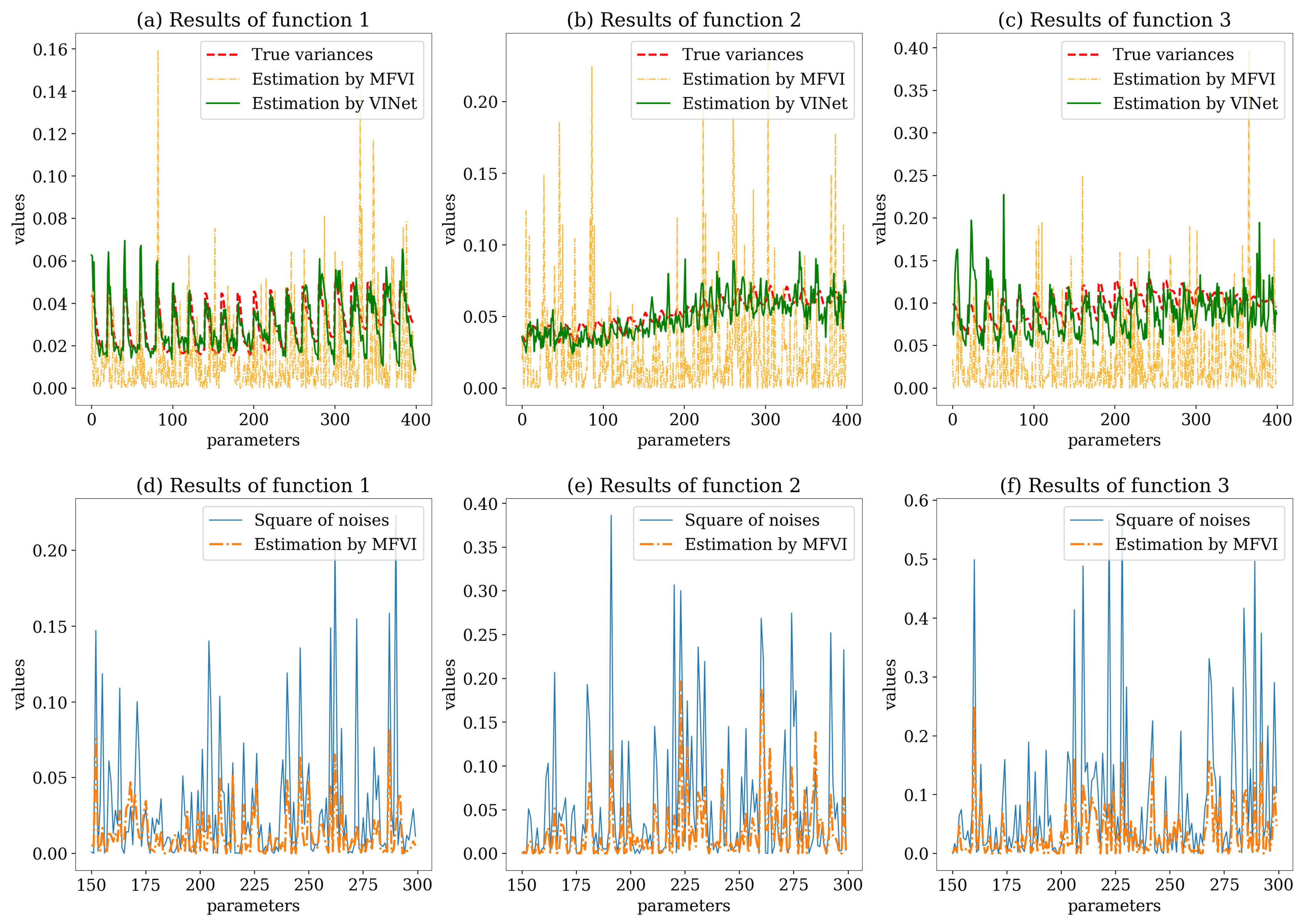}\\
	\vspace*{-0.6em}
	\caption{Comparison of the true noise variances and estimated noise variances by the MFVI and VINet (output of SNet).
		(a) (b) (c) Noise variances of all parameters for three different testing examples.
		(d) (e) (f) Noise variances estimated by MFVI and the square real noises of 
		some parameters. }\label{NoiseCompareType2TestData}
\end{figure}

Let us first compare the estimates of the noise variances given by the MFVI and VINet. Because type 2 test datasets deviate from the training datasets (more challenging than type 1 ones), we only show the results of type 2 test datasets to illustrate the effectiveness of the proposed VINet. For interested readers, we exhibit the results of type 1 test datasets in the Appendix. In Figure \ref{NoiseCompareType2TestData}, we depict the estimated noise variances by the MFVI and VINet and the background true noise variances for three randomly selected testing examples. In the first line of Figure \ref{NoiseCompareType2TestData}, the estimated variances of all 400 measurement points are given, with the dashed red line representing the true noise variances, the dash-dotted light orange line representing the estimation given by MFVI, and the solid green line representing the estimation given by the SNet in VINet. Obviously, the SNet provides evidently more accurate estimates of the noise variances for all of the three type 2 test datasets.  
In the second line of Figure \ref{NoiseCompareType2TestData}, we use the solid blue line and dash-dotted orange line to represent the square real noises, i.e., $(\bm{d} - \bm{d}_c)^2$, and the estimations given by MFVI, respectively. 
Obviously, the MFVI method only catches some features of the square real noises, but not the true noise variances $0.1^2\bm{d}_c^2$. From this figure, we see that the real noise variances $0.1^2\bm{d}_c^2$ can be visually much different from one particular realization $(\bm{d} - \bm{d}_c)^2$. Hence, it is a difficult task to estimate $0.1^2\bm{d}_c^2$ from the noisy data. Based on one piece of noisy data, the MFVI method can hardly give a reasonable estimate. However, by learning from the training datasets, the SNet captures the main features of the true noise variances.

\begin{figure}[thbp]
	\centering
	\includegraphics[width=0.8\textwidth]{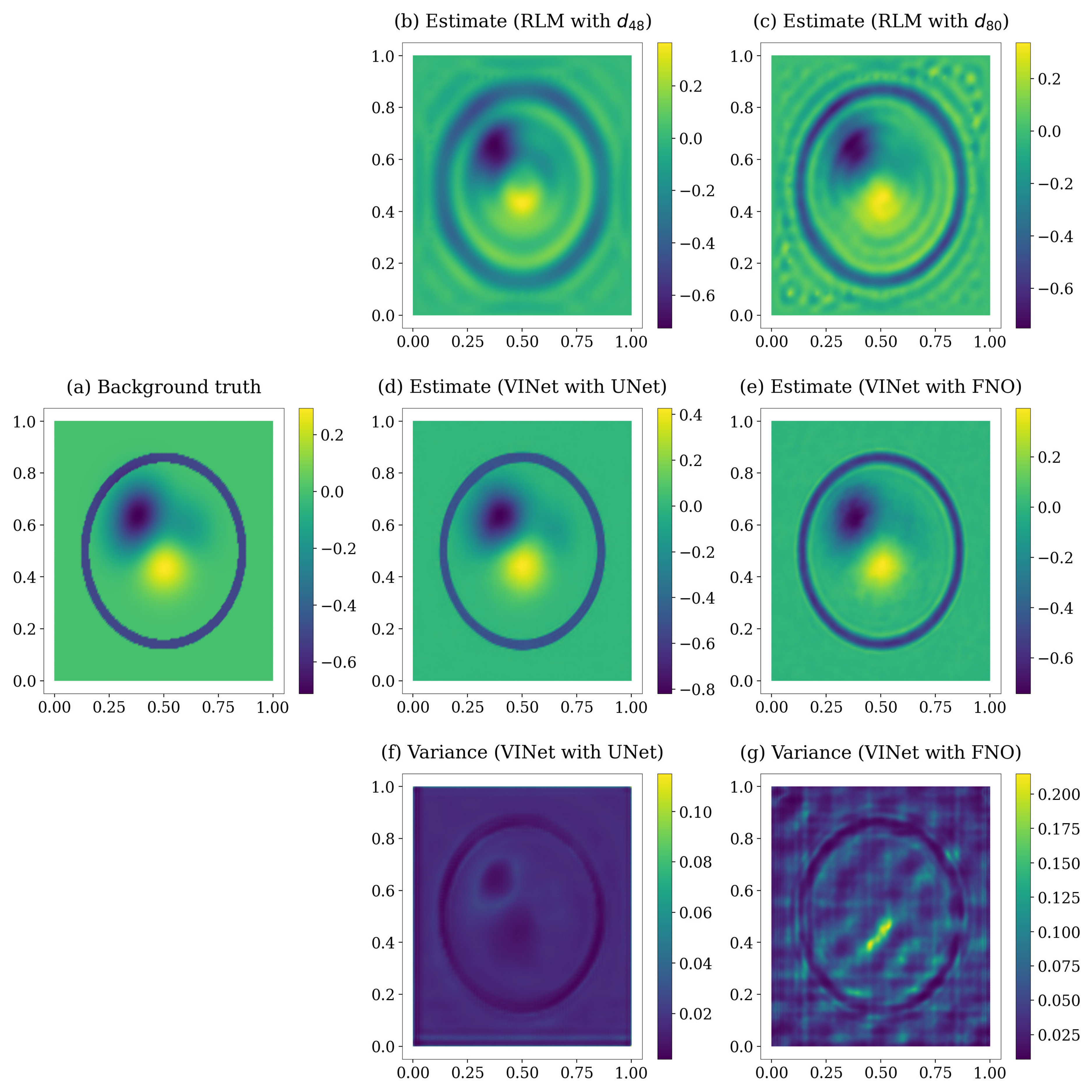}\\
	\vspace*{-0.6em}
	\caption{(a) Background true function randomly selected from type 2 testing datasets; 
	(b) Estimate obtained by the TSVD method; (c) Estimate obtained by the MFVI method; 
	(d) Estimate obtained by the VINet with U-ENet; (e) Estimate obtained by the VINet with FNO-ENet; 
	(f) Variance function obtained by the VINet with U-ENet; 
	(g) Variance function obtained by the VINet with FNO-ENet.}\label{InvResultsType2Test}
\end{figure}

Secondly, let us show the background truth and estimates obtained by the TSVD, MFVI, VINet with U-ENet, and VINet with FNO-UNet for one randomly selected data from the type 2 testing datasets. As before, the results of type 1 testing datasets are provided in the Appendix. In Figure \ref{InvResultsType2Test}, we exhibit the true function in sub-figure (a) and the estimated functions by the TSVD, MFVI, VINet with U-ENet, and VINet with FNO-ENet in sub-figures (b)-(e), respectively. Obviously, the TSVD and MFVI methods can only provide a rough estimate of the true parameter. Comparatively, no matter whether using U-ENet or FNO-ENet, the proposed VINet captures more detailed structures of the true parameter. In addition, we show the estimated variance functions of VINet in sub-figures (f) and (g). 
For this example, compared with the VINet with U-ENet, the VINet with FNO-ENet gives a less accurate estimate of the dark part (lower right corner) of the true parameter. The value of the darkest part in sub-figure (e) is approximately $0.5$. However, the value of the darkest parts in sub-figures (a) and (d) is below $0$, which reflects that the VINet with FNO-ENet provides a relatively less accurate estimate for the darkest parts compared with the VINet with U-ENet. The variance functions shown in sub-figures (f) and (g) just characterize this inaccuracy, i.e., for the region of darkest parts, the values of the variance function in sub-figure (g) are nearly three times higher than those in sub-figure (f).

\begin{table}[thbp]
	\caption{Average relative errors of the estimates obtained by the TSVD, MFVI, VINet (U-ENet), and VINet (FNO-ENet)
		with mesh sizes $n=\{30\times30, 50\times50, 100\times100, 150\times150, 200\times200, 250\times250, 300\times300 \}$
		for the type 2 test datasets.}
	\begin{center}
		\begin{spacing}{1.5}
			\begin{tabular}{c|ccccc}
				\Xhline{1.1pt}
				   & TSVD & MFVI & VINet (U-ENet) & VINet (FNO-ENet) \\
				\hline
				Relative error ($n=30^{2}$) & 0.1484 & 0.3287 & 0.0881 & 0.0802  \\
				\hline
				Relative error ($n=50^{2}$) & 0.1439 & 0.3191 & 0.0940 &  0.0849	\\
				\hline
				Relative error ($n=100^{2}$) & 0.1448 & 0.3218 & 0.0980 &  0.0881	\\
				\hline 
				Relative error ($n=150^{2}$) & 0.1455 & 0.3279 & 0.0989 &  0.0888	\\
				\hline
				Relative error ($n=200^{2}$) & 0.1456 & 0.3614 & 0.0996 &  0.0894	\\
				\hline 
				Relative error ($n=250^{2}$) & 0.1482 & 0.3535 & 0.0997 &  0.0894	\\
				\hline 
				Relative error ($n=300^{2}$) & 0.1482 & 0.3520 & 0.0999 & 0.0896	\\
				\Xhline{1.1pt}
			\end{tabular}
		\end{spacing}
	\end{center}\label{ResolutionCompareType2Ex1}
\end{table}

\begin{table}[thbp]
	\caption{Average computing times of the TSVD, MFVI, VINet (U-ENet), and VINet (FNO-ENet)
		with mesh sizes $n=\{30\times30, 50\times50, 100\times100, 150\times150, 200\times200, 250\times250, 300\times300 \}$
		for the type 2 test datasets. The computing time is measured in seconds.}
	\begin{center}
		\begin{spacing}{1.5}
			\begin{tabular}{c|ccccc}
				\Xhline{1.1pt}
				& TSVD & MFVI & VINet (U-ENet) & VINet (FNO-ENet) \\
				\hline 
				Average time ($n=30^{2}$) & 0.48 & 25.69 & 0.0296 & 0.0122  \\
				\hline 
				Average time ($n=50^{2}$) & 1.10 & 52.72 & 0.0296 &  0.0125 \\
				\hline 
				Average time ($n=100^{2}$) & 5.65 & 330.72 & 0.0298 &  0.0203 \\
				\hline 
				Average time ($n=150^{2}$) & 15.21 & 655.51 & 0.0299 &  0.0336 \\
				\hline 
				Average time ($n=200^{2}$) & 31.71 & 1092.62 & 0.0305 &  0.0436 \\
				\hline 
				Average time ($n=250^{2}$) & 55.39 & 2071.76 & 0.0305 &  0.0594 \\
				\hline 
				Average time ($n=300^{2}$) & 86.82 & 2921.92 & 0.0311 & 0.0729  \\
				\Xhline{1.1pt}
			\end{tabular}
		\end{spacing}
	\end{center}\label{TimeCompareType2Ex1}
\end{table}

At last, let us provide a detailed comparison of these methods with various discrete meshes. In Table \ref{ResolutionCompareType2Ex1}, we list the average relative errors of the estimates obtained by the TSVD, MFVI, VINet with U-ENet, and VINet with FNO-ENet for type 2 test datasets when the mesh sizes are equal to $30\times30$, $50\times50$, $100\times100$, $150\times150$, $200\times200$, $250\times250$, and $300\times300$. From Table \ref{ResolutionCompareType2Ex1}, we see that all four methods could provide stable estimates with different discretized dimensions. It should be mentioned that the VINets with both U-ENet and FNO-ENet were trained when the discrete mesh size was $n=100\times100$. The relative errors of the proposed VINets were much smaller than those of the TSVD and MFVI methods. The TSVD method could only provide an estimate of the function parameter and the MFVI method gave the estimates with the highest relative errors. From Figure \ref{InvResultsType2Test}, we actually found that the two VINets provided more details compared with the TSVD and MFVI methods. In Table \ref{TimeCompareType2Ex1}, we also provided the recorded computing times on CPU since the TSVD and MFVI methods were tested on CPU. If we ran VINet on GPU, the running time would be nearly ten times faster than that shown in Table \ref{TimeCompareType2Ex1}. Obviously, in the inference stage, the two VINets were hundreds or thousands of times faster than the TSVD and MFVI methods. For the VINets with U-ENet and FNO-ENet, we took $1789.89$ and $1536.15$ seconds (approximately $30$ and $27$ minutes) to train the neural networks, respectively. The training time plus the inference time was even smaller than that of the MFVI with mesh sizes $n = \{250\times 250, 300\times300\}$.

\subsection{An inverse source problem}\label{HelmholtzModel}

\subsubsection{Basic settings}\label{HelmholtzModelSec1}

In this subsection, we consider an inverse source problem for the Helmholtz equation with
multi-frequency data which was studied in \citep{Bao2015IP}\footnote{Deep learning based methods are employed for statistical inverse problems under the finite-dimensional setting
for CT image reconstruction, PET image reconstruction, and MRI reconstruction \citep{Arridge2019ActaNumerica}.
For these problems, the finite-dimensional formulation is widely employed.
However, to the best of our knowledge, there seems little such types of work focused on multi-frequency inverse source
problems governed by the Helmholtz equation (one of the typical and important mathematical models usually defined on infinite-dimensional space \citep{Bao2015IP}). So we compare the performance of the proposed method with the classical recursive linearization method (RLM), which should be sufficient to reveal the mechanism and superiority of the proposed method.}. Specifically speaking, we consider the following Helmholtz equation:
\begin{align}
\begin{split}
\Delta w + \kappa^2 w & =u \quad\text{in }\mathbb{R}^2, \\
\partial_r w - i\kappa w & = o(r^{-1/2}) \quad\text{as }r=|x|\rightarrow\infty,
\end{split}
\end{align}
where $\kappa$ is the wavenumber, $w$ is the acoustic field, and $u$ is the acoustic source supported in an open bounded domain $\Omega = [0,1]^2$. To simulate the problem defined on the infinity domain $\mathbb{R}^2$, we adopt the perfectly matched layer (PML) technique and reformulate the scattering problem in a bounded domain. Let $\bm{x} = (x, y) \in\mathbb{R}^2$. Denote $D$ as a rectangle containing $\Omega$, and let $d_1$ and $d_2$ be the thickness of the PML layers in the $x$ and $y$ axes, respectively. Let $\partial D$ be the boundary of $D$.
Define $s_1(x) = 1+i\sigma_1(x)$ and $s_2(y)=1+i\sigma_2(y)$, where $\sigma_1, \sigma_2$ are positive continuous functions satisfying $\sigma_1(x)=0, \sigma_2(y)=0$ in $\Omega$. For more details of the PML, please refer to \citep{Bao2010IP}.

Following the general PML technique, we can deduce the truncated PML problem: 
\begin{align}
\begin{split}
& \nabla\cdot(s\nabla w) + \kappa^2 s_1s_2 w = u \quad\text{in }D, \\
& w = 0\quad\text{on }\partial D,
\end{split}
\end{align}
where the diagonal matrix $s=\text{diag}(s_2(y)/s_1(x), s_1(x)/s_2(y))$. The forward operator related to $\kappa$ is defined by the Helmholtz equation $H_{\kappa}(u) = (w(\bm{x}_1),\ldots,w(\bm{x}_{N_d}))^T$ with $\{\bm{x}_i\}_{i=1}^{N_d}\in\partial\Omega$ and $u\in\mathcal{H}_u:=L^2(\Omega)$. Since we consider the multi-frequency case, i.e., a series of wavenumbers $0 < \kappa_1 < \kappa_2 < \cdots < \kappa_{N_f} < \infty$ are considered, the forward operator has the following form:
\begin{align}
Hu = (H_{\kappa_1}(u),\ldots,H_{\kappa_{N_f}}(u)) \in \mathbb{R}^{N_d\times N_{N_f}}.
\end{align}
Similar to the simple smoothing model, we generate the training and testing data based on a fine mesh with the number of grid points equal to $500\times 500$. In the inverse computational stage, the related functions are discretized by the finite element method on a series of regular meshes with the number of grid points 
$n = \{30\times30, 100\times100, 150\times150, 200\times200, 250\times250, 300\times300\}$. 
Concerning the operators $\mathcal{C}_{\epsilon_0}$ and $\mathcal{C}_p$, we adopt the same settings as those in Subsection \ref{SimpleModel} for the simple smoothing model.

For the CECInv layer, we employ the classical recursive linearization method (RLM) reviewed in \citep{Bao2015IP} with $\kappa = 10,20,30,40,48$. The PDEs are discretized on a regular mesh with the number of grid points equal to $30\times30$, which makes it possible to directly compute the inverse of all the sparse matrices generated by the finite element method. Obviously, for a source function, the data generated by PDEs on a fine mesh with grid points equal to $500\times 500$ are different from the data obtained by PDEs on a coarse mesh, especially when the wavenumber is large. For this example, we take the data generated by PDEs on the fine mesh with random Gaussian noise as the noisy input data. The accurate data are recognized as the data generated by solving PDEs with the coarse mesh, i.e., a mesh with size $n=30\times30$. Hence, the DNet and SNet actually learn the modeling error and random noise together in this example. Please refer to \citep{Calvetti2018IP,Jia2019IP} for more relevant references on model error learning algorithms.

Before presenting the numerical results, we clarify how the training and testing datasets are constructed. Similar to the strategies employed in \citep{Jia2019IP}, we introduce the following formula:
\begin{align}\label{utildeSample}
\tilde{u}(x,y) := \sum_{k = 1}^{3}a_{k}^{1} \exp\bigg(-a_{k}^{2}(x - a_{k}^{3})^{2}
- a_{k}^{4}(y - a_{k}^{5})^{2}\bigg),
\end{align}
where
\begin{align*}
& a_{k}^{1} \sim U[-1,1], \quad a_{k}^{2}, a_{k}^{4} \sim U[50,60], \quad a_{k}^{3}, a_{k}^{5} \sim U[0.3,0.7],
\end{align*}
with $k=1,2,3$ and $U[a,b]$ being the uniform distribution between $a$ and $b$. To demonstrate the effectiveness of the proposed method, we consider the following nonsmooth source function:
\begin{align}\label{ex1fun}
\begin{split}
u(x,y) = 
\begin{cases}
\tilde{u}(x,y) \quad \text{for }|\bm{x} - \bm{x}_0|<r_1, \\
-0.5 \,\,\,\,\,\quad \text{for }r_2 \leq |\bm{x} - \bm{x}_0| \leq r_1, \\
0.0 \,\,\,\,\qquad \text{for }|\bm{x} - \bm{x}_0| > r_2,
\end{cases}
\end{split}
\end{align}
where $\bm{x}_0 = (0.5,0.5)^T$, $r_1\sim U[0.32, 0.4]$ and $r_2\sim U[0.42, 0.48]$. Based on Eq. (\ref{ex1fun}), we generate $1000$ source functions for training and $100$ source functions for testing. For testing the generalization properties, we also test the performance of the trained model on the following source function:
\begin{align}\label{funMatlab1}
\begin{split}
u^{\dag}(x,y) = u(x,y) + \tilde{u}(x,y),
\end{split}
\end{align}
where $u$ and $\tilde{u}$ are two functions generated according to Eq. (\ref{ex1fun}) and Eq. (\ref{utildeSample}), respectively. The functions defined by (\ref{funMatlab1}) and (\ref{ex1fun}) exhibit similar properties but the function defined in (\ref{funMatlab1}) can never be generated from Eq. (\ref{ex1fun}), which mimics the following situation: the constructed learning examples usually cannot contain completely correct information, i.e., prior knowledge may exhibit bias. In the following, we call the $100$ test source functions generated from Eq. (\ref{ex1fun}) and Eq. (\ref{funMatlab1})  as type 1 and type 2 test datasets, respectively.  One thing should be emphasized is that the training datasets can be constructed by employing more sophisticated methods and are not restricted to the presented simple parametric form.

Similar to the simple smooth example, the ENet has been chosen to be either U-Net or FNO. If we choose ENet to be a U-Net, we will employ the same settings as in Subsection \ref{SimpleModel} and we will still call the ENet U-ENet in the following. When choosing ENet to be a FNO, we also take the same network structures but set the maximal number of Fourier modes to be $32$ for each Fourier integral operator, and the original inputs are lifted to a higher dimensional representation with a dimension equal to $64$. As before, we call ENet with FNO structure FNO-ENet in this section. Different from the simple smooth example, the current CECInv layer uses the measured data of 
$5$ different wavenumbers, i.e., $\kappa = 10,20,30,40,48$. Due to the small number of wavenumbers, 
it is better to employ all of the $5$ estimated functions obtained by different wavenumbers. 
In addition, the variance function of the data obtained from DNet needs to be propagated to the ENet. 
Hence, the input layer of ENet for the current case has $6$ channels (the $6$ functions are projected on the same 
mesh) and the output layer contains $2$ channels, i.e., the mean function and function parameter in the covariance operator.

The training procedure is similar to that in Subsection \ref{SimpleModel}. We will first train DNet and SNet, and then train ENet and SNet with the trained DNet. For the first stage of training, we used the Adam optimizer to train for 80 epochs with an initial learning rate of 0.0001, which is halved every 10 epochs. For training U-ENet and FNO-ENet, we also used the Adam optimizer to train for 25 epochs with an initial learning rate of 0.0001, which is halved every 10 epochs.

\subsubsection{Numerical results}\label{HelmholtzModelSec2}

\begin{figure}[t]
	\centering
	\includegraphics[width=0.9\textwidth]{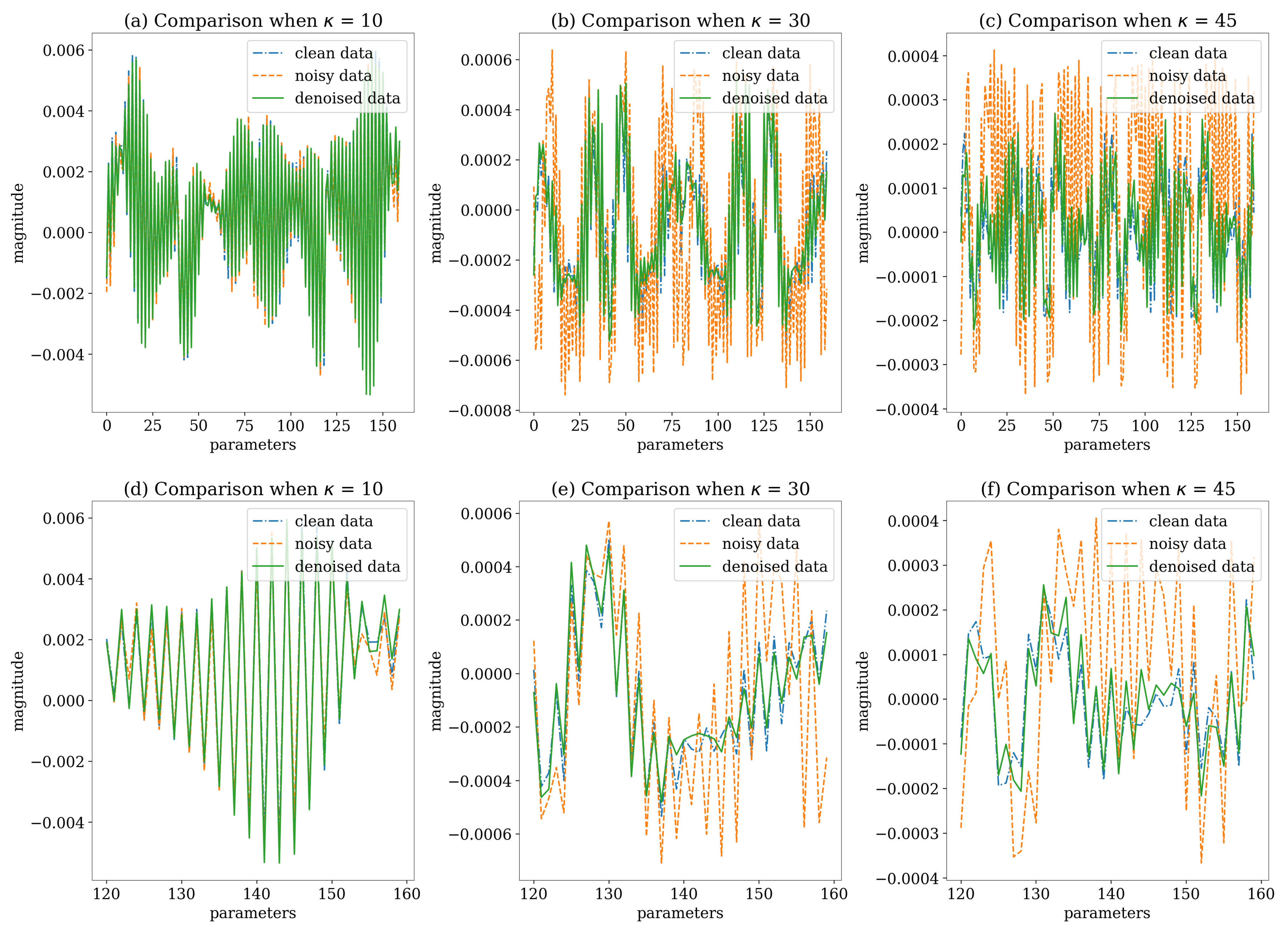}\\
	\vspace*{-0.6em}
	\caption{Noisy data (data generated by an accurate forward solver with random Gaussian noise), clean data (data generated by a rough forward solver), and denoised data (the output of the DNet) for a randomly selected source function $u$ from type 2 test datasets are shown in (a) data when the wavenumber $\kappa=10$; (b) data when the wavenumber $\kappa=30$; (c) data when the wavenumber $\kappa=45$; (d), (e), and (f) are partly enlarged versions of (a), (b), and (c), respectively.}\label{d_invert_test2}
\end{figure}

We compare the results of the classical recursive linearization method (RLM) \citep{Bao2015IP}, VINet with U-ENet, and VINet with FNO-ENet. The RLM uses multiple frequency data, i.e., data obtained at a series of wavenumbers. Low frequency data make the RLM stable with respect to the initialization (accurate initial function is usually hard to determine) and the high frequency data ensure the accuracy of the RLM. For our experiments, the data for each wavenumber was added $5\%$ random Gaussian noise. For the RLM, we employed two types of data, i.e., data $\bm{d}_{48}$ obtained when the wavenumbers are  $\kappa = 1,\cdots,48$, and data $\bm{d}_{80}$ obtained when the wavenumbers are $\kappa=1,\cdots,80$. For the VINets with U-ENet or FNO-ENet, we used $1000$ training data pairs that contain the $1000$ training source functions generated by formula (\ref{ex1fun}) and the data $\bm{d}_{48}$ generated by solving PDEs on a fine mesh, i.e., mesh with size $n=500\times 500$. For training DNet, we also needed the data $\bm{\tilde{d}}_{48}$ generated by solving Helmholtz equations on a mesh with  mesh size $n = 30\times30$, which is a cheap computational procedure due to the small mesh size (all sparse linear solvers transform to matrix multiplications). The DNet in this example actually learns to remove two types of noises:  $5\%$ random noises and the discrepancies of different discrete models at measurement points. Lastly, we should emphasize that the VINets with U-ENet and FNO-ENet were all trained when the mesh size was $n=150\times150$. All of the results for other meshes were generated by VINets through the rescale layers of U-ENet or the FNO-ENet directly.

\begin{figure}[t]
	\centering
	\includegraphics[width=0.9\textwidth]{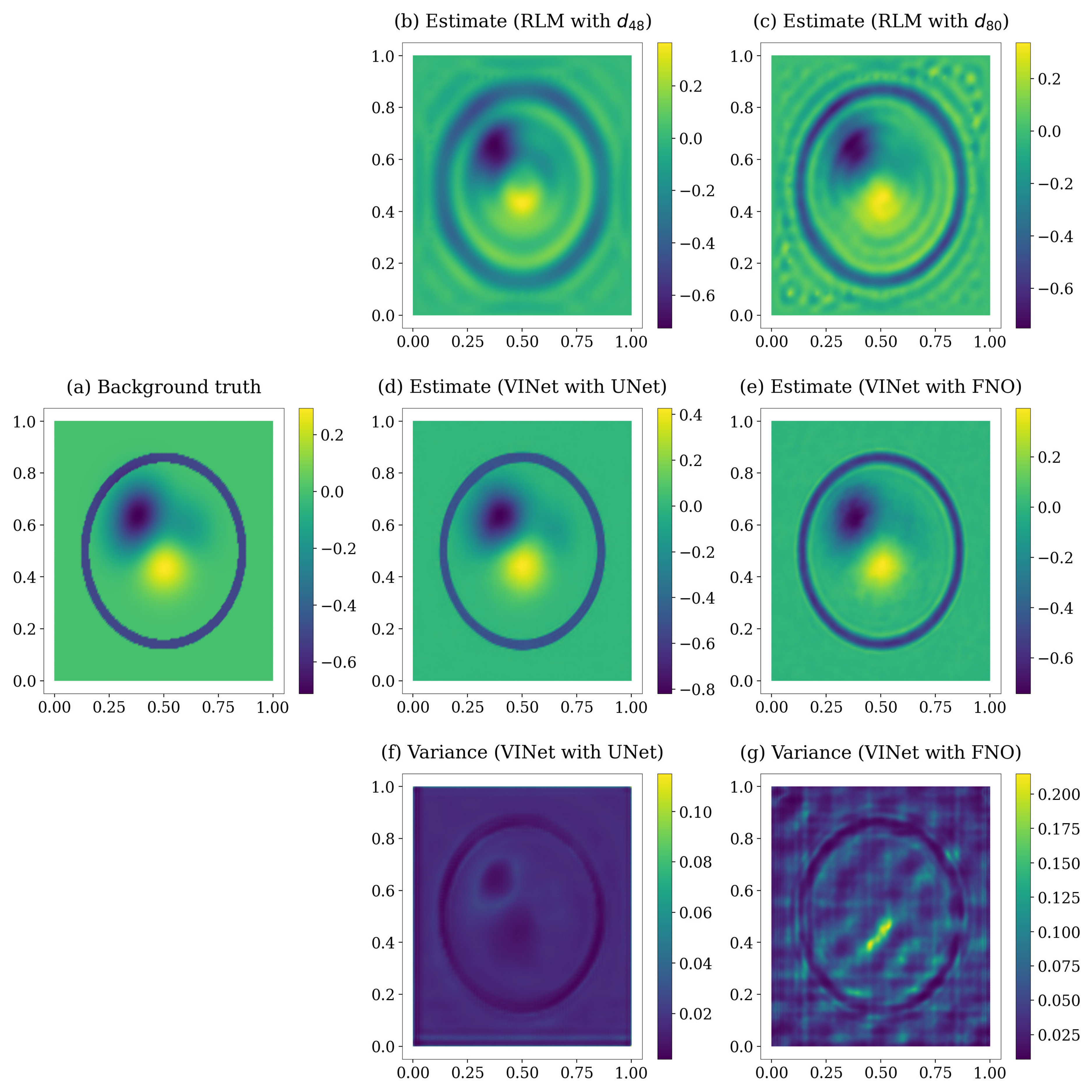}
	\vspace*{-0.6em}
	\caption{(a) One of the background true function in type 2 test datasets;
		(b) Estimate obtained by RLM with $\bm{d}_{48}$; (c) Estimate obtained by RLM with $\bm{d}_{80}$; 
		(d) Estimate obtained by VINet with U-ENet; (e) Estimate obtained by VINet with FNO-ENet; 
		(f) Estimated variance function by VINet with U-ENet; 
		(g) Estimated variance function by VINet with FNO-ENet.}\label{estUtype2}
\end{figure}

As in Subsection \ref{SimpleModel2}, we only provide the results for type 2 testing datasets (which have a different generation mechanism than the training datasets) and leave the results of type 1 testing datasets (which have the same generation mechanism as the training datasets) in the Appendix. Firstly, let us provide some visual comparisons to acquire some intuitive impressions about the estimates of different methods. In Figure \ref{d_invert_test2}, we show the results given by DNet with input data $\bm{d}_{48}$ that is randomly selected from the type 2 test datasets. In all of the sub-figures of Figure \ref{d_invert_test2}, we use the dashed blue line, dashed orange line, and solid green line to represent the clean data (solution of the rough forward solver with mesh size $n=30\times30$), noisy data (solution of the accurate forward solver with mesh size $n=500\times500$), and the output of the DNet, respectively. Sub-figures (a) (b) and (c) provide comparisons when the wavenumbers are chosen to be $10$, $30$, and $45$, respectively. To provide more details, we exhibit sub-figures (d) (e) and (f) that give the partly enlarged versions of sub-figures (a) (b) and (c), respectively. When the wavenumber is small, i.e., $\kappa = 10$, we see that the clean data, noisy data, and output of DNet exhibit no significant differences. The reason is that we only need a rough forward solver when the wavenumbers are small. As the wavenumbers increase, we need a larger mesh to obtain an accurate forward solver. From sub-figures (b), (c), (e), and (f), we can clearly find that the noisy data deviate significantly from the clean data and the outputs of DNet are visually similar to the clean data, which reflects that DNet has learned the sophisticated model error mechanism between the accurate and rough forward solvers.  

In Figure \ref{estUtype2}, we show the results of one randomly selected test example from the type 2 test datasets. In sub-figure (a), we exhibit the background true source function. In sub-figures (b) and (c), we provide the estimated source functions of RLM with data $\bm{d}_{48}$ and $\bm{d}_{80}$, respectively. It is obvious that the classical RLM can hardly handle the discontinuity, so the sharp circle cannot be accurately estimated even when the highest wavenumber is set to $80$. In sub-figures (d) and (e), we see that the VINets with U-ENet and FNO-ENet both yield evidently more accurate estimates compared to the classical RLM. In sub-figures (f) and (g), we show the variance functions obtained from the VINets with U-ENet and FNO-ENet. The high values of variance functions appear at the changing parts of the source function, especially at the discontinuous boundary, which is reasonable since discontinuous boundaries are intuitively harder to estimate.

\begin{table}[t]
	\caption{Average relative errors of the estimates obtained by the RLM (using data $\bm{d}_{48}$ and $\bm{d}_{80}$), 
	VINet with U-ENet, and VINet with FNO-ENet with mesh sizes
	$n = \{30\times30, 100\times100, 150\times150, 200\times200, 250\times250, 300\times300\}$ for the type 2 test datasets.}
	\begin{center}
		\begin{spacing}{1.5}
			\begin{tabular}{c|ccccc}
				\Xhline{1.1pt}
				& RLM($\bm{d}_{48}$) & RLM($\bm{d}_{80}$) & VINet(U-ENet) & VINet(FNO-ENet) \\
				\hline
				Relative error ($n=30^{2}$) & 0.3670 & 0.3515 & - & -  \\
				\hline
				Relative error ($n=100^{2}$) & 0.2727 & 0.2763 & 0.0693 &  0.0969	\\
				\hline 
				Relative error ($n=150^{2}$) & 0.2623 & 0.2414 & 0.0630 &  0.0401	\\
				\hline
				Relative error ($n=200^{2}$) & 0.2596 & 0.2246 & 0.0625 &  0.0410	\\
				\hline 
				Relative error ($n=250^{2}$) & 0.2581 & 0.2191 & 0.0630 &  0.0425	\\
				\hline 
				Relative error ($n=300^{2}$) & 0.2575 & 0.2149 & 0.0645 & 0.0450	\\
				\Xhline{1.1pt}
			\end{tabular}
		\end{spacing}
	\end{center}\label{resultsHelCompareType2}
\end{table}

\begin{table}[thbp]
	\caption{Average computing times of the RLM (using data $\bm{d}_{48}$ and $\bm{d}_{80}$), 
		VINet with U-ENet, and VINet with FNO-ENet with mesh sizes
		$n = \{40\times40, 100\times100, 150\times150, 200\times200, 250\times250, 300\times300\}$ for the type 2 test datasets.
		The computing time is measured in seconds.}
	\begin{center}
		\begin{spacing}{1.5}
			\begin{tabular}{c|ccccc}
				\Xhline{1.1pt}
				& RLM($\bm{d}_{48}$) & RLM($\bm{d}_{80}$) & VINet(U-ENet) & VINet(FNO-ENet) \\
				\hline 
				Average time ($n=30^{2}$) & 0.71 & 1.26 & - & -  \\
				\hline 
				Average time ($n=100^{2}$) & 17.27 & 29.09 & 0.1929 &  0.3882 \\
				\hline 
				Average time ($n=150^{2}$) & 53.21 & 83.32 & 0.1795 &  0.4011 \\
				\hline 
				Average time ($n=200^{2}$) & 126.01 & 210.33 & 0.1799 &  0.4313 \\
				\hline 
				Average time ($n=250^{2}$) & 245.93 & 407.71 & 0.1800 &  0.4697 \\
				\hline 
				Average time ($n=300^{2}$) & 495.78 & 827.83 & 0.1801 & 0.5295  \\
				\Xhline{1.1pt}
			\end{tabular}
		\end{spacing}
	\end{center}\label{TimeCompareType2Ex2}
\end{table}

Next, let us give a detailed comparison of RLM and VINet with various discrete meshes. In Table \ref{resultsHelCompareType2}, we show the averaged relative errors on different meshes obtained by RLM with $\bm{d}_{48}$, RLM with $\bm{d}_{80}$, VINet with U-ENet, and VINet with FNO-ENet, respectively. From the table, it can be observed that the discrete level indeed influences the performance of the classical RLM; the relative errors of the results obtained by RLM are larger than $1\%$ when the discrete mesh size is smaller than $n=200\times200$. Using higher frequency data $\bm{d}_{80}$, the inversion results indeed have lower relative errors. However, due to the limited measured data (we only assume the measurement points are equally spaced along the boundary with a number equal to $80$), we see that the inverse accuracy increases slowly when more wavenumbers are employed. When the VINets with U-ENet and FNO-ENet are employed, the relative errors are around $6\%$ and $4\%$, respectively. Both types of VINets provide stable inversion results for different meshes, except the VINet with FNO-ENet when the mesh size is equal to $100\times100$. This can be rationally explained by the fact that the function parameter $u$ has a discontinuous sharp boundary which leads to larger numerical errors when the discrete dimension is small (compared with the frequency truncate level of FNO). When the function parameter $u$ has better regularity, there seems to be no such phenomenon, as seen in Table \ref{ResolutionCompareType2Ex1}.

At last, in Table \ref{TimeCompareType2Ex2}, we give the computing times based on averages of the type 2 test datasets. Obviously, the inference times of both types of VINet are significantly shorter than those of the classical RLM. It should be mentioned that all of the computation in the inference stage is performed on the CPU (not the GPU) since RLM seems more suitable for running on the CPU. If we calculate the inference of VINet on GPU, the inference times of VINet will be around 0.015 seconds for both types of VINets. For training the VINets with U-ENet and FNO-ENet, we take 33950 and 39687 seconds (approximately 9.4 and 11 hours) respectively. The training stage needs only to be calculated once, and the learned computationally cheap model can be repeatedly employed for similar inversion tasks. Therefore, using large training times to exchange much slower and more accurate inversion results is reasonable.

\section{Conclusion}\label{Conclusion}

In this paper, we have studied inverse problems governed by PDEs with non-i.i.d. noise assumptions under the Bayesian analysis framework. We have provided a rigorous illustration of Bayes' formula under the non-i.i.d. noise setting with the conjugate prior measure. To construct an efficient approximate sampling algorithm, we have deduced the mean-field based variational inference (MFVI) method with infinite-dimensional parameters. However, this deduced algorithm is time-consuming as a large amount of computationally expensive PDEs need to be solved. To address this issue, we have further constructed a parametric form of posterior measures based on an intuitive analysis of the approximate measure proposed in the MFVI method and proposed a general inference framework for inverse problems named VINet. By introducing measure equivalence conditions, we have generalized the evidence lower bound argument to the infinite-dimensional setting, which yields a computable training method. We have also provided a possible parametric strategy of the parametric posterior measure, which reflects the equivalent measure assumption that can be satisfied by designing parametric strategies appropriately. Finally, we have provided a possible specific network structure that yields a concrete VINet. The VINet has been applied to two inverse problems that are governed by a simple elliptic equation and the Helmholtz equation. For both inverse problems, the VINet was able to learn prior information of interested functions and noises from training examples efficiently and provide proper estimates quickly.

The current VINet framework is proposed based on the analysis of the MFVI method under the conjugate prior measure, which restricts the form of the posterior measures. Especially for nonlinear inverse problems, the VINet can only manage linearized problems and provide an approximate probability measure. As for the MFVI, the approximate probability measure may be inaccurate for some highly nonlinear problems. Hence, generalizing the semi-conjugate VI methods to the infinite-dimensional setting and constructing the corresponding DNN-based inference method is worthy of investigation in future work. For the current preliminary studies on numerical aspects, we only employ the U-Net with a rescale layer and FNO, which is not designed for particular inverse problems. For specific inverse problems such as inverse scattering problems, special neural networks \citep{Khoo2019SISC} may be employed to construct the INet in our general statistical inference framework. For the CECInv layer, we only choose a rough discretized classical inverse algorithm for simplicity. It is interesting to choose the CECInv layer as some classical inverse methods with surrogates, e.g., polynomial approximations \citep{Li2014SISC} and Fourier neural operator \citep{Zongyi2020ICLR}. We will conduct further investigations on these issues in our future research.

\section{Appendix}\label{appendixSection}

\subsection{Bayesian well-posedness theory}\label{BayesWellPosed}

For inverse problems, it is necessary to find a probability measure $\mu^{\bm{d}}$ on $\mathcal{H}$ (a separable Hilbert space), which is known as the posterior probability measure, and is specified by its density with respect to a prior probability measure $\mu_{0}$. Here, we assume the data $\bm{d}\in\mathbb{R}^{N_d}$ as in the main text.
The Bayes' formula on the Hilbert space is defined by
\begin{align}\label{BayesianFormula1}
	\frac{d\mu^{\bm{d}}}{d\mu_{0}}(x) = \frac{1}{Z_{\mu^{\bm{d}}}}\exp\big( -\Phi(x;\bm{d}) \big),
\end{align}
where $\Phi(\cdot;\bm{d}): \mathcal{H}\rightarrow \mathbb{R}$ is a continuous function, and $\exp\big( -\Phi(\cdot;\bm{d}) \big)$ is integrable with respect to $\mu_{0}$. The constant $Z_{\mu^{\bm{d}}}$ is chosen to ensure that $\mu^{\bm{d}}$ is a probability measure. Under the current setting, Assumption 1 and Theorems 15--16 in \citep{Dashti2014} can be stated as follows.

\begin{assumption}\label{AssumptionStuart}
(Assumption 1 in \citep{Dashti2014})
Assume that $\Phi \in C(\mathcal{H} \times \mathbb{R}^{N_d}; \mathbb{R})$. Further assume that there are functions $M_i : \mathbb{R}^{+} \times \mathbb{R}^{+} \rightarrow \mathbb{R}^{+}$, $i = 1,2$, which are monotonic and non-decreasing separately in each argument, and with $M_2$ strictly positive. For all $x \in \mathcal{H}$ and $\bm{d}, \bm{d}_1, \bm{d}_2 \in B_{\mathbb{R}^{N_d}}(0,r) := \{\bm{d} \in \mathbb{R}^{N_d} \, : \, \|\bm{d}\|_{\mathbb{R}^{N_d}} < r\}$, the following holds:
\begin{align*}
\Phi(x; \bm{d}) & \geq -M_1(r, \|x\|_{\mathcal{H}}), \\
|\Phi(x;\bm{d}_1) - \Phi(x;\bm{d}_2)| & \leq M_2(r,\|x\|_{\mathcal{H}})\|\bm{d}_1 - \bm{d}_2\|_{\mathbb{R}^{d}}. 
\end{align*}
\end{assumption}

\begin{theorem}
(Theorem 15 in \citep{Dashti2014})
Let Assumption \ref{AssumptionStuart} hold. Assume that $\mu_0(\mathcal{H}) = 1$ and that $\mu_0(\mathcal{H}\cap B) > 0$ for some bounded set $B$ in $\mathcal{H}$. Additionally, assume that for every fixed $r > 0$, $\exp\left( M_1(r, \|x\|_{\mathcal{H}}) \right) \in L_{\mu_0}^1(\mathcal{H};\mathbb{R})$, where $L_{\mu_0}^1(\mathcal{H};\mathbb{R})$ represents $\mathcal{H}$-valued integrable functions under the measure $\mu_0$. Then, for every $\bm{d}\in\mathbb{R}^{N_d}$, $Z_{\mu^{\bm{d}}}$ is positive and finite, and the probability measure $\mu^{\bm{d}}$ given by (\ref{BayesianFormula1}) is well-defined.
\end{theorem}

\begin{theorem}
(Theorem 16 in \citep{Dashti2014})
Let Assumption \ref{AssumptionStuart} hold. Assume that $\mu_0(\mathcal{H}) = 1$ and that $\mu_0(\mathcal{H}\cap B) > 0$ for some bounded set $B$ in $\mathcal{H}$. Additionally, assume that, for every fixed $r > 0$, 
\begin{align*} 
\exp\left( M_1(r, \|x\|_{\mathcal{H}}) \right) \left( 1+M_2(r, \|x\|_{\mathcal{H}})^2 \right) \in L_{\mu_0}^1(\mathcal{H};\mathbb{R}). 
\end{align*}
Then, there is a $C=C(r) > 0$ such that, for all $\bm{d},\bm{d}' \in B_{\mathbb{R}^{N_d}}(0,r)$,
\begin{align*}
	d_{\text{Hell}}(\mu^{\bm{d}}, \mu^{\bm{d}'}) \leq C \|\bm{d} - \bm{d}'\|_{\mathbb{R}^{N_d}}, 
\end{align*}
where $d_{\text{Hell}}(\cdot,\cdot)$ is the Hellinger distance between two probability measures.
\end{theorem}

\subsection{The infinite-dimensional variational inference theory}\label{AppendixGeneTheoryVI}

In this section, we provide a brief introduction to the infinite-dimensional variational inference theory with the mean-field assumption. For interested readers, we suggest reading the paper \citep{Jia2020SIAM}, which comprehensively describes the full ideas and theories of this topic. It should be emphasized that we have slightly improved the statement of Theorem 11 in \citep{Jia2020SIAM} to make the theory more applicable to practical problems. See Theorem \ref{theoremGeneral} below for details.

Remembering the Bayes' formula (\ref{BayesianFormula1}) shown in Subsection \ref{BayesWellPosed} and 
denoting $\Phi(x) := \Phi(x;\bm{d})$, 
the variational inference method can be formulated as solving the following optimization problem:
\begin{align}\label{geneMiniPro1}
\argmin_{\nu\in\mathcal{A}}D_{\text{KL}}(\nu||\mu^{\bm{d}}),
\end{align}
where $\mathcal{A}\subset\mathcal{M}(\mathcal{H})$ is a set of measures.
Here $\mathcal{M}(\mathcal{H})$ is the set of Borel probability measures on $\mathcal{H}$.

For a fixed positive constant $M$, we assume the variable $x=(x_1, \ldots, x_M)$. Restricted to the settings in the main text, the constant $M=2$ and $x=(u,\bm{\sigma})$. Generally, the Hilbert space $\mathcal{H}$ and set $\mathcal{A}$ can be specified as
\begin{align}
\mathcal{H} = \prod_{j=1}^{M}\mathcal{H}_{j}, \qquad
\mathcal{A} = \prod_{j=1}^{M}\mathcal{A}_{j},
\end{align}
where $\mathcal{H}_{j}, j=1,\cdots,M$ are separable Hilbert spaces and $\mathcal{A}_{j}\subset\mathcal{M}(\mathcal{H}_{j})$. Denote $\nu:=\prod_{i=1}^{M}\nu^{i}$ by a probability measure satisfying $\nu(dx) = \prod_{i=1}^{M}\nu^{i}(dx)$.
Based on these assumptions, we can rewrite the minimization problem (\ref{geneMiniPro1})  as
\begin{align}\label{mipro1}
\argmin_{\nu^{i}\in\mathcal{A}_{i}}D_{\text{KL}}\bigg(\prod_{i=1}^{M}\nu^{i}\big|\big|\mu^{\bm{d}}\bigg)
\end{align}
for suitable sets $\mathcal{A}_{i}$ with $i=1,2,\cdots,M$.
Now, we introduce the approximate probability measure $\nu$ given in (\ref{geneMiniPro1}) which is equivalent to $\mu_0$, defined by
\begin{align}
\frac{d\nu}{d\mu_{0}}(x) = \frac{1}{Z_{\nu}}\exp\big( -\Phi_{\nu}(x) \big).
\end{align}
Comparing with the finite-dimensional case, we may assume that the potential $\Phi_{\nu}(x)$ can be decomposed as
\begin{align}
\exp\left(-\Phi_{\nu}(x)\right) = \prod_{j=1}^{M}\exp\left(-\Phi_{\nu}^{j}(x_{j})\right),
\end{align}
where $x = (x_{1},\cdots,x_{M})$. As illustrated in \citep{Jia2020SIAM}, this approach does not incorporate the parameters in the prior probability measure into the hierarchical Bayes' model. Therefore, the following assumption needs to be introduced.

\begin{assumption}\label{assumption1}
Let us introduce a reference probability measure
\begin{align}
\mu_{r}(dx) = \prod_{j=1}^{M}\mu_{r}^{j}(dx_{j}),
\end{align}
which is equivalent to the prior probability measure with the following relation:
\begin{align}
\frac{d\mu_{0}}{d\mu_{r}}(x) = \frac{1}{Z_{0}}\exp(-\Phi^{0}(x)).
\end{align}
For each $j\in\{1,2,\cdots,M\}$, there exists a predefined continuous function $a_j(\epsilon, x_j)$ satisfying $\mathbb{E}^{\mu_r^j}[a_j(\epsilon,\cdot)] < \infty$, where $\epsilon\in(0, \epsilon_0^j)$ with $\epsilon_0^j$ being a small positive number and $x_j\in\mathcal{H}_j$. Additionally, we assume that the approximate probability measure $\nu$ is equivalent to the reference measure $\mu_{r}$ and the Radon--Nikodym derivative of $\nu$ with respect to $\mu_{r}$ is given by 
\begin{align}\label{defP2}
\frac{d\nu}{d\mu_{r}}(x) = \frac{1}{Z_{r}}\exp\bigg( -\sum_{j=1}^{M}\Phi_{j}^{r}(x_{j}) \bigg).
\end{align}
\end{assumption}

Following Assumption \ref{assumption1}, we know that the approximate measure can be decomposed as
$\nu(dx) = \prod_{j=1}^{M}\nu^{j}(dx_j)$ with
\begin{align}\label{aaa1}
\frac{d\nu^{j}}{d\mu_{r}^{j}} = \frac{1}{Z_{r}^{j}}\exp\big( -\Phi_{j}^{r}(x_{j}) \big),
\end{align}
where $Z_{r}^{j} = \mathbb{E}^{\mu_{r}^{j}}\big(\exp\big( -\Phi_{j}^{r}(x_{j}) \big)\big)$ ensures that $\nu^{j}$
is a probability measure. For $j=1, 2, \dots, M$, let $\mathcal{Z}_j$ be defined as a Hilbert space embedded in $\mathcal{H}_j$. Then, for $j=1,2,\cdots,M$, we introduce
\begin{align*}
& \qquad\qquad\,\,
\text{R}^1_j = \bigg\{
\Phi_j^r \, \Big| \, \sup_{1/N \leq \|x_j\|_{\mathcal{Z}_j} \leq N} \Phi_j^r(x_j) < \infty
\text{ for all }N>0
\bigg\},   \\
& \text{R}^2_j =
\bigg\{\Phi_j^r \, \Big| \, \int_{\mathcal{H}_j}\exp\left( -\Phi_{j}^{r}(x_j) \right)
\max(1, a_j(\epsilon,x_j)) \mu_r^{j}(dx_j) < \infty, \text{ for }\epsilon\in[0,\epsilon_0^j) \bigg\},
\end{align*}
where $\epsilon_0^j$ and $a_j(\cdot,\cdot)$ are defined as in Assumption \ref{assumption1}.
With these preparations, we can define $\mathcal{A}_{j}$ ($j=1,2,\cdots,M$) as follows:
\begin{align}\label{defAj}
\mathcal{A}_{j} = \bigg\{ \nu^{j}\in\mathcal{M}(\mathcal{H}_{j}) \, \bigg| \,
\begin{array}{l}
\nu^{j} \text{ is equivalent to } \mu_{r}^{j} \text{ with }
(\ref{aaa1}) \text{ holding true, }  \\
\text{and }
\Phi_{j}^{r}\in \text{R}^1_j \cap \text{R}^2_j
\end{array}
\bigg\}.
\end{align}

Now, we state the main theorem that yields practical iterative algorithms.

\begin{theorem}\label{theoremGeneral}
Assume that the approximate probability measure in problem (\ref{mipro1}) satisfies Assumption \ref{assumption1}.
For $j=1,2,\cdots,M$, we define $T_{N}^{j} = \{ x_j \, |\, 1/N\leq \|x_j\|_{\mathcal{Z}_{j}} \leq N \}$, where $N$ is a positive constant. Each reference measure $\mu_{r}^{j}$ is assumed to satisfy $\sup_{N} \mu_{r}^{j}(T_{N}^{j}) = 1$. In addition, we assume
\begin{align}\label{assMain}
\sup_{x_i\in T_N^i}
\int_{\prod_{j\neq i} \mathcal{H}_j} \left( \Phi^0(x) + \Phi(x) \right)  1_A(x)
\prod_{j\neq i}\nu^{j}(dx_j) < \infty
\end{align}
and
\begin{align}\label{assMain1}
\int_{\mathcal{H}_i}\!\exp\bigg(\!\! -\int_{\prod_{j\neq i}\mathcal{H}_j}\!\!\!(\Phi^0(x)+\Phi(x))1_{A^c}(x)\prod_{j\neq i}\nu^{j}(dx_j)\!\bigg)M_i(x)\mu_{r}^{i}(dx_i) < \infty,
\end{align}
where $A:=\{ x \,|\, \Phi^0(x)+\Phi(x) \geq 0 \}$,
and $M_i(x):=\max\left(1, a_i(\epsilon, x_i)\right), i,j=1,2,\cdots,M$.
Then, problem (\ref{mipro1}) admits a solution $\nu = \prod_{j=1}^{M}\nu^{j}\in\mathcal{M}(\mathcal{H})$ satisfying
\begin{align}
\frac{d\nu}{d\mu_{r}} \propto \exp\bigg( -\sum_{i=1}^{M}\Phi_{i}^{r}(x_{i}) \bigg),
\end{align}
where
\begin{align}\label{mainThSpec1}
\Phi_{i}^{r}(x_{i}) & = \int_{\prod_{j\neq i}\mathcal{H}_{j}}\!\! \bigg(\Phi^{0}(x) +
\Phi(x) \bigg) \prod_{j\neq i}\nu^{j}(dx_{j}) + \text{Const}
\end{align}
and
\begin{align}
\nu^{i}(dx_{i}) & \propto \exp\big( -\Phi_{i}^{r}(x_{i}) \big)\mu_{r}^{i}(dx_{i}).
\end{align}
\end{theorem}

It should be pointed out that this theorem and the meanings of ${\rm R}_j^1$ and ${\rm R}_j^2$ ($j=1,\ldots,M$) are slightly different from the statements given in \citep{Jia2020SIAM}. The proof of Theorem 11 in \citep{Jia2020SIAM} can be taken step by step to prove Theorem \ref{theoremGeneral}. In fact, the present version can be seen as an improved version of Theorem 11 in \citep{Jia2020SIAM}, which can be verified more easily for practical problems.

\subsection{Proof details}\label{AppendixProof}
\noindent\textbf{Proof of Theorem \ref{keyVItheorem}}

\begin{proof}
	To prove the theorem, it suffices to verify the conditions given in Theorem \ref{theoremGeneral}.
	Specifically speaking, we need to verify the following four inequalities:
	\begin{align*}
	T_1 = \sup_{u\in T_{N}^{u}} \int_{(\mathbb{R}^{+})^{N_d}}\Phi(u,\bm{\sigma})&\nu^{\bm{\sigma}}(d\bm{\sigma}) < \infty, \\
	T_2 = \sup_{\bm{\sigma}\in T_{N}^{\bm{\sigma}}} \int_{\mathcal{H}_u} \Phi(u,\bm{\sigma})\nu^{\bm{\sigma}}& \nu^{u}(du) < \infty,  \\
	T_3 = \int_{(\mathbb{R}^{+})^{N_d}}\exp\left( -\int_{\mathcal{H}_u}\Phi(u,\bm{\sigma}) \nu^{u}(du) \right)&\max(1, a(\epsilon,\bm{\sigma}))\mu_0^{\bm{\sigma}}(d\bm{\sigma}) < \infty, \\
	T_4 = \int_{\mathcal{H}_u} \exp\left(-\int_{(\mathbb{R}^{+})^{N_d}} \Phi(u,\bm{\sigma}) \nu^{\bm{\sigma}}(d\bm{\sigma}) \right)
	&\max(1, \|u\|_{\mathcal{H}_u}^2)\mu_0^u(du) < \infty.
	\end{align*}
	For term $T_1$, we have
	\begin{align*}
	T_1 \leq & \sup_{u\in T_{N}^{u}}\frac{1}{2}\int_{(\mathbb{R}^+)^{N_d}} \left(\sum_{k=1}^{N_d}\frac{1}{\sigma_k}\|\bm{d}-Hu\|_2^2 +
	\sum_{k=1}^{N_d} \log\sigma_k \right) \nu^{\bm{\sigma}}(d\bm{\sigma}) \\
	\leq & \left(1+\sup_{u\in T_{N}^u}\|\bm{d}-Hu\|_2^2 \right)\int_{(\mathbb{R}^+)^{N_d}}
	\sum_{k=1}^{N_d}\left( \frac{1}{\sigma_k} + \log\sigma_k \right)\nu^{\bm{\sigma}}(d\bm{\sigma}) \\
	\leq & C \int_{(\mathbb{R}^+)^{N_d}}
	\sum_{k=1}^{N_d}\left( \exp\left( \frac{\epsilon'}{\sigma_k} \right) + \sigma_k^{\epsilon'} \right)\nu^{\bm{\sigma}}(d\bm{\sigma})	\\
	\leq & C \int_{(\mathbb{E}^{+})^{N_d}}\exp\left(-\Phi_{\bm{\sigma}}(\bm{\sigma})\right)\max(1,a(2\epsilon',\bm{\sigma}))\mu_0^{\bm{\sigma}}(d\bm{\sigma}) < \infty,
	\end{align*}
	where $\epsilon'<\epsilon_0/2$ and $C$ is a generic constant that may be different from line to line.
	
For term $T_2$, we have
	\begin{align*}
	T_2 \leq & \sup_{\bm{\sigma}\in T_{N}^{\bm{\sigma}}} \frac{1}{2}\sum_{k=1}^{N_d}\left( \frac{1}{\sigma_k} + \log\sigma_k \right)
	\int_{\mathcal{H}_u}\left(1 + \|\bm{d}-Hu\|_2^2\right)\nu^u(du).
	\end{align*}
	Considering $\nu^u\in\mathcal{A}_u$, we get
	\begin{align*}
	\int_{\mathcal{H}_u}\left(1 + \|\bm{d}-Hu\|_2^2\right) \,\nu^u(du) \leq C\int_{\mathcal{H}_u}(1+\|u\|_{\mathcal{H}_u}^2)\nu^u(du) < \infty
	\end{align*}
	and
	\begin{align*}
	\sum_{k=1}^{N_d}\left( \frac{1}{\sigma_k} + \log\sigma_k \right) \leq & \sum_{k=1}^{N_d}\left( \frac{1}{\sigma_k} + \sigma_k \right)
	\leq 2NN_d < \infty,
	\end{align*}
	which indicates $T_2 < \infty$.
	
Using the estimate
	\begin{align*}
	\exp\left( -\int_{\mathcal{H}_u}\frac{1}{2}\left( \|d-Hu\|_{\Sigma}^2 + \sum_{k=1}^{N_d}\log\sigma_k \right) \nu^u(du) \right) \leq &
	\exp\left( -\frac{1}{2}\sum_{k=1}^{N_d}\log\sigma_k \right),
	\end{align*}
	we can estimate $T_3$ as follows:
	\begin{align*}
	T_3 \leq & \int_{(\mathbb{R}^+)^{N_d}}\prod_{k=1}^{N_d}\sigma_k^{-1/2}\max\left(1, a(\epsilon,\bm{\sigma})\right)\mu_0^{\bm{\sigma}}(d\bm{\sigma}) \\
	\leq & \sqrt{N_d}\int_{(\mathbb{R}^+)^{N_d}} \sum_{k=1}^{N_d}\sigma_{k}^{-1/2}\max(1, a(\epsilon,\bm{\sigma}))\mu_0^{\bm{\sigma}}(d\bm{\sigma})<\infty,
	\end{align*}
	where the last inequality follows from the properties of inverse Gamma distribution.
	Notice that
	\begin{align*}
	& \exp\left( -\frac{1}{2}\int_{(\mathbb{R}^+)^{N_d}}\|\bm{d}-Hu\|_2^2 + \sum_{k=1}^{N_d}\log\sigma_k \nu^{\bm{\sigma}}(d\bm{\sigma})  \right)  \\
	& \qquad\qquad\qquad\qquad\qquad\qquad
	\leq \exp\left( -\frac{1}{2}\sum_{k=1}^{N_d}\int_{(\mathbb{R}^+)^{N_d}}\log\sigma_k \nu^{\bm{\sigma}}(\bm{\sigma}) \right) < \infty.
	\end{align*}
	
Then, for term $T_4$, we have 
	\begin{align}
	T_4 \leq C \int_{\mathcal{H}_u}\max(1, \|u\|_{\mathcal{H}_u}^2)\mu_0^u(du) < \infty.
	\end{align}
The proof is completed by noting the general theory shown in Subsection \ref{AppendixGeneTheoryVI} in the Appendix. 
\end{proof}

\noindent\textbf{Proof of Theorem \ref{vanishNoiseTheo}}
\begin{proof}
	Let us firstly give the explicit form of the operator $H^*$. Taking $\bm{f}=(f_1,\ldots,f_{N_d})^T \in \mathbb{R}^{N_d}$
	and $g=\sum_{k=1}^{\infty}g_ke_k\in\mathcal{H}_u$, we have
	\begin{align*}
	\langle H^*\bm{f}, g\rangle_{L^2} = \langle \bm{f}, Hg\rangle_{\ell^2} = \sum_{k=1}^{\infty}\lambda_k g_k \langle \bm{f}, Se_k \rangle_{\ell^2}
	= \sum_{k=1}^{\infty}\lambda_k g_k \sum_{i=1}^{N_d} f_i e_k(x_i),
	\end{align*}
	which implies
	\begin{align}
	H^*\bm{f} = \sum_{k=1}^{\infty} \left( \sum_{i=1}^{N_d}f_i\lambda_k e_k(x_i) \right) e_k(x).
	\end{align}
	Then, we have
	\begin{align}
	\begin{split}
	H^*\Sigma_{\text{inv}}^* & H\bar{u}_0 = \sum_{k=1}^{\infty} \lambda_k \bar{u}_{0k} H^*\Sigma_{\text{inv}}^* Se_k \\
	= & \sum_{k=1}^{\infty}\lambda_k \bar{u}_{0k} H^* \left( e_k(x_1)\alpha_1^*/\beta_1^*, \ldots, e_k(x_{N_d})\alpha_{N_d}^*/\beta_{N_d}^* \right)^T  \\
	= & \sum_{j=1}^{\infty}\left( \sum_{i=1}^{N_d}\frac{\alpha_i^*}{\beta_i^*}\sum_{k=1}^{\infty}\lambda_k\bar{u}_{0k}e_k(x_i)\lambda_j e_j(x_i) \right)e_j(x).
	\end{split}
	\end{align}
	With these preparations, we find that
	\begin{align}
	\begin{split}
	& (H^*\Sigma_{\text{inv}}^*H + \mathcal{C}_0^{-1})^{-1}\mathcal{C}_0^{-1}\bar{u}_0 =  \\
	& \qquad\qquad\qquad \sum_{j=1}^{\infty}\frac{\alpha_j^{-1}}{\sum_{i=1}^{N_d}\frac{\alpha_i^*}{\beta_i^*}\sum_{k=1}^{\infty}\lambda_k\frac{\bar{u}_{0k}}{\bar{u}_{0j}}e_k(x_i)\lambda_j e_j(x_i) + \alpha_j^{-1}}\bar{u}_{0j}e_j(x),
	\end{split}
	\end{align}
	which obviously implies the desired conclusion.
\end{proof}

\noindent\textbf{Proof of Theorem \ref{errorSimpleAnalysis}}
\begin{proof}
	Based on Eq. (\ref{ForE}), we have
	\begin{align}
	\begin{split}
	\bar{u}_p = & \mathcal{C}_p \left[
	(H^*\Sigma_{\text{inv}}^*H + \mathcal{C}_0^{-1})u^{\dag} + H^*\Sigma_{\text{inv}}^*\bm{\epsilon}
	+ \mathcal{C}_0^{-1}(\bar{u}_0 - u^{\dag})
	\right] \\
	= & u^{\dag} + \mathcal{C}_p\left[
	H^*\Sigma_{\text{inv}}^*\bm{\epsilon} + \mathcal{C}_0^{-1}(\bar{u}_0 - u^{\dag})
	\right],
	\end{split}
	\end{align}
	which yields 
	\begin{align}\label{simpleEst1}
	\begin{split}
	\mathbb{E}_0\left[\|\bar{u}_p - u^{\dag}\|_{\mathcal{H}}\right] \leq &
	\mathbb{E}_0\left[ \|\mathcal{C}_pH^*\Sigma_{\text{inv}}^*\bm{\epsilon}\|_{\mathcal{H}_u} \right] +
	\|\mathcal{C}_p\mathcal{C}_0^{-1}(\bar{u}_0 - u^{\dag})\|_{\mathcal{H}_u}	 \\
	\leq & \mathbb{E}_0\left[ \|\mathcal{C}_pH^*\Sigma_{\text{inv}}^*\bm{\epsilon}\|_{\mathcal{H}_u}^2 \right]^{1/2}
	+ \|\bar{u}_0 - u^{\dag}\|_{\mathcal{H}_u}  \\
	= & \text{tr}\left( \mathcal{C}_pH^*\Sigma_{\text{inv}}^*\Sigma \Sigma_{\text{inv}}^*H \mathcal{C}_p^* \right)^{1/2}
	+ \|\bar{u}_0 - u^{\dag}\|_{\mathcal{H}_u}.
	\end{split}
	\end{align}
	Now, we have
	\begin{align}
	\begin{split}
	\mathbb{E}_0[\|\bar{u}_p(\bm{d};W_I) - u^{\dag}\|_{\mathcal{H}_u}] \leq &
	\mathbb{E}_0[\|\bar{u}_p(\bm{d};W_I) - \bar{u}_p\|_{\mathcal{H}_u}] +
	\mathbb{E}_0[\|\bar{u}_p - u^{\dag}\|_{\mathcal{H}_u}],
	\end{split}
	\end{align}
	which yields the required estimate by inserting the above inequality (\ref{simpleEst1}).
\end{proof}

\noindent\textbf{Proof of Theorem \ref{equivalentTheorem}}
\begin{proof}
	Because the operator $A$ is assumed to be a self-adjoint and positive definite, it holds  for appropriate $f$ that
	\begin{align}\label{po1}
	\langle (\underline{a} + \delta A^{\alpha/2})f, f\rangle \leq \langle (a+\delta A^{\alpha/2})f, f\rangle \leq
	\langle (\bar{a} + \delta A^{\alpha/2})f, f\rangle,
	\end{align}
	where we denote $a := a(\bm{d};W_I)$ and omit all of the identity operators. For concision, we omit the identity operators in the proof when there are no ambiguities from the context.
	Obviously, we also have
	\begin{align}\label{po2}
	\langle (\underline{a} + \delta A^{\alpha/2})^2f, f\rangle \leq \langle (a+\delta A^{\alpha/2})^2f, f\rangle \leq
	\langle (\bar{a} + \delta A^{\alpha/2})^2f, f\rangle,
	\end{align}
	when choosing appropriate $f$ in each estimate.
	According to Theorem 2.25 (Feldman-Hajek theorem) in \citep{Prato2014book}, we only need
	to prove the following two results:
	\begin{enumerate}
		\item $\mathcal{C}_{\epsilon_0}^{1/2}\mathcal{H} = \mathcal{C}_p^{1/2}\mathcal{H}=\mathcal{H}_0,$
		\item the operator
		$T:=(\mathcal{C}_{p}^{-1/2}\mathcal{C}_{\epsilon_0}^{1/2})(\mathcal{C}_{p}^{-1/2}\mathcal{C}_{\epsilon_0}^{1/2})^*-\text{Id}$
		is a Hilbert-Schmidt operator on $\bar{\mathcal{H}}_0$.
	\end{enumerate}
	First, we prove (1). Taking arbitrary $u\in\mathcal{C}_{\epsilon_0}^{1/2}\mathcal{H}$ gives
	\begin{align}
	u = \mathcal{C}_{\epsilon_0}^{1/2}v \quad \text{for some }v\in\mathcal{H}.
	\end{align}
	The variable can be rewritten as $u = \mathcal{C}_p^{1/2}\mathcal{C}_p^{-1/2}\mathcal{C}_{\epsilon_0}^{1/2}v$.
	Considering (\ref{po2}), we have
	\begin{align}\label{est11}
	\begin{split}
	\|\mathcal{C}_p^{-1/2}\mathcal{C}_{\epsilon_0}^{1/2}v\|_{\mathcal{H}}^2 = &
	\|(a+\delta A^{\alpha/2})(\epsilon_0^{-1}+\delta A^{\alpha/2})^{-1}v\|_{\mathcal{H}}^2 \\
	\leq & \|(\bar{a}+\delta A^{\alpha/2})(\epsilon_0^{-1}+\delta A^{\alpha/2})^{-1}v\|_{\mathcal{H}}^2 \\
	= & \sum_{k=1}^{\infty}\left(\frac{\bar{a}+\delta\lambda_k^{\alpha/2}}{\epsilon_0^{-1}+\delta\lambda_k^{\alpha/2}}\right)^2 v_k^2 \\
	\leq & C \|v\|_{\mathcal{H}}^2 < \infty,
	\end{split}
	\end{align}
	which implies $u\in\mathcal{C}_p^{1/2}\mathcal{H}$. Conversely, we assume $u\in\mathcal{C}_p^{1/2}\mathcal{H}$ with
	$u = \mathcal{C}_p^{1/2}v$ for some $v\in\mathcal{H}$. Similarly, we can rewrite
	$u=\mathcal{C}_{\epsilon_0}^{1/2}\mathcal{C}_{\epsilon_0}^{-1/2}\mathcal{C}_p^{1/2}v$. In the following, we give estimates of
	$\mathcal{C}_{\epsilon_0}^{-1/2}\mathcal{C}_p^{1/2}v$ which is more complex than the estimate (\ref{est11}).
	Obviously, we have
	\begin{align}\label{estPar1}
	\|\mathcal{C}_{\epsilon_0}^{-1/2}\mathcal{C}_p^{1/2}v\|_{\mathcal{H}}^2 = & \sum_{k,\ell =1}^{\infty}v_k v_{\ell}
	\langle \mathcal{C}_{\epsilon_0}^{-1/2}\mathcal{C}_p^{1/2}e_k, \mathcal{C}_{\epsilon_0}^{-1/2}\mathcal{C}_p^{1/2}e_{\ell} \rangle_{\mathcal{H}}
	= \sum_{k,\ell =1}^{\infty}v_k v_{\ell}B_{k,\ell}.
	\end{align}
	For the term $B_{k,\ell}$, we have
	\begin{align}\label{estPar2}
	\begin{split}
	B_{k,\ell} = & \langle (a-\epsilon_0^{-1})\mathcal{C}_p^{1/2}e_k, (a-\epsilon_0^{-1})\mathcal{C}_p^{1/2}e_{\ell} \rangle_{\mathcal{H}}
	+ 2\langle (\epsilon_0^{-1}-a)\mathcal{C}_{p}^{1/2}e_k, e_{\ell}\rangle + \langle e_k, e_{\ell} \rangle \\
	= & B_{k,\ell}^{1} + B_{k,\ell}^{2} + B_{k,\ell}^3.
	\end{split}
	\end{align}
	For the term $B_{k,\ell}^{1}$, we have
	\begin{align}
	\begin{split}
	B_{k,\ell}^{1} \leq & C \|(a+\delta A^{\alpha/2})^{-1}e_k\|_{\mathcal{H}} \|(a+\delta A^{\alpha/2})^{-1}e_{\ell}\|_{\mathcal{H}} \\
	\leq & C \|(\underline{a}+\delta A^{\alpha/2})^{-1}e_k\|_{\mathcal{H}} \|(\underline{a}+\delta A^{\alpha/2})^{-1}e_{\ell}\|_{\mathcal{H}} \\
	= & C \frac{1}{\underline{a} + \delta\lambda_{k}^{\alpha/2}} \frac{1}{\underline{a} + \delta\lambda_{\ell}^{\alpha/2}}.
	\end{split}
	\end{align}
	Hence, we further obtain
	\begin{align}\label{xuyao11}
	\begin{split}
	\sum_{k,\ell=1}^{\infty}B_{k,\ell}^{1} \leq	C \left[\sum_{k=1}^{\infty}\left(v_k^2 + \left(\frac{1}{\underline{a}+\delta\lambda_k^{\alpha/2}}\right)^2 \right)\right]^2 < \infty.
	\end{split}
	\end{align}
	For the term $B_{k,\ell}^{2}$, we have the following estimate
	\begin{align}\label{xuyao12}
	\begin{split}
	\sum_{k,\ell=1}^{\infty}v_kv_{\ell}B_{k,\ell}^{2} = & 2\sum_{k}^{\infty}v_k \langle (\epsilon_0^{-1}-a)\mathcal{C}_{p}^{1/2}e_k, \sum_{\ell=1}^{\infty}v_{\ell}e_{\ell}\rangle_{\mathcal{H}} \\
	\leq & C \sum_{k=1}^{\infty}\left(v_k^2 + \left(\frac{1}{\underline{a}+\delta\lambda_k^{\alpha/2}}\right)^2 \right) \|v\|_{\mathcal{H}} < \infty.
	\end{split}
	\end{align}
	For the term $B_{k,\ell}^{3}$, we find
	\begin{align}\label{xuyao13}
	\sum_{k,\ell=1}^{\infty}v_kv_{\ell}B_{k,\ell}^{2} \leq \|v\|_{\mathcal{H}}^2 < \infty.
	\end{align}
	Inserting estimates from (\ref{xuyao11}) to (\ref{xuyao13}) into (\ref{estPar1}), we arrive at 
	\begin{align}
	\|\mathcal{C}_{\epsilon_0}^{-1/2}\mathcal{C}_p^{1/2}v\|_{\mathcal{H}}^2 < \infty.
	\end{align}
	Hence, the proof of (1) is completed.
	Before going further, we denote
	$\langle\cdot,\cdot\rangle_{\mathcal{H}_0}:=\langle\cdot,\cdot\rangle_{\mathcal{C}_{\epsilon_0}^{1/2}\mathcal{H}} = \langle (\epsilon_0+\delta A^{\alpha/2})\cdot, (\epsilon_0+\delta A^{\alpha/2})\cdot \rangle_{\mathcal{H}}$.
	For proving (2), we introduce $\tilde{e}_j = \frac{1}{\epsilon_0^{-1}+\delta\lambda_j^{\alpha/2}}e_j$
	for $j=1,2,\ldots$, which is an orthonormal basis on $\mathcal{C}_{\epsilon_0}^{1/2}\mathcal{H}$.
	Following simple calculations, we obtain 
	\begin{align}
	\begin{split}
	\sum_{j=1}^{\infty}\langle T\tilde{e}_j, T\tilde{e}_j \rangle_{\mathcal{H}_0} = &
	\sum_{j=1}^{\infty}\langle (\mathcal{C}_p^{-1/2}\mathcal{C}_{\epsilon_0}\mathcal{C}_p^{-1/2}-\text{Id})\tilde{e}_j,
	(\mathcal{C}_p^{-1/2}\mathcal{C}_{\epsilon_0}\mathcal{C}_p^{-1/2}-\text{Id})\tilde{e}_j \rangle_{\mathcal{H}_0} \\
	= & \sum_{j=1}^{\infty} (I_j^1 + I_j^2 + I_j^3),
	\end{split}
	\end{align}
	where
	\begin{align*}
	& I_j^1 = \langle (\mathcal{C}_p^{-1/2}\mathcal{C}_{\epsilon_0}\mathcal{C}_p^{-1/2} - \mathcal{C}_{\epsilon_0}^{1/2}\mathcal{C}_p^{-1/2})\tilde{e}_j,  (\mathcal{C}_p^{-1/2}\mathcal{C}_{\epsilon_0}\mathcal{C}_p^{-1/2} - \mathcal{C}_{\epsilon_0}^{1/2}\mathcal{C}_p^{-1/2})\tilde{e}_j \rangle_{\mathcal{H}_0},
	\end{align*}
	\begin{align*}
	& I_j^2 = 2\langle (\mathcal{C}_p^{-1/2}\mathcal{C}_{\epsilon_0}\mathcal{C}_p^{-1/2} - \mathcal{C}_{\epsilon_0}^{1/2}\mathcal{C}_p^{-1/2})\tilde{e}_j,
	(\mathcal{C}_{\epsilon_0}^{1/2}\mathcal{C}_p^{-1/2} - \text{Id})\tilde{e}_j \rangle_{\mathcal{H}_0},
	\end{align*}
	\begin{align*}
	& I_j^3 = \langle (\mathcal{C}_{\epsilon_0}^{1/2}\mathcal{C}_p^{-1/2} - \text{Id})\tilde{e}_j,
	(\mathcal{C}_{\epsilon_0}^{1/2}\mathcal{C}_p^{-1/2} - \text{Id})\tilde{e}_j \rangle_{\mathcal{H}_0}.
	\end{align*}
	Since $I_j^2$ can be bounded by $I_j^1$ combined with $I_j^3$, we only focus on the estimates of $I_j^1$ and $I_j^3$ in the following.
	
	\textbf{Estimate of $I_j^1$}:
	Because
	\begin{align}
	\begin{split}
	\mathcal{C}_p^{-1/2}\mathcal{C}_{\epsilon_0}^{1/2} - \text{Id} = &
	(a+\delta A^{\alpha/2})(\epsilon_0^{-1} + \delta A^{\alpha/2})^{-1} - \text{Id} \\
	= & \left( a+\delta A^{\alpha/2} - \epsilon_0^{-1} - \delta A^{\alpha/2} \right) (\epsilon_0^{-1} + \delta A^{\alpha/2})^{-1} \\
	= & (a - \epsilon_0^{-1})(\epsilon_0^{-1} + \delta A^{\alpha/2})^{-1},
	\end{split}
	\end{align}
	we have
	\begin{align}\label{xuyao1}
	\begin{split}
	I_j^1 = & \langle (\mathcal{C}_p^{-1/2}\mathcal{C}_{\epsilon_0}^{1/2} - \text{Id})\mathcal{C}_{\epsilon_0}^{1/2}\mathcal{C}_p^{-1/2}\tilde{e}_j, (\mathcal{C}_p^{-1/2}\mathcal{C}_{\epsilon_0}^{1/2} - \text{Id})\mathcal{C}_{\epsilon_0}^{1/2}\mathcal{C}_p^{-1/2}\tilde{e}_j \rangle_{\mathcal{H}_0} \\
	= & \langle (a-\epsilon_0^{-1})B\tilde{e}_j, (a-\epsilon_0^{-1})B\tilde{e}_j \rangle_{\mathcal{H}_0} \\
	= & \langle \mathcal{C}_{\epsilon_0}^{-1/2}[(a-\epsilon_0^{-1})B\tilde{e}_j],
	\mathcal{C}_{\epsilon_0}^{-1/2}[(a-\epsilon_0^{-1})B\tilde{e}_j] \rangle_{\mathcal{H}} \\
	\leq & 2 I_j^{11} + 2 I_j^{12},
	\end{split}
	\end{align}
	where $B\tilde{e}_j := \mathcal{C}_{\epsilon_0}\mathcal{C}_p^{-1/2}\tilde{e}_j$ and
	\begin{align*}
	& I_j^{11} = \langle (\epsilon_0^{-1}+\delta A^{\alpha/2})[aB\tilde{e}_j], (\epsilon_0^{-1}+\delta A^{\alpha/2})[aB\tilde{e}_j] \rangle_{\mathcal{H}}, \\
	& I_j^{12} = \epsilon_0^{-2}\langle (\epsilon_0^{-1}+\delta A^{\alpha/2})B\tilde{e}_j, (\epsilon_0^{-1}+\delta A^{\alpha/2})B\tilde{e}_j \rangle_{\mathcal{H}}.
	\end{align*}
	For the term $I_j^{11}$, we have
	\begin{align}
	\begin{split}
	I_j^{11} = & \langle (\epsilon_0^{-1}+\delta A^{\alpha/2})\sum_{k=1}^{\infty}\langle aB\tilde{e}_j, e_k \rangle_{\mathcal{H}} e_k,
	(\epsilon_0^{-1}+\delta A^{\alpha/2})\sum_{\ell=1}^{\infty}\langle aB\tilde{e}_j, e_{\ell} \rangle_{\mathcal{H}} e_{\ell}
	\rangle_{\mathcal{H}}   \\
	= & \sum_{k=1}^{\infty} \langle aB\tilde{e}_j, e_k\rangle_{\mathcal{H}}^2 \left( \frac{1}{\epsilon_0^{-1} + \delta\lambda_k^{\alpha/2}} \right)^2 \\
	= &  \left(\frac{1}{\epsilon_0^{-1} + \delta\lambda_j^{\alpha/2}}\right)^2
	\sum_{k=1}^{\infty} \left( \frac{1}{\epsilon_0^{-1} + \delta\lambda_k^{\alpha/2}} \right)^2
	\langle e_j, \mathcal{C}_p^{-1/2}\mathcal{C}_{\epsilon_0}(ae_k)\rangle_{\mathcal{H}}^2.
	\end{split}
	\end{align}
	Because
	\begin{align}
	\begin{split}
	& \sum_{j=1}^{\infty}\left(\frac{1}{\epsilon_0^{-1} + \delta\lambda_j^{\alpha/2}}\right)^2
	\langle e_j, (a+\delta A^{\alpha/2})(\epsilon_0^{-1} + \delta A^{\alpha/2})^{-2}(ae_k) \rangle_{\mathcal{H}}^2 \\
	& \qquad
	\leq C \|ae_k\|_{\mathcal{H}}^2\sum_{j=1}^{\infty}\left(\frac{1}{\epsilon_0^{-1} + \delta\lambda_j^{\alpha/2}}\right)^2 < \infty,
	\end{split}
	\end{align}
	we obtain
	\begin{align}\label{xuyao2}
	\sum_{j=1}^{\infty}I_j^{11} \leq C \|ae_k\|_{\mathcal{H}}^2 \sum_{j=1}^{\infty}\left(\frac{1}{\epsilon_0^{-1} + \delta\lambda_j^{\alpha/2}}\right)^2
	\sum_{k=1}^{\infty}\left(\frac{1}{\epsilon_0^{-1} + \delta\lambda_k^{\alpha/2}}\right)^2 < \infty.
	\end{align}
	For the term $I_j^{12}$, we have
	\begin{align}\label{zongIn1}
	\begin{split}
	I_j^{12} = & \langle (\epsilon_0^{-1} + \delta A^{\alpha/2})^{-1}(a+\delta A^{\alpha/2})\tilde{e}_j,
	(\epsilon_0^{-1} + \delta A^{\alpha/2})^{-1}(a+\delta A^{\alpha/2})\tilde{e}_j \rangle_{\mathcal{H}} \\
	\leq & 2
	\langle (\epsilon_0^{-1} + \delta A^{\alpha/2})^{-1}(a\tilde{e}_j), (\epsilon_0^{-1} + \delta A^{\alpha/2})^{-1}(a\tilde{e}_j) \rangle_{\mathcal{H}} \\
	& \qquad\qquad
	+ 2\delta^2\langle
	(\epsilon_0^{-1} + \delta A^{\alpha/2})^{-1}A^{\alpha/2}\tilde{e}_j, (\epsilon_0^{-1} + \delta A^{\alpha/2})^{-1}A^{\alpha/2}\tilde{e}_j
	\rangle_{\mathcal{H}}.
	\end{split}
	\end{align}
	For the first term on the right-hand side of the above inequality, we have
	\begin{align}\label{zongIn2}
	\begin{split}
	& \langle (\epsilon_0^{-1} + \delta A^{\alpha/2})^{-1}(a\tilde{e}_j), (\epsilon_0^{-1} + \delta A^{\alpha/2})^{-1}(a\tilde{e}_j) \rangle_{\mathcal{H}} \\
	= & \langle (\epsilon_0^{-1} + \delta A^{\alpha/2})^{-1}\sum_{k=1}^{\infty}\langle a\tilde{e}_j, e_k\rangle e_k,
	(\epsilon_0^{-1} + \delta A^{\alpha/2})^{-1}\sum_{\ell=1}^{\infty}\langle a\tilde{e}_j, e_\ell\rangle e_\ell \rangle_{\mathcal{H}} \\
	= & \sum_{k=1}^{\infty} \langle a\tilde{e}_j, e_k\rangle_{\mathcal{H}}^2 \left( \frac{1}{\epsilon_0^{-1} + \delta\lambda_k^{\alpha/2}} \right)^2 \\
	\leq & \bar{a}^2 \left( \frac{1}{\epsilon_0^{-1} + \delta\lambda_j^{\alpha/2}} \right)^2
	\sum_{k=1}^{\infty}\left( \frac{1}{\epsilon_0^{-1} + \delta\lambda_k^{\alpha/2}} \right)^2.
	\end{split}
	\end{align}
	For the second term on the right-hand side of (\ref{zongIn1}), we have
	\begin{align}\label{zongIn3}
	\begin{split}
	\langle
	(\epsilon_0^{-1} + \delta A^{\alpha/2})^{-1}A^{\alpha/2}\tilde{e}_j, (\epsilon_0^{-1} + \delta A^{\alpha/2})^{-1}A^{\alpha/2}\tilde{e}_j
	\rangle_{\mathcal{H}} \leq C \left( \frac{1}{\epsilon_0^{-1} + \delta\lambda_j^{\alpha/2}} \right)^2.
	\end{split}
	\end{align}
	Substituting estimates (\ref{zongIn2}) and (\ref{zongIn3}) into (\ref{zongIn1}), we obtain
	\begin{align}\label{xuyao3}
	\sum_{j=1}^{\infty}I_j^{12} \leq C \sum_{j=1}^{\infty}\left( \frac{1}{\epsilon_0^{-1} + \delta\lambda_j^{\alpha/2}} \right)^2 < \infty.
	\end{align}
	Combining estimates (\ref{xuyao1}), (\ref{xuyao2}), and (\ref{xuyao3}) leads to 
	\begin{align}
	\sum_{j=1}^{\infty}I_j^1 < \infty.
	\end{align}
	
	\textbf{Estimate of $I_j^3$}: For the term $I_j^3$, we find that
	\begin{align}
	\begin{split}
	I_j^3 = & \langle (\epsilon_0^{-1}+\delta A^{\alpha/2})^{-1}[(a-\epsilon_0^{-1})\tilde{e}_j],
	(\epsilon_0^{-1}+\delta A^{\alpha/2})^{-1}[(a-\epsilon_0^{-1})\tilde{e}_j] \rangle_{\mathcal{H}_0} \\
	= & \langle (a-\epsilon_0^{-1})\tilde{e}_j, (a-\epsilon_0^{-1})\tilde{e}_j \rangle_{\mathcal{H}} \\
	\leq & C \left( \frac{1}{\epsilon_0^{-1} + \delta \lambda_j^{\alpha/2}} \right)^2,
	\end{split}
	\end{align}
	which implies
	\begin{align}
	\sum_{j=1}^{\infty}I_j^3 < \infty.
	\end{align}
	Combining estimates of $I_j^1$ and $I_j^3$, we obtain
	\begin{align}
	\sum_{j=1}^{\infty}\langle T\tilde{e}_j, T\tilde{e}_j \rangle_{\mathcal{H}_0} < \infty,
	\end{align}
	which indicates that $T$ is a Hilbert--Schmidt operator. Hence, the proof is completed.
\end{proof}

\subsection{More numerical results of simple smoothing model}\label{SubSecAppendixNumSimpleSmooth}

\begin{figure}[thbp]
	\centering
	\includegraphics[width=0.8\textwidth]{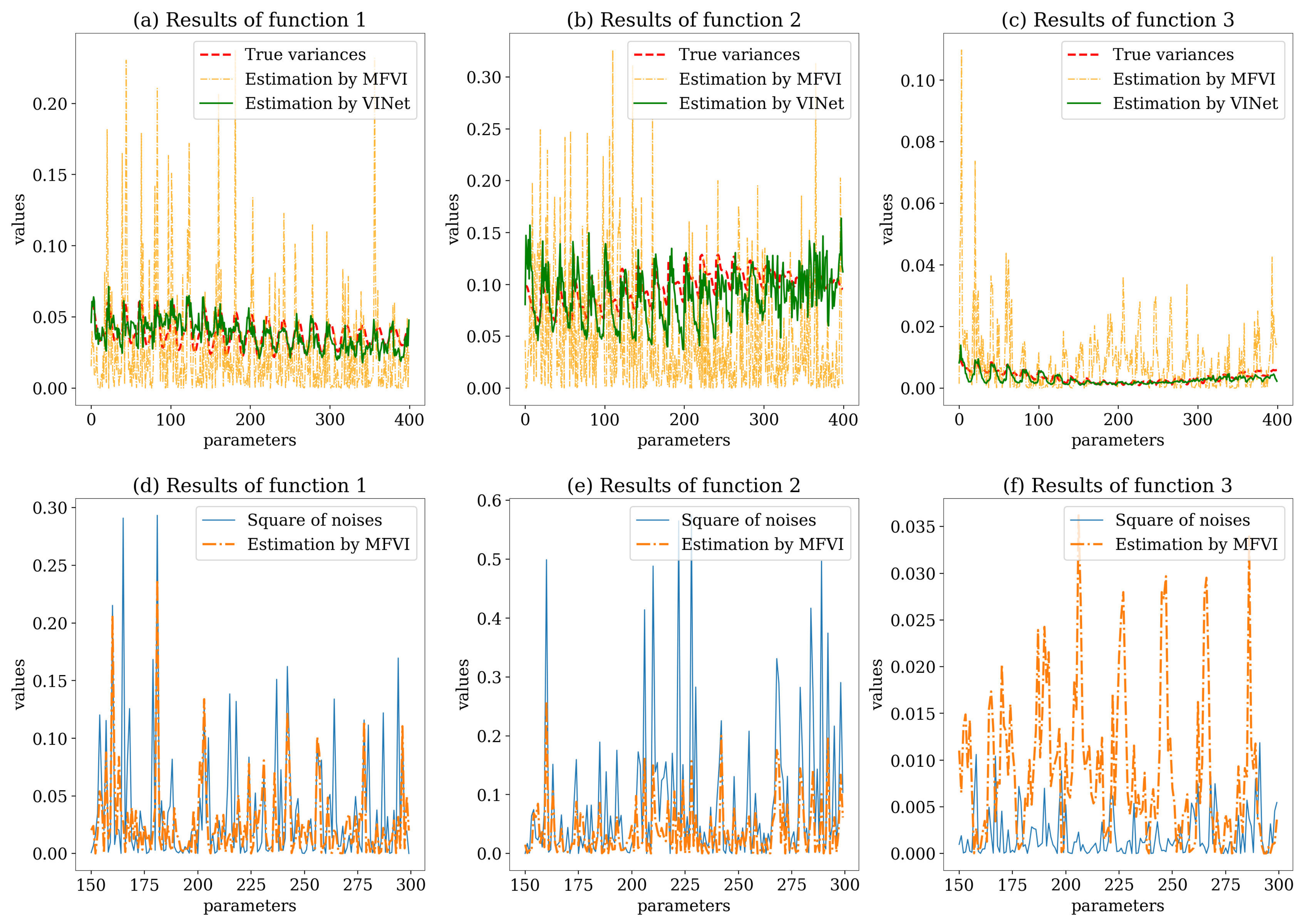}\\
	\vspace*{-0.6em}
	\caption{Comparison of the true noise variances and estimated noise variances by the MFVI and VINet (output of SNet).
		(a) (b) (c) Noise variances of all parameters for three different testing examples.
		(d) (e) (f) Noise variances estimated by MFVI and the square real noises of 
		some parameters. }\label{NoiseCompareType1TestData}
\end{figure}

\begin{figure}[thbp]
	\centering
	\includegraphics[width=0.8\textwidth]{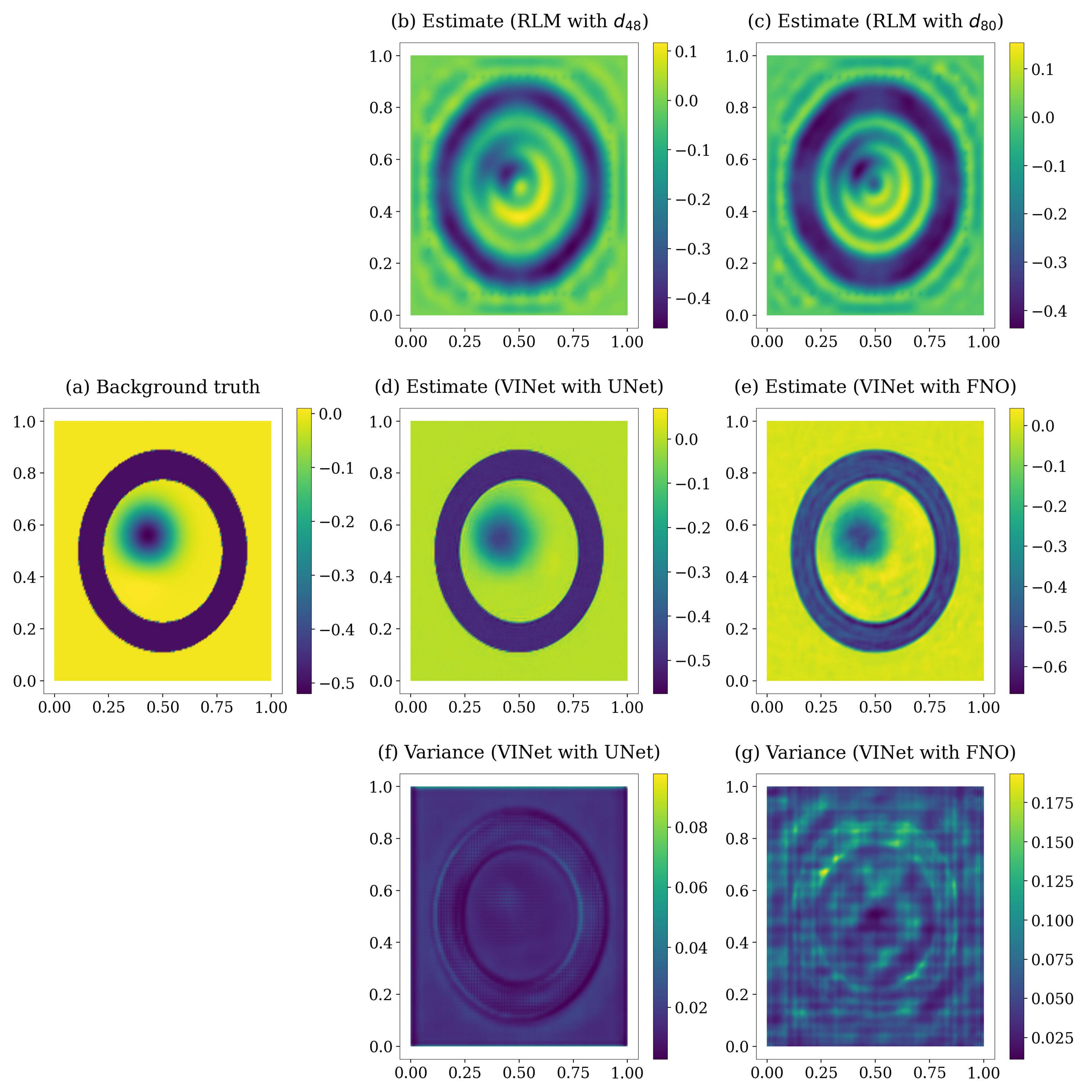}\\
	\vspace*{-0.6em}
	\caption{(a) Background true function randomly selected from type 1 testing datasets; 
		(b) Estimate obtained by the TSVD method; (c) Estimate obtained by the MFVI method; 
		(d) Estimate obtained by the VINet with U-ENet; (e) Estimate obtained by the VINet with FNO-ENet; 
		(f) Variance function obtained by the VINet with U-ENet; 
		(g) Variance function obtained by the VINet with FNO-ENet.}\label{InvResultsType1Test}
\end{figure}

In this subsection, let us provide the results for TSVD, MFVI, and VINet obtained by using the type 1 test 
datasets generated from the same probability measure as the training datasets. 

In Figure \ref{NoiseCompareType1TestData}, we depict the estimated noise variances by the MFVI and VINet and the background true noise variances for three randomly selected testing examples. In the first line of Figure \ref{NoiseCompareType1TestData}, the estimated variances of all $400$ measurement points are given with the dashed red line representing the true noise variances, the dashed-dotted light orange line representing the estimation given by MFVI, and the solid green line representing the estimation given by the SNet. As in the main test, the SNet provides accurate estimates of the noise variances for all of the three type 1 test datasets. In the second line of Figure \ref{NoiseCompareType1TestData}, we use the solid blue line and dashed-dotted orange line to represent the square real noises, i.e., $(\bm{d} - \bm{d}_c)^2$, and the estimations given by MFVI, respectively. Obviously, the MFVI method only captures some features of the square real noises but not the true noise variances.

In Figure \ref{InvResultsType1Test}, we show the background truth and estimates obtained by the TSVD, MFVI, VINet with U-ENet, and VINet with FNO-ENet for one randomly selected data from the type 1 testing datasets. Similar conclusions can be drawn as in the main text: the TSVD and MFVI methods can only capture some major trends of the truth, while the two types of VINets provide more detailed structures.

\begin{table}[thbp]
	\caption{Average relative errors of the estimates obtained by the TSVD, MFVI, VINet (U-ENet), and VINet (FNO-ENet)
		with mesh sizes $n=\{30\times30, 50\times50, 100\times100, 150\times150, 200\times200, 250\times250, 300\times300 \}$
		for the type 2 test datasets.}
	\begin{center}
		\begin{spacing}{1.5}
			\begin{tabular}{c|ccccc}
				\Xhline{1.1pt}
				& TSVD & MFVI & VINet (U-ENet) & VINet (FNO-ENet) \\
				\hline
				Relative error ($n=30^{2}$) & 0.1406 & 0.2301 & 0.0878 & 0.0786  \\
				\hline
				Relative error ($n=50^{2}$) & 0.1413 & 0.2328 & 0.0933 &  0.0822	\\
				\hline
				Relative error ($n=100^{2}$) & 0.1457 & 0.2499 & 0.0972 &  0.0854	\\
				\hline 
				Relative error ($n=150^{2}$) & 0.1481 & 0.2504 & 0.0982 &  0.0865	\\
				\hline
				Relative error ($n=200^{2}$) & 0.1494 & 0.2371 & 0.0989 &  0.0870	\\
				\hline 
				Relative error ($n=250^{2}$) & 0.1498 & 0.2499 & 0.0991 &  0.0871	\\
				\hline 
				Relative error ($n=300^{2}$) & 0.1499 & 0.2406 & 0.0991 & 0.0871	\\
				\Xhline{1.1pt}
			\end{tabular}
		\end{spacing}
	\end{center}\label{ResolutionCompareType1Ex1}
\end{table}

At last, in Table \ref{ResolutionCompareType1Ex1}, we provide a detailed comparison of TSVD, MFVI, VINet with U-ENet, and VINet with FNO-ENet with various discrete meshes. Similar to the results in the main text, the relative errors of the two VINets are lower than those of the TSVD and MFVI methods. All of the methods are stable with respect to different discrete meshes.

\subsection{More numerical results of inverse source problem}\label{SubSecAppendixNumInvSource}

\begin{figure}[thbp]
	\centering
	\includegraphics[width=0.8\textwidth]{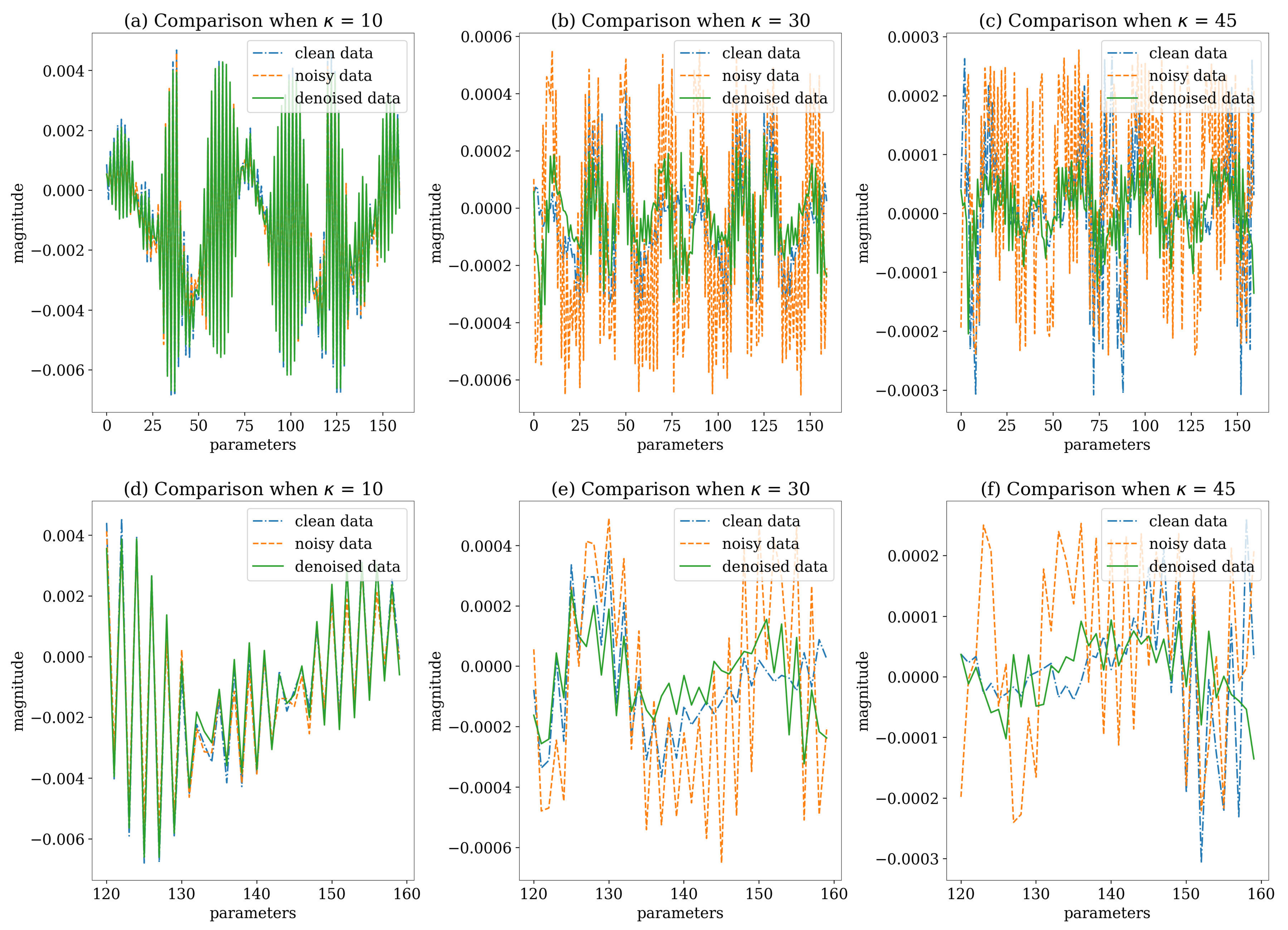}\\
	\vspace*{-0.6em}
	\caption{Noisy data (data generated by accurate forward solver with random Gaussian noise), clean data (data generated by rough forward solver), and denoised data (the output of the DNet) for a randomly selected source function $u$ from type 1 test datasets.
		(a) Data when the wavenumber $\kappa=10$;
		(b) Data when the wavenumber $\kappa=30$;
		(c) Data when the wavenumber $\kappa=45$;
		(d)(e)(f) Partly enlarged versions of (a), (b), and (c), respectively.}\label{d_invert_test1}
\end{figure}

\begin{figure}[thbp]
	\centering
	\includegraphics[width=0.8\textwidth]{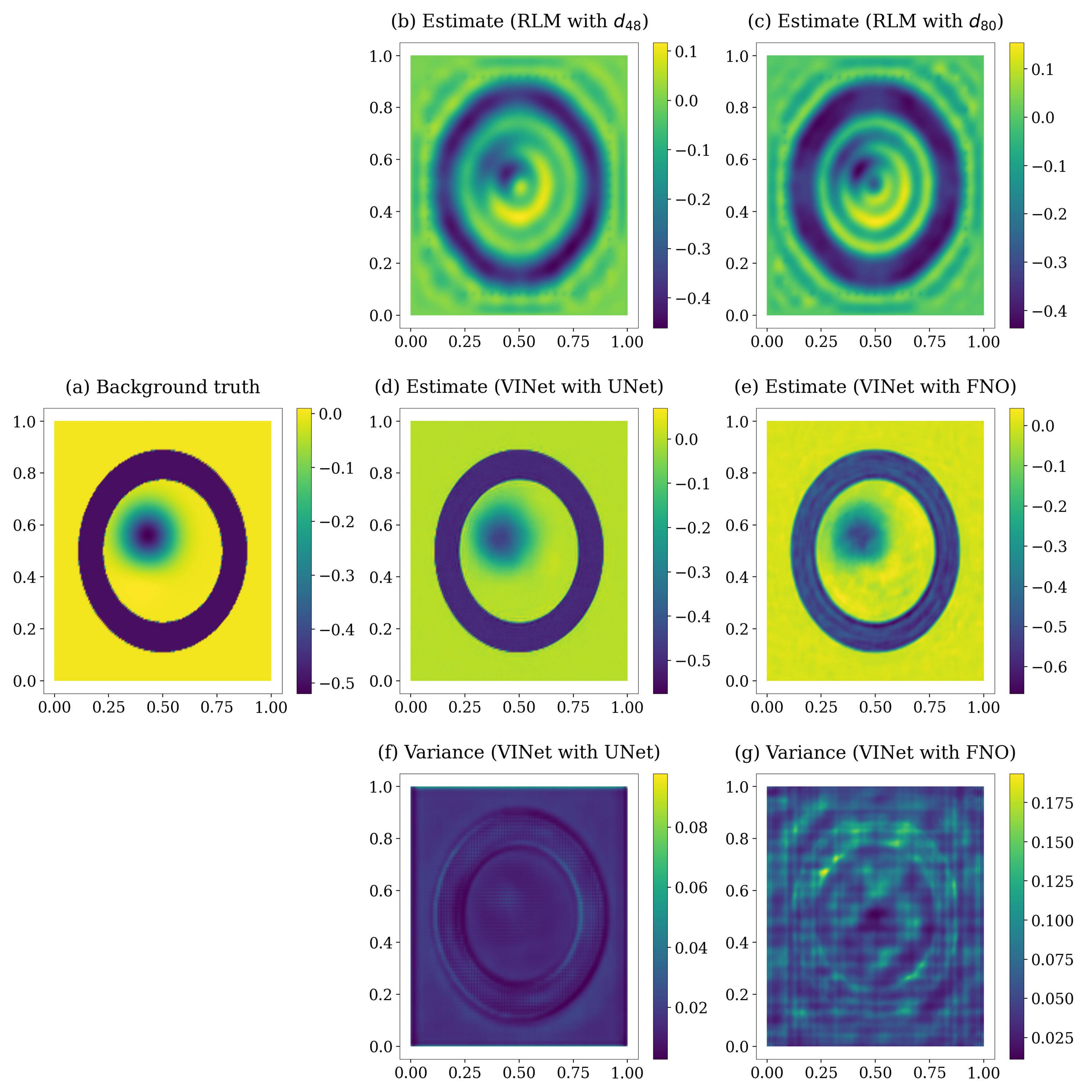}
	\vspace*{-0.6em}
	\caption{(a) One of the background true function in type 1 test datasets;
		(b) Estimate obtained by RLM with $\bm{d}_{48}$; (c) Estimate obtained by RLM with $\bm{d}_{80}$; 
		(d) Estimate obtained by VINet with U-ENet; (e) Estimate obtained by VINet with FNO-ENet; 
		(f) Estimated variance function by VINet with U-ENet; 
		(g) Estimated variance function by VINet with FNO-ENet.}\label{estUtype1}
\end{figure}

As in the main text, we first show the results of DNet. In Figure \ref{d_invert_test1}, we use the dash-dotted blue line, dashed orange line, and solid green line to represent the clean data, noisy data, and the output of the DNet, respectively. Sub-figures (a), (b), and (c) provide comparisons when the wavenumbers are chosen to be 10, 30, and 45. To provide more detail, we exhibit sub-figures (d), (e), and (f) that demonstrate the partly enlarged versions of sub-figures (a), (b), and (c), respectively. Similar to the main text, the noisy data deviate largely from the clean data and the outputs of DNet are visually similar to the clean data when the wavenumber $\kappa = 30, 45$. Hence, the DNet learns the sophisticated model error mechanism between the accurate and rough forward solvers.

In Figure \ref{estUtype1}, we exhibit the results of one randomly selected test example from the type 1 test datasets. Similar results can be obtained as in the main text. The classical RLM can hardly handle the discontinuity, even when the highest wavenumber is set to be $80$. VINets with U-ENet and FNO-ENet both yield more accurate estimations.

\begin{table}[t]
	\caption{Average relative errors of the estimates obtained by the RLM (using data $\bm{d}_{48}$ and $\bm{d}_{80}$), 
		VINet with U-ENet, and VINet with FNO-ENet with mesh sizes 
		$n = \{30\times30, 100\times100, 150\times150, 200\times200, 250\times250, 300\times300\}$ for the type 1 test datasets.}
	\begin{center}
		\begin{spacing}{1.5}
			\begin{tabular}{c|ccccc}
				\Xhline{1.1pt}
				& RLM($\bm{d}_{48}$) & RLM($\bm{d}_{80}$) & VINet(U-ENet) & VINet(FNO-ENet) \\
				\hline
				Relative error ($n=30^{2}$) & 0.3349 & 0.3250 & - & -  \\
				\hline
				Relative error ($n=100^{2}$) & 0.2310 & 0.2321 & 0.0474 &  0.1163	\\
				\hline 
				Relative error ($n=150^{2}$) & 0.2114 & 0.1932 & 0.0426 &  0.0443	\\
				\hline
				Relative error ($n=200^{2}$) & 0.2087 & 0.1728 & 0.0406 &  0.0449	\\
				\hline 
				Relative error ($n=250^{2}$) & 0.2084 & 0.1653 & 0.0405 &  0.0462	\\
				\hline 
				Relative error ($n=300^{2}$) & 0.2100 & 0.1624 & 0.0452 & 0.0489	\\
				\Xhline{1.1pt}
			\end{tabular}
		\end{spacing}
	\end{center}\label{resultsHelCompareType1}
\end{table}

Finally, we give a detailed comparison for RLM and VINet with various discrete meshes. In Table \ref{resultsHelCompareType1}, we show the average relative errors on different meshes obtained by RLM with $\bm{d}_{48}$, RLM with $\bm{d}_{80}$, VINet with U-ENet, and VINet with FNO-ENet, respectively. We can see that the discrete level indeed influences the performance of the classical RLM as illustrated in the main text. By comparing Tables 5 and \ref{resultsHelCompareType1}, we see that the relative errors of VINet with FNO-ENet are both around $4\%$, but the relative errors of VINet with U-ENet are around $4\%$ and $6\%$, respectively. Hence, under the current setting, the VINet with U-ENet seems to have weaker generalization properties compared with the VINet with FNO-ENet.


\section*{Acknowledgments}
The first author was supported by the NSFC (Grant Nos. 11871392, 12090020, and 12090021).
The third author was supported in part by the NSF grant DMS-2208256.
The fourth author was supported by the NSFC (Grant Nos. 11690011 and U1811461) and the Macao Science and Technology Development Fund grant No. 061/2020/A2.


\bibliographystyle{plain}
\bibliography{references}

\end{document}